\documentclass{amsart}
\newcommand{\mysec}{\section}
\newcommand{\UCSection}{Section}
\newcommand{\LCsection}{section}
\newcommand{\mychap}{\UCSection}
\newcommand{\referee}{Referee~2}

\usepackage{amssymb,amsmath,verbatim,amstext,eucal}
\usepackage[all]{xy}
\usepackage{amsrefs}
\usepackage{mathtools}
\usepackage{xcolor} 
\usepackage{fullpage}

\numberwithin{subsection}{section}
\numberwithin{equation}{section}

\newcommand{\numberby}{equation}

\newtheorem{lemma}[\numberby]{Lemma}

\newtheorem{proposition}[\numberby]{Proposition}
\newtheorem{theorem}[\numberby]{Theorem}
\newtheorem{conjecture}[\numberby]{Conjecture}
\newtheorem{corollary}[\numberby]{Corollary}

\newcommand{\FLEX}{\relax}
\newcommand{\flex}[1]{\renewcommand{\FLEX}{#1}}
\newtheorem{flexthm}[\numberby]{\FLEX}
\newenvironment{flexstate}[2]{\flex{#1}\begin{flexthm}[#2]}{\end{flexthm}}

\theoremstyle{definition}
\newtheorem{definition}[\numberby]{Definition}
\newtheorem{example}[\numberby]{Example}

\newenvironment{remark}[1]{\refstepcounter{\numberby}%
\vskip 6pt \par {\sc #1\ \thelemma .}}{\vskip 6pt \par}

\newenvironment{remark*}[1]{\par \vskip 5pt \noindent 
{\bf #1.}}{\vskip 5pt \par}

\newcommand{\bh}{\ensuremath{{\mathcal B}({\mathcal H})}}

\newcommand{\cstar}{\hbox{$C^*$}}
\newcommand{\cstaralg}{$C^*$-algebra}

\providecommand{\dual}[1]{\ensuremath{#1^{\#}}}

\newcommand{\dom}{\operatorname{dom}}
\newcommand{\dperp}{{\perp\perp}}
\newcommand{\ds}{\displaystyle}
\newcommand{\dstext}[1]{\quad\text{#1}\quad}
\newcommand{\eps}{\ensuremath{\varepsilon}}

\newcommand{\Id}{\operatorname{Id}}
\newcommand{\id}{{\operatorname{id}}}

\newcommand{\indlim}{\varinjlim}

\newcommand{\idealin}{\unlhd}

\newcommand{\norm}[1]{\left\lVert#1\right\rVert}

\providecommand{\proof}{\noindent {\em Proof.}\,\,}
\providecommand{\qed}%
{\hfill \vrule height5pt width4pt depth1pt \vspace{+2.00ex}}

\newcommand{\ran}{\operatorname{range}}

\newcommand{\sot}{\textsc{sot}}

\newcommand{\spn}{\operatorname{span}}
\newcommand{\supp}{\operatorname{supp}}
\newcommand{\tnorm}[1]{{\left\vert\kern-0.25ex\left\vert\kern-0.25ex\left\vert #1 
   \right\vert\kern-0.25ex\right\vert\kern-0.25ex\right\vert}}

\newcommand{\unit}[1]{#1^{(0)}}

\newcommand{\bbC}{{\mathbb{C}}}

\newcommand{\bbN}{{\mathbb{N}}}

\newcommand{\bbP}{{\mathbb{P}}}

\newcommand{\bbR}{{\mathbb{R}}}

\newcommand{\bbT}{{\mathbb{T}}}

  \newcommand{\A}{{\mathcal{A}}}
  
  \newcommand{\B}{{\mathcal{B}}}
  \newcommand{\C}{{\mathcal{C}}}
  
  \newcommand{\D}{{\mathcal{D}}}
  
  \newcommand{\E}{{\mathcal{E}}}
  \newcommand{\F}{{\mathcal{F}}}
  
\renewcommand{\H}{{\mathcal{H}}}
  \newcommand{\I}{{\mathcal{I}}}
  \newcommand{\J}{{\mathcal{J}}}
\renewcommand{\L}{{\mathcal{L}}}
  \newcommand{\M}{{\mathcal{M}}}
  \newcommand{\N}{{\mathcal{N}}}
\renewcommand{\O}{{\mathcal{O}}}
\renewcommand{\P}{{\mathcal{P}}}
  \newcommand{\Q}{{\mathcal{Q}}}
  
\renewcommand{\S}{{\mathcal{S}}}

  \newcommand{\Y}{{\mathcal{Y}}}
  \newcommand{\Z}{{\mathcal{Z}}}
\newcommand{\fA}{{\mathfrak{A}}}

\newcommand{\fC}{{\mathfrak{C}}}

\newcommand{\fI}{{\mathfrak{I}}}
\newcommand{\fJ}{{\mathfrak{J}}}

\newcommand{\fM}{{\mathfrak{M}}}

\newcommand{\fP}{{\mathfrak{P}}}
\newcommand{\fQ}{{\mathfrak{Q}}}

\newcommand{\fr}{{\mathfrak{r}}}
\newcommand{\fS}{{\mathfrak{S}}}
\newcommand{\fs}{{\mathfrak{s}}}

\newcommand{\fU}{{\mathfrak{U}}}

\providecommand{\ideal}{\operatorname{ideal}}
\providecommand{\ropen}{\operatorname{\textsc{Ropen}}}
\providecommand{\rideal}{\operatorname{\textsc{Rideal}}}

\providecommand{\iip}{ideal intersection property}

\DeclareMathOperator{\PsExp}{\operatorname{PsExp}}

\providecommand{\Eigone}{\E^1}
\newcommand{\fg}{\mathfrak{g}}
\newcommand{\pd}{pseudo-Cartan inclusion}
\newcommand{\Pd}{Pseudo-Cartan Inclusion}
\newcommand{\mac}{\breve}
\newcommand{\mintp}{\otimes_{\min}}
\newcommand{\ann}{\operatorname{Ann}}
\newcommand{\Ann}{\operatorname{Ann}}
\newcommand{\wnd}{weakly non-degenerate}
\newcommand{\rcom}{\operatorname{\textsc{rCom}}}
\newcommand{\udc}[1]{(#1^c)^{\sim}}

\newcommand{\ms}{\hbox{\,:\hskip -1pt :\,}}

\newcommand{\pad}{\vskip 4pt}
\long\def\dproof #1\enddproof{\pad\noindent \it Proof. \rm #1 \hfill$\diamondsuit$ \pad}

\newcounter{ccnt}
\long\def\dclaim#1%
\enddclaim{\stepcounter{ccnt}\pad \noindent\textsc{Claim \theccnt.} \it
  #1 \rm}

\long\def\dclaimnn#1%
\enddclaimnn{\pad \noindent\textsc{Claim.} \it
    #1 \rm}

 \newcounter{cct}
  \long\def\dstate#1#2%
  {\stepcounter{cct}\pad
    {\sc #1 \thecct.} \it #2\rm\pad}

   \newcommand{\oar}{\underrightarrow}

   \newcommand{\fix}[1]{\operatorname{fix} #1}
\newcommand{\Mod}{\text{Mod}}

\makeindex

\begin{document}

\title[Pseudo-Cartan Inclusions]{Pseudo-Cartan Inclusions}


\author[D.R. Pitts]{David R. Pitts}
\address[D.R. Pitts]{
  Department of Mathematics\\
  University of Nebraska-Lincoln\\
  Lincoln\\
  NE 68588-0130 U.S.A.}
\email{dpitts2@unl.edu}

\date{September 10, 2026}




\begin{abstract} A \pd\ is a regular inclusion having a Cartan
  envelope.  Unital \pd s were classified in~\cite{PittsStReInII}; we
  extend this classification to include the non-unital case.  The
  class of \pd s coincides with the class of regular inclusions having
  the faithful unique pseudo-expectation property and can also be
  described using the \iip.  We describe the twisted groupoid
  associated with the Cartan envelope of a \pd.  These results
  significantly extend previous results obtained for the unital
  setting.

  We explore properties of \pd s and the relationship between a \pd\
  and its Cartan envelope. For example, if $\D\subseteq \C$ is a \pd\
  with Cartan envelope $\B\subseteq \A$, then $\C$ is simple if and
  only if $\A$ is simple.  Also every regular $*$-automorphism of $\C$
  having $\D$ as an invariant subspace 
  uniquely extends to a $*$-automorphism of $\A$.  We show that the
  inductive limit of \pd s with suitable connecting maps is a \pd, and
  the minimal tensor product of \pd s is a \pd.  Further, we describe
  the Cartan envelope of \pd s arising from these constructions.  We
  conclude with some applications and a few open questions.
\end{abstract}

\newcommand{\tpad}{\hskip20pt}

\maketitle

\bigskip
\centerline{\sc Contents}

\contentsline {section}{\tocsection {}{1}{Introduction}}{2}{}%
\contentsline {section}{\tocsection {}{2}{Foundations}}{4}{}%
\contentsline {subsection}{\tocsubsection {}{\tpad2.1}{Terminology and Notation}}{4}{}%
\contentsline {subsection}{\tocsubsection {}{\tpad2.2}{Weakly Non-Degenerate Inclusions}}{7}{}%
\contentsline {subsection}{\tocsubsection {}{\tpad2.3}{Pseudo-Expectations}}{11}{}%
\contentsline {subsection}{\tocsubsection {}{\tpad2.4}{The Ideal Intersection Property}}{15}{}%
\contentsline {section}{\tocsection {}{3}{Cartan Envelopes}}{16}{}%
\contentsline {section}{\tocsection {}{4}{Pseudo-Cartan Inclusions and their Cartan Envelopes}}{30}{}%
\contentsline {subsection}{\tocsubsection {}{\tpad4.1}{Definition of Pseudo-Cartan Inclusions and Examples}}{30}{}%
\contentsline {subsection}{\tocsubsection {}{\tpad4.2}{The AUP and Unitization of Cartan Envelopes}}{33}{}%
\contentsline {subsection}{\tocsubsection {}{\tpad4.3}{Some Properties Shared by a Pseudo-Cartan Inclusion and its Cartan Envelope}}{36}{}%
\contentsline {subsection}{\tocsubsection {}{\tpad4.4}{Constructing Pseudo-Cartan Inclusions from Cartan Inclusions}}{38}{}%
\contentsline {section}{\tocsection {}{5}{The Twisted Groupoid of the Cartan Envelope}}{40}{}%
\contentsline {subsection}{\tocsubsection {}{\tpad5.1}{Reprising the Unital Case}}{40}{}%
\contentsline {subsection}{\tocsubsection {}{\tpad5.2}{The Non-Unital Case}}{41}{}%
\contentsline {section}{\tocsection {}{6}{Constructions and Properties of Pseudo-Cartan Inclusions}}{43}{}%
\contentsline {subsection}{\tocsubsection {}{\tpad6.1}{Mapping Results}}{43}{}%
\contentsline {subsection}{\tocsubsection {}{\tpad6.2}{Inductive Limits}}{50}{}%
\contentsline {subsection}{\tocsubsection {}{\tpad6.3}{Minimal Tensor Products}}{53}{}%
\contentsline {section}{\tocsection {}{7}{Applications}}{55}{}%
\contentsline {section}{\tocsection {}{8}{Questions}}{59}{}%
\contentsline {section}{\tocsection {Appendix}{A}{}}{60}{}%
\contentsline {section}{\tocsection {Appendix}{}{References}}{62}{}%

\mysec[Introduction]{Introduction}\label{1intro}
 A famous theorem of Gelfand states that an abelian \cstaralg\ $\C$ may be
recognized as an algebra of continuous functions: $\C\simeq C_0(X)$,
where $X=\hat\C$ is the space of non-zero multiplicative linear
functionals on $\C$ equipped with the weak-$*$ topology.
In some cases, non-commutative \cstaralg s can also be represented as
functions on spaces.  For example,
Str{\u{a}}til{\u{a} and Voiculescu~\cite{StratilaVoiculescuReAFAl},
  showed that for a unital AF-algebra $\A=\indlim\A_n$ (each
  $\A_n$ is a finite dimensional \cstaralg), there exists a maximal
  abelian \cstar-subalgebra $\B\subseteq \A$ such that for every $n$,
  $\B\cap\A_n$ is maximal abelian in $\A_n$.  Str{\u{a}}til{\u{a} and
    Voiculescu used $\B$ to construct a ``coordinate system'' $X$ for $\A$ which
    behaves in much the same way as a system of matrix units for
    $M_n(\bbC)$; furthermore, the AF-algebra $\A$ can be viewed as a collection
    of functions on $X$.   Shortly after
    Str{\u{a}}til{\u{a} and Voiculescu's work, Feldman and
Moore~\cite{FeldmanMooreErEqReII} showed that if $\D$ is a $W^*$-Cartan
subalgebra of the (separably acting) von Neumann algebra  $\M$ and $\D\simeq L^\infty(X,\mu)$,
then there is a ``measured'' equivalence
relation $R$ on $X\times X$ along with a 2-cocycle $\sigma$ on $R$ such that
$\M$ is, loosely speaking, isomorphic to an algebra of functions, $M(R,\sigma)$, on $R$ with product
given by a ``generalized matrix multiplication''; further, the
isomorphism of $\M$ onto $M(R,\sigma)$ carries $\D$ onto an algebra
of functions supported on $\{(x,x): x\in X\}$.

Efforts to adapt the ideas of Feldman and Moore to \cstaralg s
 began with work of Renault in~\cite{RenaultGrApC*Al}, continued with
 Kumjian's work on $C^*$-diagonals in~\cite{KumjianOnC*Di}, and
culminated with Renault's definition of a Cartan MASA in a \cstaralg\
found in~\cite{RenaultCaSuC*Al}.  Given a Cartan MASA $\D$ in the
\cstaralg\ $\C$, the Kumjian-Renault theory produces a Hausdorff,
\'etale, and effective groupoid $G$, along with  a central groupoid extension
$\unit G\times \bbT\hookrightarrow \Sigma\twoheadrightarrow G$; this
 is the \cstar-algebraic analog of 
Feldman-Moore's equivalence relation $R$ and 2-cocycle $\sigma$.
Moreover, there is a unique isomorphism of $\C$ onto a
\cstaralg, $C^*_r(\Sigma,G)$, which carries $\D$ onto
$C_0(\unit{G})$. The algebra $C^*_r(\Sigma,G)$ can be viewed as an
algebra of functions on $\Sigma$ with a product generalizing matrix
multiplication, so here again, $\C$ may be viewed as an algebra of functions on a
space.  Conversely, every such central extension can be used to
construct a Cartan MASA in a \cstaralg, so the theory gives a complete
description of Cartan inclusions $\D\subseteq \C$.  These results are
presented in~\cite{RenaultCaSuC*Al}; for a less terse treatment,
see~\cite{SimsGroupoidsBook}.  While Renault and Kumjian worked in the
context of separable \cstaralg s, Raad showed in~\cite{RaadGeReThCaSu}
that this assumption can be removed.

Despite the success of the Kumjian-Renault theory, there are desirable
settings in the \cstaralg\ context where the axioms for a Cartan MASA
in the algebra $\C$ are not satisfied; examples include the virtual
Cartan inclusions introduced in~\cite[Definition~1.1]{PittsStReInI} and
the weak Cartan inclusions defined
in~\cite[Definition~2.11.5]{ExelPittsChGrC*AlNoHaEtGr}.  For 
these classes, the inclusion $\D\subseteq\C$ is regular
and $\D$ is abelian.    In such
settings, it is
possible to follow the Kumjian-Renault procedure to construct a
groupoid $G$ and a central groupoid extension as above.  However, when
this is done, the resulting groupoid $G$ may fail to be Hausdorff,
see~\cite[Theorem~4.4]{PittsStReInII}.  One approach to dealing with
such situations is to accept the fact that non-Hausdorff groupoids may
arise, and develop a theory which includes them.  This is the approach taken in
\cite{ExelPittsChGrC*AlNoHaEtGr}.

An alternate approach to studying regular inclusions of the form
$\D\subseteq \C$ with $\D$ abelian, is to seek embeddings of
$\D\subseteq \C$ into a Cartan inclusion $\B\subseteq \A$   
so that the dynamics of $\D\subseteq \C$ are
embedded into the dynamics of $\B\subseteq \A$.  (Such embeddings are
called regular embeddings.)  The idea is that the coordinates from the
larger inclusion 
$\B\subseteq \A$ might then be used to study $\C$.   In~\cite{PittsStReInI}, we gave a
characterization of those unital regular inclusions which can be
embedded into a \cstar-diagonal or a Cartan inclusion.
However, the \cstar-diagonal $\B\subseteq \A$ obtained from the construction
in~\cite{PittsStReInI} may bear little relation to the inclusion
$\D\subseteq \C$, and the
embedding of $\D\subseteq \C$ into $\B\subseteq \A$ is
not sufficiently rigid to transfer information from $\B\subseteq \A$
to $\D\subseteq \C$.

Part of the motivation for~\cite{PittsStReInII} was to remedy this
defect.  In~\cite{PittsStReInII}, the notion of the Cartan envelope is
defined and two characterizations of those unital regular inclusions
having a Cartan envelope are given.  The Cartan envelope described
in~\cite{PittsStReInII} is unique when it exists and is a minimal
Cartan pair generated by the original inclusion $\D\subseteq \C$.
Also, in~\cite[Section~7]{PittsStReInII}, a description of the
twisted groupoid associated to the Cartan envelope for
$\D\subseteq \C$ is given.

Many regular inclusions are not unital.  Examples include inclusions
from graph algebras (or higher rank graph algebras), from crossed
product constructions, or arise from tensor product constructions
(such as those in Section~\ref{sec:tenprod} below).  It is
therefore of interest to extend results on Cartan envelopes to the
non-unital setting.

Three primary goals of the present paper are: first, to establish results on
Cartan envelopes which apply  in the non-unital context; second,
to give results showing that the Cartan envelope of a regular
inclusion carries information about  the original inclusion; and third, to
give constructions of \pd s along with their Cartan envelopes.    Extending the characterization of regular inclusions having a Cartan
envelope to include non-unital regular inclusions is
unexpectedly subtle.
On first glance, one might expect that
passing from the unital setting to the non-unital case would simply be
a matter of adjoining units and then applying the results on Cartan
envelopes for unital
inclusions from~\cite{PittsStReInII}.  However, the unitization
process need not preserve normalizing elements (see
Definition~\ref{innore}\eqref{innore3}), nor does it preserve regular
mappings (Definition~\ref{innore}\eqref{innore7}). For a concrete
example, suppose $\fJ\idealin C_0(\bbR)$ is the ideal of functions in
$C_0(\bbR)$ which vanish on $[-1,1]$.  Any $f\in C_0(\bbR)$ is a
normalizer for the inclusion $\fJ\subseteq C_0(\bbR)$, yet $f$ is not
a normalizer for the unitized inclusion
$\tilde\fJ \subseteq (C_0(\bbR))^{\sim}$ unless $|f|$ is constant on
$[-1,1]$.  In particular, this shows that while the identity mapping
on $C_0(\bbR)$ is a regular map, its extension to $(C_0(\bbR))^{\sim}$
(that is, $\Id_{(C_0(\bbR))^{\sim}}$) is not regular.  Furthermore,
\cite[Example~3.1]{PittsNoApUnInC*Al} gives an example of an inclusion
$\D\subseteq \C$ which is not regular, but whose unitization
$\tilde\D\subseteq \tilde\C$ is a Cartan inclusion.  Finally, a vexing
issue we have been unable to resolve is whether regularity for an
inclusion is preserved when units are adjoined.  The ill behavior of
regularity properties under the unitization process leads to serious
technical difficulties when considering Cartan envelopes for
non-unital inclusions. Using tools found in~\cite{PittsNoApUnInC*Al},
we develop techniques to overcome those obstacles and are able to
characterize those regular inclusions having a Cartan envelope; we
shall describe the characterization in a moment.

Our work on Cartan envelopes yields an alternate set of axioms for a Cartan MASA
$\D\subseteq \C$.    Recall that an inclusion
$\B\subseteq \A$ of \cstaralg s has the \iip\ (iip) if every non-zero
ideal of $\A$ has non-trivial intersection with $\B$.   We also use
$\B^c=\A\cap \B'$ for the relative commutant of $\B$ in $\A$.  The
columns of the following table give Renault's axioms for a Cartan MASA
and 
(equivalent) alternate axioms for a Cartan MASA. 
\begin{center}
  \begin{tabular}{|l|l|} \hline
    \bf Cartan Inclusions:  & \bf Cartan Inclusions:\\
    \bf Renault's Axioms     & \bf Alternate Axioms
    \\
    \hline
   $\D\subseteq \C$ is regular& $\D\subseteq \C$ is  regular\\\hline
    $\D$ is a MASA in $\C$ & 
                             $\D$ and  $\D^c$ are abelian, and both \\
                                              & the inclusions
                                                $\D\subseteq \D^c$ and
                                                $\D^c\subseteq \C$\\
    &
                                                 have the \iip\\
                                                  \hline
    $\exists$ a \textit{faithful} conditional expectation &
                                                                       $\exists$
                                                                       a 
                                                           conditional 
                                                           expectation
                                                         \\
    $E:\C\rightarrow \D$&$E:\C\rightarrow \D$
    \\
    \hline
  \end{tabular}
\end{center}
(Renault's original definition of Cartan MASA also assumes  $\D$
contains an approximate identity for $\C$,
but~\cite[Theorem~2.4]{PittsNoApUnInC*Al} shows that assumption can be removed.)
By keeping the first two entries of the right-hand column but not
requiring the third, we obtain a class of inclusions which we call
\textit{\pd s}.

The characterization theorem, Theorem~\ref{!pschar}, extends
the characterizations of the existence of a Cartan envelope found for
unital inclusions in~\cite{PittsStReInII} to include the non-unital
case.  Theorem~\ref{!pschar} shows the equivalence of the following:
the inclusion $\D\subseteq \C$ (again with $\D$ abelian) has a Cartan
envelope; $\D\subseteq \C$ is a \pd; and $\D\subseteq \C$ has the
faithful unique pseudo-expectation property (see
Definition~\ref{f!pseDef}).  That an inclusion is a \pd\ if and only
if it has the faithful unique pseudo-expectation property is the
reason for the term \pd.
A description of the Kumjian-Renault twist
$\unit{G}\times\bbT\hookrightarrow \Sigma\twoheadrightarrow G$ for the
Cartan envelope of $\D\subseteq \C$ is given in
Theorem~\ref{gtwCE}.  This description exhibits $\Sigma$ as a family of
linear functionals on $\C$ and $G$ is the quotient of $\Sigma$ by a
certain action of $\bbT$.

The class of \pd s contains the
following classes: Cartan
inclusions; the virtual Cartan inclusions of~\cite{PittsStReInI}; the
weak Cartan inclusions introduced in~\cite{ExelPittsChGrC*AlNoHaEtGr};
and inclusions of the form $J\subseteq C_0(X)$, where $X$ is a locally
compact Hausdorff space and $J$ is an essential ideal in $C_0(X)$.

 \mychap~\ref{prelim} sets terminology, notation, and also gives a number of
foundational results needed throughout the paper.  These results concern \wnd\ inclusions, pseudo-expectations,
and the \iip.   Some of the results in \mychap~\ref{prelim}, such as 
Theorem~\ref{f!ps4rin}, illustrate the technical challenges which can
arise when verifying that
desirable properties of an inclusion are preserved under the
unitization process.

After presenting our results on existence and uniqueness of
Cartan envelopes in \mychap~\ref{mainr}, we explore properties of \pd
s and relationships between a \pd\ and its Cartan envelope.
\mychap~\ref{Sec: pd} contains the definition of  \pd s 
along with a few properties shared by $\C$ and $\A$ when $\D\subseteq
\C$ is a \pd\ having Cartan envelope $\B\subseteq \A$.  For example, Proposition~\ref{simple} shows $\C$ is simple if and
only if $\A$ is simple.   While a
\pd\  can be (minimally) embedded into a Cartan inclusion,
Proposition~\ref{pdcons} gives a reverse process: it shows how to construct a family of \pd s
from a given Cartan MASA $\B\subseteq \A$, and also
examines the Cartan envelope for \pd s  arising from this
construction.  

\mychap~\ref{tgpC} constructs the Kumjian-Renault
groupoid model for the Cartan envelope of a given \pd\
$\D\subseteq \C$ starting with  a family of linear functionals on $\C$.
In particular, this allows us to explicitly see the embedding of
$\D\subseteq \C$ into its \cstar-envelope $\B\subseteq\A$ by realizing $\C$ as an
algebra of functions (with a convolution product as multiplication).

\mychap~\ref{const} contains a rigidity property of the Cartan
envelope and explores some permanence properties of \pd s and their
Cartan envelopes.  The rigidity property is Proposition~\ref{autoext}:
it shows that given a \pd\ $\D\subseteq \C$ with Cartan envelope
$\B\subseteq \A$, any regular, $\D$-preserving $*$-automorphism of $\C$ extends
uniquely to a $*$-automorphism of $\A$.  We show in Theorem~\ref{indlim} that an
inductive limit of \pd s with \text{suitable} connecting maps is again
a \pd\  and the Cartan envelope of such an inductive limit is
the inductive limit of the Cartan envelopes of the approximating \pd
s.  The proof of Theorem~\ref{indlim} uses an important mapping
property of Cartan envelopes, Theorem~\ref{cmap}.  Theorem~\ref{gentp}
shows that the minimal tensor product of two \pd s is a \pd, and the
Cartan envelope of their minimal tensor product is the minimal tensor
product of their Cartan envelopes.  Because quotients of \pd s by
regular ideals are the subject of a forthcoming paper (by Brown,
Fuller, Reznikoff and the author), we do not consider quotients here.

We provide some applications of our work in \mychap~\ref{apps}. We
show that when $\D\subseteq \C$ is a unital \pd, the \cstar-envelope
of a Banach algebra $\D\subseteq \A\subseteq \C$ is the
\cstar-subalgebra of $\C$ generated by $\A$ and that $\D$ norms $\C$
in the sense of~\cite{PopSinclairSmithNoC*Al}.  We combine these
results to show that if, for $i=1,2$: $\D_i\subseteq \C_i$ are
 \pd s such that their unitizations $\tilde\D_i\subseteq \tilde\C_i$
 are also \pd s; $\A_i$ are Banach algebras with
$\D_i \subseteq \A_i\subseteq \C_i$; and $u:\A_1\rightarrow \A_2$ is
an isometric isomorphism, then $u$ uniquely extends to a
$*$-isomorphism of $C^*(\A_1)$ onto $C^*(A_2)$.  While these
applications have antecedents in the literature, our results are
significant generalizations.  We include a few open questions in \mychap~\ref{Sec:Ques}.

We require~\cite[Theorem~6.9]{PittsStReInII} at two points in our
work: in Section~\ref{ruc}\eqref{ruc3}, and in the
proof of Lemma~\ref{GNS0}.  While \cite[Theorem~6.9]{PittsStReInII} is correctly
stated in~\cite{PittsStReInII},  due to an error in the statement
of~\cite[Lemma~2.3]{PittsStReInII}, the proof
of~\cite[Theorem~6.9]{PittsStReInII} is insufficient.   The correct statement of~\cite[Lemma~2.3]{PittsStReInII} and a
complete proof of~\cite[Theorem~6.9]{PittsStReInII} may be found in
Appendix~\ref{AppendixA}. 

We thank Adam Fuller for several useful comments and Alex Kumjian for
a  comment regarding the history of the adaptation of the
Feldman-Moore theory to the $C^*$-algebraic context.

The author received helpful reports from two different referees, and
is grateful to both referees.  \referee 's report was particularly
thorough: it contained dozens of comments, identified gaps in some
proofs, gave suggestions for streamlining arguments (e.g.\ the
proofs of Lemma~\ref{pseprop},
Theorem~\ref{indlimEnv}\eqref{indlim3.14}, and Theorem~\ref{gentp}), and provided
numerous other helpful remarks.  We are indebted to \referee\ for their interest
and careful reading of the manuscript; their comments greatly improved the
paper and exposition.   Of course, any mistakes in the
paper are solely the responsibility of the author.

\mysec[Foundations]{Foundations}\label{prelim}
\numberwithin{equation}{subsection}

\subsection{Terminology and Notation}
For a normed linear space $X$, we shall use $\dual{X}$ to denote the collection of
all bounded linear functionals on $X$.

Let $\A$ be a Banach algebra.
When $\A$ is unital, we denote its unit by $I$ or $I_\A$.  Unless
explicitly stated otherwise, when referring to an ideal $J$ in $\A$,
we will always assume $J$ is closed and two-sided.  When $J$ is an
ideal in $\A$, we write $J\idealin \A$.  The ideal $J$ is called an
\textit{essential ideal} if for any non-zero ideal $K\idealin A$,
$J\cap K\neq \{0\}$.

We  follow the notation and conventions regarding unitization of a
\cstaralg\ found in \cite[II.1.2.1]{BlackadarOpAl}.  For a \cstaralg\
$\A$, let $\A^\dag:=\A\times \bbC$ with the usual operations and norm which make
it into a unital \cstaralg.
Define $\tilde\A$ and a
$*$-monomorphism $u_\A:\A\rightarrow \tilde\A$ by
  \begin{align}
    \tilde\A&:=
              \begin{cases} \A& \text{if $\A$
      is unital}\\ \A^\dag& \text{if $\A$ is not
      unital;}\end{cases} \label{unitizedef}\\
    \intertext{and for $x\in\A$,}
  u_\A(x)&:=\begin{cases} x & \text{if $\A$
      is unital}\\  (x,0) &\text{if $\A$ is not
      unital.}
  \end{cases} \label{unitizedef1}
  \end{align}
  We shall consider the trivial algebra $\{0\}$ to be unital.
Also, when there is no danger of confusion, we will identify $\A$ with
$u_\A(\A)$, and regard $\A\subseteq \tilde\A$.

Suppose $\A_1$ and $\A_2$ are 
\cstaralg s.
 Given a
 bounded
linear map $\psi:\A_1\rightarrow \A_2$, the map
$\tilde\psi:\tilde \A_1\rightarrow \tilde\A_2$ given by
\begin{equation}\label{stext}
  \tilde\psi(x)=\begin{cases} u_{\A_2}(\psi(x)) & \text{if $\A_1$ is unital;}\\
    u_{\A_2}(\psi(a))+\lambda I_{\tilde\A_2}& \text{if $\A_1$ is not unital and
      $x=(a,\lambda)\in \tilde\A_1$}\end{cases}
\end{equation}
will be called the \textit{standard extension \index{Standard
    extension} of $\psi$} to $\tilde\A_1$.  We will frequently use the
following property of the standard extension without comment:
\begin{equation}\label{stext1}
  u_{\A_2}\circ\psi=\tilde\psi\circ u_{\A_1}.
\end{equation}
It is not
generally the case that $\tilde\psi(I_{\tilde\A_1})=I_{\tilde A_2}$.
 When $\psi$ is
a $*$-homomorphism, so is $\tilde\psi$.  

\begin{flexstate}{Observation}{} \label{ccp} If $\psi$ is contractive and
  completely positive, then $\tilde\psi$ is also contractive and
  completely positive.
\end{flexstate}
\begin{proof}
When $\A_1$ is not unital, this follows from~\cite[Proposition
2.2.1]{BrownOzawaC*AlFiDiAp} applied to $u_{\A_2}\circ\psi$, and from the definition of $\tilde\psi$ when $\A_1$ is 
unital.
\end{proof}

Fundamental to
our study are inclusions, which we usually consider as a \cstaralg\
$\B$ contained in another \cstaralg\ $\A$.
At times, particularly in this section, it will be helpful to explicitly
mention the embedding of $\B$ into $\A$.   While this makes notation more
cumbersome, the benefit of clarity outweighs the notational burdens. 
\begin{definition}\label{incdef}
  \begin{enumerate}
    \item \label{incdef1} An \textit{inclusion}\index{Inclusion} is a triple
  $(\A,\B, f)$, where $\A$ and $\B$ are \cstaralg s and
  $f:\B\rightarrow \A$ is a $*$-monomorphism.  We include the
  possibility that $\B=\{0\}$.
  (Starting in Section~\ref{mainr}, we will always assume
  $\B$ is abelian, but until then we do not make that restriction.)

\item\label{incdef2} Let $(\A_1,\B_1,f_1)$ be an
  inclusion.  An inclusion $(\A_2,\B_2,f_2)$ together with a
  $*$-monomorphism $\alpha: \A_1\rightarrow \A_2$ such that
  $\alpha(f_1(\B_1))\subseteq f_2(\B_2)$ will be called an
  \textit{expansion of $(\A_1,\B_1,f_1)$\index{Expansion}} and will be denoted by
  $(\A_2,\B_2, f_2\ms\alpha)$.  Because $f_2$ is one-to-one,
  $(\A_2,\B_2,f_2\ms\alpha)$ is an expansion of $(\A_1,\B_1,f_1)$ if and
  only if $\alpha:\A_1\rightarrow \A_2$ is a $*$-monomorphism and
  there exists a $*$-monomorphism
  $\underline\alpha :\B_1\rightarrow\B_2$ such that
  \begin{equation}\label{itiso} \alpha\circ
    f_1=f_2\circ\underline\alpha.
      \end{equation}

 \item\label{incdef3}  We will call the two inclusions
  $(\A_1,\B_1,f_1)$ and $(\A_2,\B_2,f_2)$ \textit{isomorphic
    inclusions} if there
is a $*$-isomorphism $\alpha: \A_1\rightarrow
  \A_2$ such that $\alpha( f_1(\B_1))=f_2(\B_2)$.

\item \label{innore2} The inclusion $(\A,\B,f)$ is called a
  \textit{unital inclusion}\index{Inclusion!unital} if $\A$ and $\B$
  are unital and $f(I_\B)=I_\A$.
\item\label{incdef4} Given the inclusion $(\A,\B,f)$,  
the standard extension of $f$ to $\tilde
f:\tilde\B\rightarrow\tilde\A$ described in \eqref{stext} is a $*$-monomorphism 
 because $\B$ is an
essential ideal in $\tilde\B$.  Thus $(\tilde\A,\tilde\B, \tilde f)$ is an
inclusion.  We will
call $\tilde f$ 
\textit{the standard embedding} of $\tilde\B$ into $\tilde\A$.
(Caution:
$(\tilde\A,\tilde\B,\tilde f)$ need not be a unital inclusion.)

Because
\begin{equation}\label{ud1}\tilde f\circ u_\B=u_\A\circ
  f,\end{equation} $(\tilde\A,\tilde\B, \tilde f\ms u_\A)$ is an expansion of
$(\A,\B,f)$.

\item \label{incdef5} For $1\leq i<j\leq 3$, let $(\B_j, \B_i, f_{ji})$ be inclusions.
  If $f_{31}=f_{32}\circ f_{21}$, we will say the two inclusions $(\B_2,\B_1, f_{21})$
  and $(\B_3,\B_2,f_{32})$ are \textit{intermediate
    inclusions}\index{Inclusion!intermediate}  for
  $(\B_3,\B_1,f_{31})$.  \begin{equation}\label{intCD}
    \xymatrix{ & \B_3\\ \B_1\ar[ur]^{f_{31}}\ar[r]_{f_{21}}
      &\B_2\ar[u]_{f_{32}}\,\,.}
  \end{equation}
\end{enumerate}

\end{definition}

\begin{remark}{Notation}\label{idst}
Let $(\A,\B,f)$ be an inclusion.
  \begin{itemize}
\item  When there is little danger of
confusion, we will often identify $\B$ with its
image $f(\B)$, so that $\B\subseteq \A$.  When this identification is made, we
will suppress $f$ and call $(\A,\B)$ an inclusion.

  \item When we write
$(\tilde\A,\tilde\B)$, we always mean $(\tilde\A,\tilde\B,\tilde
f)$. Since $\tilde f$ is not necessarily a unital
mapping,  $(\tilde\A,\tilde\B)$ need not be a unital inclusion in
the sense of Definition~\ref{incdef}\eqref{innore2}.

\item We will sometimes write
$(\A,\B,\subseteq)$ for inclusions when  $\B$ is a \cstar-subalgebra
of $\A$.

\item If $(\A_2,\B_2, f_2\ms \alpha)$ is an expansion of $(\A_1,\B_1,f_1)$ and
$\B_i$ has been identified with $f_i(\B_i)$ as above, we will again
suppress writing 
$f_i$ and say $(\A_2,\B_2\ms \alpha)$ is an expansion of $(\A_1,\B_1)$.

\item If  $(\B_2,\B_1,f_{21})$ and $(\B_3,\B_2, f_{32})$ are
inclusions which are intermediate for $(\B_3,\B_1,f_{31})$, we will sometimes
say that 
\textit{$\B_2$ is intermediate to $(\B_3,\B_1)$} and may also indicate
this using the notation $\B_1\subseteq \B_2\subseteq \B_3$.
\end{itemize}
\end{remark}

There are: various types of inclusions; objects associated to an
inclusion; and  properties an inclusion may have.  We describe
some of them here.
\begin{remark}{Definitions, Notations, and Comments} \label{innore}  Let
$(\A,\B,f)$ be an inclusion.
\begin{enumerate}
\item\label{innore-1}  We will use $\rcom(\A,\B,f)$ for the
  \textit{relative commutant}\index{Relative commutant}
  of $f(\B)$ in $\A$, that is,
  \[\rcom(\A,\B,f):=\{a\in \A: af(b)=f(b)a \text{ for all } b\in
\B\}.\]   When $\B$ is identified with $f(\B)$, we will frequently use
the notation $\B^c$ or $\rcom(\A,\B)$ instead of $\rcom(\A,\B,f)$.
\item\label{innore0} We say that $(\A,\B,f)$ has the 
  \textit{\iip}\index{Ideal intersection property} if every non-zero
  ideal in $\A$ has non-trivial intersection with $f(\B)$.  We will
  also use the term \textit{essential
    inclusion}\index{Inclusion!essential} as a synonym for an
  inclusion with the \iip.  Likewise, when $(\A,\B, f)$ is an
  essential inclusion, we will sometimes call the map $f$ an
  \textit{essential map}.
The literature contains several synonyms for the \iip, for example
  in~\cite{KwasniewskiMeyerApAlExPrUnPsEx}, an inclusion with the
  \iip\ is said to \textit{detect
    ideals}, and in \cite{PittsZarikianUnPsExC*In}, such an
  inclusion is called \textit{$C^*$-essential}.

  The reader is cautioned that our use of `essential' here conflicts with
  usage elsewhere. 
  
\item \label{innore1} $(\A,\B,f)$ has the \textit{approximate unit
    property}\index{Approximate unit property (AUP)}
  (abbreviated AUP) if there is an approximate unit
$(u_\lambda)$ for $\B$ such that $f(u_\lambda)$ is an approximate unit
  for $\A$.  
  In  remarks following~\cite[Definition~2.1]{ExelOnKuC*DiOpId}, Exel notes
that the Krein-Milman Theorem
and~\cite[Lemma~2.32]{AkemannShultzPeC*Al} yield the following
characterization of the approximate unit property.
\begin{flexstate}{Fact}{see \cite{ExelOnKuC*DiOpId}}
  \label{AUPChar} An
  inclusion $(\A,\B,f)$ has the approximate unit property if and only if no pure
  state of $\A$ annihilates $f(\B)$.
\end{flexstate}
\item \label{innore3} A \textit{normalizer}\index{Normalizer} for $(\A,\B,f)$ is an element of the set,
  \[\N(\A,\B,f):=\{v\in \A: vf(\B) v^*\cup v^*f(\B) v\subseteq
    f(\B)\}.\]  If $\B$ is identified with $f(\B)$, we will write
    $\N(\A,\B)$ rather than $\N(\A,\B,f)$.
\item\label{innore3.5} Closely related to normalizers are
  intertwiners.  An element $w\in \A$ is called an
  \textit{intertwiner}\index{Intertwiner} if the sets $f(\B)w$ and $wf(\B)$ coincide.  We will write $\I(\A,\B,f)$
  (or $\I(\A,\B)$ when $f$ is suppressed) for the set of all
  intertwiners.
\item \label{innore4} $(\A,\B,f)$ is said to be a \textit{regular
  inclusion}\index{Inclusion!regular} if
$\overline{\spn} \,\N(\A,\B,f)=\A$.
\item \label{innore5} If $f(\B)$ is maximal abelian in $\A$, we will call
  $(\A,\B,f)$  a \textit{MASA
    inclusion}.\index{Inclusion!MASA}
\item \label{innore6}$(\A,\B,f)$ is a \textit{Cartan inclusion},\index{Inclusion!Cartan} also
  called a \textit{Cartan pair}, if it is a regular MASA inclusion and
  there exists a faithful conditional expectation
  $\Delta:\A\rightarrow f(\B)$.  Cartan inclusions necessarily have the
  AUP,~\cite[Theorem~2.6]{PittsNoApUnInC*Al}. 
  
  The conditional expectation $\Delta$ is unique;
  see~\cite[Corollary~5.10]{RenaultCaSuC*Al}
  or~\cite[Theorem~3.5]{PittsStReInI}.  Also, $\Delta$ is invariant
  under $\N(\A,\B,f)$ in the sense that for every $v\in \N(\A,\B,f)$ and
  $x\in \A$,
  \begin{equation}\label{exinv}
    \Delta(v^*xv)=v^*\Delta(x) v.
  \end{equation}
  While we do not have an explicit reference for~\eqref{exinv}, it is
  certainly known;  for example, it follows from Renault's (twisted) groupoid model for a
  Cartan pair.

  A remark on the term `invariant' to describe a conditional
  expectation satisfying~\eqref{exinv} seems appropriate here.  By
  considering the conditional expectation onto $C_0(X)$ for
  inclusions of the form $(C_0(X)\rtimes_r\Gamma,
  C_0(X))$ where $X$ is a locally compact Hausdorff space and $\Gamma$ is a discrete group, one might
  choose 
  the term `equivariant' instead of `invariant'.  However,
  `invariant' is used in the literature, see for
  example,~\cite[Definition~3.14]{BrownFullerPittsReznikoffReIdIdInQu},
  so we have elected to use `invariant'.
\begin{remark}{Comment}\label{CariffUCar}
  We will repeatedly use the fact that the regular inclusion $(\A,\B,f)$
  is a Cartan inclusion if and only if $(\tilde \A,\tilde \B,\tilde f)$ is a
  Cartan inclusion, see \cite[Proposition~3.2]{PittsNoApUnInC*Al}.
  Example~3.1 of~\cite{PittsNoApUnInC*Al} shows the necessity of
  the regularity hypothesis on $(\A,\B,f)$.
\end{remark}
\item \label{innore7} If $(\A_1,\B_1,f_1)$ is another inclusion, a $*$-homomorphism $\alpha:
\A\rightarrow \A_1$ is a \textit{regular homomorphism}\index{Regular
  homomorphism}
if
\[\alpha(\N(\A,\B,f))\subseteq \N(\A_1,\B_1,f_1).\]
We will sometimes indicate this by saying
$\alpha:(\A,\B,f)\rightarrow (\A_1,\B_1,f_1)$ is a regular
homomorphism.
\end{enumerate}
\end{remark}

When $\B_2$ is intermediate to $(\B_3,\B_1)$, it may happen that
$\tilde\B_2$ is not intermediate to $(\tilde\B_3,\tilde\B_1)$.  Here
is an example of this behavior (the example also provides an example of an inclusion
whose ``unitization'' is not unital).
\begin{remark}{Example}\label{uincng}
Let $S_0$ and $S_1$ be isometries on a Hilbert space such that
$S_0S_0^*+S_1S_1^*=I$, let $\S$ be the inverse semigroup of partial
isometries consisting of all finite products of elements of the set $\{S_0, S_1,
S_0^*, S_1^*\}$, and let
$\D$ be the \cstaralg\ generated by the projections in $\S$.
There is a unique multiplicative
linear functional $\tau:\D\rightarrow \bbC$ on $\D$ such that for
every $n\in \bbN$, $\tau(S_1^nS_1^n{}^*)=1$. (Indeed,
$\D$ is isomorphic to the continuous functions on the Cantor middle
thirds set, and under this isomorphism, $\tau$ corresponds to
evaluation at 1.)  Let
\[\B_3=\B_1=\ker\tau, \quad \B_2=\D,\] and let 
$\sigma:\D\rightarrow \ker\tau$ be the map \[\sigma(d)=S_0d S_0^*.\] Put 
\[f_{31}:=\sigma^2|_{\B_1},\quad
  f_{21}:=\sigma|_{\B_1}\dstext{and}f_{32}:=\sigma.\]  Since
$f_{31}=f_{32}\circ f_{21}$, 
$(\B_2,\B_1,f_{21})$ and $(\B_3,\B_2, f_{32})$ are  intermediate
inclusions for 
$(\B_3,\B_1, f_{31})$.
Let  $q$ denote the quotient map of $\tilde \B_3$ onto $
  \tilde\B_3/\B_3\simeq \bbC$.   Since $\B_1$ and $\B_3$ are not
  unital, 
  $q\circ \tilde f_{31}\neq 0$.  On the other hand, since $\B_2$ is
  unital, 
  $q\circ \tilde
    f_{32}\circ \tilde f_{21}=0$.    Therefore,   $\tilde f_{31}\neq \tilde
    f_{32}\circ \tilde f_{21}$, so $(\tilde\B_2,\tilde\B_1,\tilde
    f_{21})$ and $(\tilde \B_3, \tilde\B_2, \tilde f_{32})$ are  not
    intermediate inclusions for 
    $(\tilde\B_3,\tilde\B_1, \tilde f_{31})$.  Finally, note that $(\B_2,\B_1,
    f_{21})$ is an example of an inclusion where $(\tilde\B_2, \tilde\B_1,\tilde
    f_{21})$ is not a unital inclusion.
\end{remark}

\subsection{Weakly Non-Degenerate Inclusions}

The approximate unit property for an inclusion $(\A,\B,f)$ has
been used in the literature as a definition of non-degeneracy, see for
example~\cite[Definition~1.4]{CrytserNagySiCrEtGrC*Al}
or~\cite[Definition~2.3.1(iii)]{ExelPittsChGrC*AlNoHaEtGr}.  However, when $J\idealin \A$ is a proper
and essential ideal, the inclusion $(\A, J,\subseteq)$ cannot have the AUP, so
using the AUP as a non-degeneracy condition excludes such examples.
For our purposes, the following weaker version of non-degeneracy is
more appropriate.
\begin{definition}\label{newnondeg} Let $(\A,\B,f)$ be  an inclusion.
  \begin{enumerate}
  \item Define the \textit{annihilator} of $\B$ in $\A$ to be the set,
  \[\ann(\A,\B,f):=\{a\in \A: af(b)=0=f(b)a \text{ for all } b\in
    \B\}.\] Notice that $\ann(\A,\B,f)$ is  a \cstar-subalgebra of $\A$
  and also  an ideal of
  $\rcom(\A,\B,f)$. 
 \item  An inclusion
   $(\A,\B,f)$ such that $\ann(\A,\B,f)=\{0\}$ will be called a
   \textit{\wnd\ inclusion}.\index{Inclusion!weakly non-degenerate}   

\end{enumerate}
As usual, when $\B$ is identified with $f(\B)\subseteq \A$, we will
write $\ann(\A,\B)$ instead of $\ann(\A,\B,f)$ and will say $(\A,\B)$
is \wnd\ when $(\A,\B,f)$ is \wnd.
\end{definition}

\begin{remark}{Remark} We chose the term ``\wnd'' instead of ``non-degenerate'' to
avoid confusion with existing literature.  Because of the similarity
between the terms ``\wnd'' and ``non-degenerate'' we will always use the
term ``approximate unit property'' instead of ``non-degenerate.''
Finally, to avoid possible confusion, for a representation $\pi$ of a \cstaralg\ $\A$ on a Hilbert
space $\H$, we will sometimes explicitly write
$\overline{\pi(\A)\H}=\H$ instead of using the term ``non-degenerate representation.''
\end{remark}

Here are some conditions ensuring an inclusion is \wnd.  
\begin{lemma}\label{->ndeg}  Let $(\A,\B)$ be an inclusion.
  \begin{enumerate}
    \item\label{->ndeg0}  If $(\A,\B)$ has the AUP, then $(\A,\B)$ is \wnd.
\item\label{->ndeg0.5}  If $\B\idealin \A$ is an essential ideal, then
  $(\A,\B)$ is \wnd.  In particular, when $\A$ is not unital, $(\tilde
  \A,\A, u_\A)$ is \wnd.
    \item\label{->ndeg1}  If $\B$ is contained in the center of $\A$  and $(\A,\B)$
    has the \iip, then $(\A,\B)$ is \wnd.  
  \item\label{->ndeg2}  Suppose $E:\A\rightarrow \B$ is a faithful
  conditional expectation.  Then $(\A,\B)$ is \wnd.
\item\label{->ndeg3}  Suppose  $(\C,\D)$ is an inclusion such
        that $\D$ is abelian and $(\D^c,\D)$ has the
        \iip. 
        Then $(\C,\D)$ is \wnd.   In particular, any MASA inclusion is \wnd.
      \end{enumerate}
  \end{lemma}
  \begin{proof}
    \eqref{->ndeg0} Let $(u_\lambda)\subseteq \B$ be an approximate
    unit for $\A$.  If $x\in \ann(\A,\B)$, then $x=\lim xu_\lambda=0$.

\eqref{->ndeg0.5}  Since $\B$ is an ideal of $\A$, $\ann(\A,\B)$ is
also an ideal in $\A$.   Clearly $\Ann(\A,\B)\cap \B=\{0\}$, and as
$\B$ is an essential ideal, $\Ann(\A,\B)=\{0\}$.
    
    \eqref{->ndeg1} By hypothesis, $\B^c=\A$, so $\ann(\A,\B)$ is an
    ideal of $\A$.  As $\ann(\A,\B)$  has trivial intersection with $\B$,  the \iip\
    gives $\ann(\A,\B)=\{0\}$.  Therefore,  $(\A,\B)$ is \wnd.
    
\eqref{->ndeg2}  Let $x\in \Ann(\A,\B)$.  Since 
$x^*x\in
\Ann(\A,\B)$ and $E(x^*x)\in \B$, 
\[E(x^*x)^2=E(x^*xE(x^*x))=0,\] with the first equality following from
Tomiyama's theorem, see~\cite[Theorem~1.5.10]{BrownOzawaC*AlFiDiAp}.   This gives $x=0$ by
faithfulness of $E$.  So $(\A,\B)$   is \wnd.

\eqref{->ndeg3}
Applying part~\eqref{->ndeg1} to $(\D^c,\D)$ gives 
 $\ann(\D^c,\D)=\{0\}$.  But $\ann(\C,\D)$ is an ideal in $\D^c$
 having trivial intersection with $\D$.  So
 $\ann(\C,\D)=\{0\}$, showing $(\C,\D)$  is \wnd.
\end{proof}

In general, for an inclusion $(\A,\B,f)$, there are multiple
embeddings of 
$\tilde\B$ into $\tilde\A$
which extend $u_\A\circ f$, 
and it may happen that the image of $I_{\tilde\B}$ under the
standard embedding is not  $I_{\tilde\A}$. 
Fortunately, these behaviors cannot occur when $(\A,\B,f)$ is \wnd, as
we shall see in parts \eqref{sameunit5} and \eqref{sameunit6} of
Lemma~\ref{sameunit}
below.   Before proving Lemma~\ref{sameunit}, we
require some preparation.   The first is the following well-known fact (see
\cite[Lemma~4.8]{HamanaReEmCStAlMoCoCStAl} or \cite[(2.2)]{EffrosAsNoOr}).
\begin{lemma}\label{positive}  Suppose $\H=\H_1\oplus\H_2$ is the
  direct sum of the Hilbert spaces $\H_1$ and $\H_2$.  Let $A\in
    \B(\H_1)$ be positive and invertible, $B\in \B(\H_2,\H_1)$, and
    $0\leq C\in \B(\H_2)$.
 Then 
  \[T:=\begin{bmatrix} A&B\\ B^* & C\end{bmatrix} \in \B(\H_1\oplus
    \H_2)\] is a
  positive operator if and only if $B^*A^{-1}B\leq C$.  Moreover, if
  $T\geq 0$ and $C=0$, then $B=0$.
\end{lemma}
\begin{proof}
  Both statements follow from the factorization,
  \[\begin{bmatrix}
      A&B\\ B^*& C\end{bmatrix}=
    \begin{bmatrix} I_{\H_1} & A^{-1}B\\ 0& I_{\H_2}\end{bmatrix}^*
    \begin{bmatrix} A& 0\\ 0&
      C-(A^{-1/2}B)^*(A^{-1/2}B)\end{bmatrix}
    \begin{bmatrix} I_{\H_1} & A^{-1}B\\ 0& I_{\H_2}\end{bmatrix}.\]
\end{proof}

Among inclusions, \wnd\ inclusions are rather well-behaved.  Here are
a number of desirable properties which \wnd\ inclusions possess.  
\begin{lemma}\label{sameunit}  Suppose $(\A,\B,f)$ is a \wnd\ 
  inclusion.  The following statements hold.
  \begin{enumerate}
    \item \label{sameunit1}
    If $\B$ is unital, then $\A$ is unital and $f(I_\B)$ is the unit for $\A$.
\item\label{sameunit2}  Suppose $a\in \A$ satisfies $f(b)a=af(b)=f(b)$ for
  every $b\in \B$.  Then $\A$ is unital and $a=I_\A$.
\item \label{sameunit3}  Let $(u_\lambda)$ be an approximate unit for
  $\B$.   Suppose $a\in \A$ satisfies:  $0\leq a$,  
  $\norm{a}\leq 1$, and  $f(u_\lambda)\leq a $ for every $\lambda$.  Then $\A$
  is unital and $a=I_\A$. 
\item\label{sameunit4}  $(\tilde \A, \B, u_\A\circ f)$ is \wnd.
\item\label{sameunit5}  $(\tilde\A,\tilde\B, \tilde f)$ is a unital
  inclusion; in particular, it is \wnd.
\item\label{sameunit6}   If $\sigma: \tilde\B\rightarrow \tilde\A$ is
  a $*$-monomorphism such that $\sigma\circ u_\B=u_\A\circ f$, then
  $\sigma=\tilde f$.  
\end{enumerate}
\end{lemma}
\begin{proof}
 \eqref{sameunit1}  We may suppose $\A\subseteq \bh$.  Let $I=I_\H$, $p:=f(I_\B)$, and choose
  $a\in \A$ with $0\leq a$.   Relative to the decomposition,
  $\H=p\H+(I-p)\H$, $a$ and $a+p$ have the form,
\[a=\begin{bmatrix} pap & pa(I-p)\\ (I-p) a p& (I-p)
    a(I-p)\end{bmatrix}\dstext{and}
  a+p=\begin{bmatrix} pap+p & pa(I-p)\\ (I-p) a p& (I-p)
    a(I-p)\end{bmatrix}.\] Let  
$x:=(I-p) a(I-p)= a-ap-pa+pap$.  Then $x\in \ann(\A,\B,f)$.   As $(\A,\B,f)$ is \wnd, $(I-p)a(I-p)=0$.
Since $0\leq a\leq a+p$, applying Lemma~\ref{positive}  to
$a+p$ gives
$(I-p) ap=pa(I-p)=0$.  Therefore,
\[ap=pa=pap=a.\]   As $\A$ is the linear span of its positive elements,
we conclude that $p=I_\A$.

\eqref{sameunit2}  Let $p=a^*a$. Then for every $b\in \B$,
$pf(b)=f(b)p=f(b)$ because $\B$ is a self-adjoint subspace of $\A$.  Therefore, $p^2-p\in \ann(\A,\B,f)=\{0\}$, so $p$ is a
projection in $\A$.   Consider $\C:=f(\B)+\bbC p$.  Then $\C$ is a unital
subalgebra of $\A$, and since $\ann(\A,\C, \subseteq) \subseteq
\ann(\A,\B,f)$, $(\A,\C,\subseteq)$ is \wnd.   By
part~\eqref{sameunit1}, $p=I_\A$.  Noting that $p-a\in \ann(\A,\B,f)$,
we obtain 
$I_\A=p=a$, as desired.

\eqref{sameunit3} Once again, we assume $\A\subseteq \bh$. Since
$(f(u_\lambda))$ is a bounded increasing net of positive semi-definite
operators, it converges in the strong operator topology.  Let
$p:=\sot\lim f(u_\lambda)$.  Then $p$ is a projection and for every $b\in \B$,
\begin{equation}\label{sameunit3.1} f(b)p=pf(b)=f(b).
\end{equation}
Since $f(u_\lambda)\leq a$ for each
$\lambda$, we obtain $p\leq a$.  Decomposing $a$ with respect to
$\H=(I-p)\H\oplus p\H$ we find
$a=\begin{bmatrix} p^\perp ap^\perp &p^\perp a p\\ pap^\perp& pap\end{bmatrix}$.  As $a\geq p$, we have
$pap\geq p\geq 0$.  Since 
$\norm{pap}\leq1$, we obtain  $p=pap$.  Then
\[0\leq a-p + p^\perp=\begin{bmatrix} p^\perp ap^\perp+p^\perp&p^\perp ap\\
    pap^\perp & pap-p\end{bmatrix} =\begin{bmatrix}p^\perp
    ap^\perp+p^\perp&p^\perp a p\\ pap^\perp&0\end{bmatrix}.\]
Lemma~\ref{positive} gives $p^\perp ap=0$.  Therefore
$a=\begin{bmatrix} p^\perp ap^\perp&0\\0&p\end{bmatrix}$.  It follows that for
$b\in \B$,
\[af(b)=f(b)a=f(b).\]  An application of part~\eqref{sameunit2}
completes the proof.

\eqref{sameunit4} There is nothing to do when $\A$ is unital, so
assume $\A$ has no unit.  Let $(a,\lambda)\in \ann(\tilde \A, \B,
u_\A\circ f)$.   It suffices to show $\lambda=0$, for once this is
established, it follows that $a\in \ann(\A,\B, f)$, whence  $a=0$.  To show
$\lambda=0$, we argue by contradiction.

Suppose $\lambda\neq 0$.  Then by scaling, we may assume
$\lambda=-1$.  Note that for $b\in \B$, $u_\A(f(b))=(f(b),0)$.  Thus
\[(0,0)=(a,-1)(f(b),0)=(af(b)-f(b),0)=(f(b)a-f(b),0)=(f(b),0)(a,-1).\]
So for all $b\in \B$, $af(b)=f(b)=f(b)a$.  Part~\eqref{sameunit2} now gives $a=I_\A$,
contrary to  assumption on $\A$.  Thus, $\lambda=0$,
completing the proof.

\eqref{sameunit5} Let $a=\tilde f(I_{\tilde\B})$.  Then
$\norm{a}\leq 1$ and $a\geq 0$.  For every $b\in \B$,
$a u_\A(f(b))=u_\A (f(b)) a=u_\A(f(b))$.  Applying
part~\eqref{sameunit4},  then part~\eqref{sameunit2}, yields
$a=I_{\tilde\A}$.  Therefore $(\tilde\A,\tilde\B, \tilde f)$ is a
unital inclusion and hence is \wnd.

\eqref{sameunit6} If $\B$ is unital, so is $\A$ by
part~\eqref{sameunit1}.  Therefore, $\sigma=\sigma\circ u_\B=u_\A\circ
f=f=\tilde f$, and all is well.

Now assume $\B$ is not unital.   We claim that
$\tilde f(I_{\tilde\B})=\sigma(I_{\tilde\B})$.   For $b\in \B$,
\[\sigma(I_{\tilde\B}) u_\A(f(b))=
  \sigma(I_{\tilde\B})\sigma(u_\B(b))=\sigma(u_\B(b))=u_\A(f(b));\]
likewise 
\[u_\A(f(b))\sigma(I_{\tilde\B})=u_\A(f(b)).\]
Similarly,
$\tilde f(I_{\tilde\B}) u_\A(f(b))=u_\A(f(b))\tilde f(I_{\tilde\B})=u_\A(f(b))$.  Therefore,
\[\tilde f(I_{\tilde\B})-\sigma(I_{\tilde\B})\in \ann(\tilde\A,\B, u_\A\circ f).\]
Part~\eqref{sameunit4} gives $\tilde
f(I_{\tilde\B})=\sigma(I_{\tilde\B})$, so the claim holds.

For $(b,\lambda)\in \tilde\B$,
\[\sigma(b,\lambda)=\sigma(u_\B(b))+\lambda \tilde
  f(I_{\tilde\B})=u_\A(f(b))+\lambda \tilde f(I_{\tilde\B}) =\tilde
  f(b,\lambda),\] so $\sigma=\tilde f$.  
\end{proof}

We now show that the poor behavior for intermediate inclusions
exhibited in Example~\ref{uincng} cannot occur when $(\B_3,\B_1)$ is
\wnd:  if $\B_2$ is intermediate to the \wnd\ inclusion
$(\B_3,\B_1)$, then $\tilde\B_2$ is intermediate to
$(\tilde\B_3,\tilde\B_1)$.
\begin{corollary}\label{interinc}
  For $1\leq i < j\leq 3$, let $(\B_j, \B_i, f_{ji})$ be inclusions
  such that $f_{31}=f_{32}\circ f_{21}$.  If $(\B_3,\B_1,f_{31})$ is
  \wnd, then $\tilde f_{31}=\tilde f_{32}\circ \tilde f_{21}$.
\end{corollary}
\begin{proof}
Since $(\B_3, \B_1, f_{31})$ is \wnd, both $(\B_3,\B_2, f_{32})$ and
$(\B_2,\B_1,f_{21})$ are \wnd.
Thus by
Lemma~\ref{sameunit}\eqref{sameunit5}, \[\tilde f_{32}(\tilde
  f_{21}(I_{\tilde\B_1}))=\tilde f_{32}(I_{\tilde\B_2})=I_{\tilde\B_3}=\tilde
  f_{31}(I_{\tilde\B_1}).\]

Given $x\in \tilde\B_1$, there exist $b_1\in \B_1$ and $\lambda\in
\bbC$ so that $x=u_{\B_1}(b_1)+\lambda I_{\tilde\B_1}$.  Using 
using~\eqref{ud1},
\begin{align*}\tilde f_{32}(\tilde f_{21}(x))
  &=\tilde f_{32}(\tilde f_{21}(u_{\B_1}(b_1)))+\lambda I_{\tilde\B_3}
  = \tilde f_{32}(u_{\B_2}(f_{21}(b_1)))+\lambda  I_{\tilde\B_3}\\
  &=u_{\B_3}( f_{32}((f_{21}(b_1))))+\lambda
    I_{\tilde\B_3}=u_{\B_3}(f_{31}(b_1))+\lambda I_{\tilde\B_3}\\
  &=\tilde f_{31}(x),
\end{align*} as desired.
\end{proof}

The following extends~\cite[Corollary~3.19]{PaulsenCoBoMaOpAl} from
unital inclusions to \wnd\ inclusions.
\begin{proposition}\label{wndbi}  Let $(\A_1,\B,f_1)$ and $(\A_2,\B,f_2)$ 
be \wnd\ inclusions and suppose \[\Phi: \A_1\rightarrow \A_2\] is a contractive
completely positive
map such that $\Phi\circ f_1=f_2$.   The following statements hold.
\begin{enumerate}
\item \label{wndbi1}  The standard extension,
  $\tilde\Phi:\tilde\A_1\rightarrow \tilde\A_2$, of $\Phi$  is a unital completely
  positive map. 
\item \label{wndbi2}  For
every $x\in \A_1$ and $h,k\in \B$,
\[\Phi(f_1(h)xf_1(k))=f_2(h)\Phi(x)f_2(k).\]
\end{enumerate}
\end{proposition}
\begin{proof}  For $i=1,2$,
parts~\eqref{sameunit4} and~\eqref{sameunit5} of
Lemma~\ref{sameunit} show that $(\tilde\A_i,\B, u_{\A_i}\circ f_i)$
are \wnd\ inclusions and 
$(\tilde\A_i,\tilde\B,\tilde f_i)$ are unital inclusions.

\eqref{wndbi1}  Observation~\ref{ccp} shows $\tilde\Phi$ is completely
positive.  To see $\tilde\Phi$ is a unital map, let $(h_\lambda)$ be an approximate unit
for $\B$.   We have $u_{\A_1}(f_1(h_\lambda))=\tilde f_1(u_\B(h_\lambda))\leq I_{\tilde\A_1}$, so
\[u_{\A_2}(f_2(h_\lambda))=u_{\A_2}(\Phi(f_1(h(u_\lambda))))=\tilde\Phi(u_{\A_1}(f_1(h_\lambda)))\leq \tilde\Phi(I_{\tilde\A_1}).\]
Applying
Proposition~\ref{sameunit}\eqref{sameunit3} to
$(\tilde\A_2, \B, u_{\A_2}\circ f_2)$, we obtain
$\tilde\Phi(I_{\tilde\A_1})=I_{\tilde\A_2}$, so $\tilde\Phi$ is a
unital map.

\eqref{wndbi2}
By
  \cite[Corollary~3.19]{PaulsenCoBoMaOpAl}, $\tilde \Phi$ 
is a $\tilde
  \B$-bimodule map, that is, for $\tilde h, \tilde k\in \tilde\B$ and
$\tilde x\in \tilde\A_1$,
\[\tilde \Phi(\tilde f_1(\tilde h)\tilde x\tilde f_1(\tilde k))=\tilde f_2(\tilde
h)\tilde \Phi(\tilde x)\tilde f_2(\tilde k).\]
Thus for $h, k\in \B$ and $x\in \A_1$, (and using \eqref{stext1})
\begin{align*}
(u_{\A_2}\circ\Phi)(f_1(h)xf_1(k)))&=(\tilde \Phi\circ u_{\A_1})(f_1(h)xf_1(k))\\ &=
\tilde \Phi(\tilde f_1(u_\B( h)) u_{\A_1}(x)\tilde f_1(u_\B (k))))\\
&=\tilde f_2(u_\B(h))\tilde \Phi(u_{\A_1}(x))\tilde f_2(u_{\B}(k)) \\
&=(u_{\A_2}\circ f_2)(h)\tilde \Phi(u_{\A_1}(x)) (u_{\A_2}\circ
                                                                        f_2)(k)\\
  &=(u_{\A_2}\circ f_2)(h) (u_{\A_2}\circ \Phi)(x) (u_{\A_2}\circ f_2)(k).
\end{align*} As $u_{\A_2}$ is a $*$-monomorphism, the proposition
follows.
\end{proof}

\subsection{Pseudo-Expectations}

We shall require pseudo-expectations, and we give a brief discussion
of them here.  Hamana~\cite{HamanaInEnC*Al} showed that given a
\textit{unital} \cstar-algebra $\A$ there is a \cstar-algebra $I(\A)$
and  a one-to-one unital
$*$-homomorphism $\iota :\A \rightarrow I(A)$
such that:
\begin{itemize}
  \item $I(\A)$ is an injective object in the category of operator systems and
    unital completely positive  maps; and
\item the only unital completely positive map $\phi: I(\A)\rightarrow
I(\A)$ satisfying $\phi\circ\iota=\iota$ is the identity mapping on
$I(\A)$. 
\end{itemize}
Hamana calls the pair $(I(\A),\iota)$ an \emph{injective
  envelope}\index{Injective envelope} of $\A$.  The injective envelope of $\A$ is monotone closed and has
the following uniqueness property: if $(\I_1,\iota_1)$ and
$(\I_2,\iota_2)$ are injective envelopes for $\A$, there exists a
unique $*$-isomorphism $\alpha:\I_1\rightarrow \I_2$ such that
$\alpha\circ \iota_1=\iota_2$.

\begin{remark}{Remark}\label{abcat}  Let $\fC$ be the category of
  unital, abelian \cstaralg s and unital $*$-homomorphisms.
  Hadwin and Paulsen note that an
   object $\A$ in $\fC$ is injective if and only if $\A$ is
   injective in the category of operator systems and completely
   positive unital maps, see
   \cite[Theorem~2.4]{HadwinPaulsenInPrAnTo}.  
 \end{remark}

When the \cstaralg\ $\A$ is non-unital,  an injective
envelope for $\A$ is defined to be an injective envelope $(I(\tilde\A),\iota)$ for
$\tilde \A$.
\begin{flexstate}{Remark}{}\label{injEnN} \rm To simplify notation, we will often say $(I(\A),\iota)$
is an injective envelope of $\A$ regardless of whether $\A$ has a unit:
when $\A$ is not unital, it is to be understood that $(I(\A),\iota)$ means
$(I(\tilde A),\iota)$, where $(I(\tilde\A),\iota)$ is an
injective envelope for $\tilde A$.   When there is a need for clarity,
we will sometimes write $(I(\tilde\A),\iota)$ for an injective envelope
of $\A$. 
\end{flexstate}

We need the following fact about injective envelopes; it is due to Hamana.
Let $(I(\A),\iota)$ be an injective envelope for the \cstaralg\ $\A$,
let $J\idealin \A$ be a closed ideal of $\A$,  let $(e_\lambda)$ be an
approximate unit for $\A$, and let $P=\sup_{I(\A)_{s.a.}}
\iota(e_\lambda)$, where $I(\A)_{s.a.}$ is the partially ordered set
of  self-adjoint elements of
$I(\A)$.   By~\cite[Lemma~1.1, parts (i) and
(iii)]{HamanaCeReMoCoCStAl}, $P$ is a central projection of $I(\A)$
and $(PI(\A), \iota|_J)$ is an injective envelope for $J$.  We will
call $P$ the \textit{support projection}\index{Support projection} for $J$.

The proof of the
next lemma follows from the uniqueness property of injective envelopes.
\begin{lemma}[c.f. {\cite[Proposition~1.11]{PittsStReInI}}]\label{AbInEn}
Suppose $\A$ is an abelian \cstaralg\ and let $(I(\A),\iota)$ be an
injective envelope for $\A$.  For $i=1,2$, let  $J_i$ be closed ideals
of $\A$ with support projections $Q_i$.   If $\phi: J_1\rightarrow J_2$ is
a $*$-isomorphism, then there is a unique $*$-isomorphism
$\overline\phi: Q_1I(\A)\rightarrow Q_2I(\A)$ such that
$\overline\phi\circ \iota|_{J_1}=\iota\circ \phi$.
\end{lemma}

The following is an interesting example of a \wnd\ inclusion.
\begin{lemma}\label{ndegIB}  Suppose
  $(I(\tilde\B),\iota)$ is an injective envelope for the \cstaralg\ $\B$.
  Then $(I(\tilde\B), \B, \iota\circ u_\B)$ is a \wnd\ inclusion. \rm 
\end{lemma}
\begin{proof}
  We suppress $\iota$ and $u_\B$, so that $\B\subseteq \tilde\B\subseteq
I(\tilde\B)$.     Let $(u_\lambda)$ be an approximate unit for $\B$, and
let $p:=\sup_{I(\tilde\B)} u_\lambda$, that is, $p$ is the least upper
bound of $\{u_\lambda\}$ taken in the partially ordered set $I(\tilde\B)_{sa}$.
As $\B$ is a hereditary subalgebra of
$\tilde\B$,~\cite[Lemma~1.1(i)]{HamanaCeReMoCoCStAl} shows $p$ is a projection
and $\fP:=pI(\tilde\B)p$ is an injective envelope for $\B$.

We claim that $p=I_{I(\tilde\B)}$.  As
$\B\subseteq \fP\subseteq I(\tilde\B)$,
$p=\sup_{I(\tilde\B)}u_\lambda\leq \sup_\fP u_\lambda\leq p$ so that
\begin{equation*}\label{mysup}
  \sup_{I(\tilde\B)}u_\lambda= \sup_\fP u_\lambda.
\end{equation*}
By uniqueness of injective envelopes, there is a $*$-isomorphism
$\theta: \fP\rightarrow I(\B)$ with
$\theta|_{\tilde\B}=\id|_{\tilde\B}$.  Then
\[I_{I(\tilde\B)}=\theta(p)=\theta(\sup_{\fP}
u_\lambda)=\sup_{I(\tilde\B)} \theta(u_\lambda)=\sup_{I(\tilde\B)}
u_\lambda =p.\]   Thus the claim holds.

For
$a\in \ann(I(\tilde\B),
\B)$,~\cite[Corollary~4.10]{HamanaReEmCStAlMoCoCStAl} shows
\[a^*a=a^*pa=a^*(\sup_{I(\tilde\B)}
u_\lambda)a=\sup_{I(\tilde\B)}a^*u_\lambda a =0.\]
\end{proof}

\begin{definition}[c.f.~\cite{PittsStReInI}] \label{PSdef} Let
  $(\A,\B,f)$ be an inclusion and let $(I(\B),\iota)$ be an injective
  envelope for $\B$.  A \emph{pseudo-expectation for $(\A,\B,f)$
    relative to $(I(\B),\iota)$}\index{Pseudo-expectation} is a
  contractive and completely positive linear map
  $E:\A\rightarrow I(\B)$ such that $E\circ f=\iota\circ u_\B$.
    
When $\B$ is identified with $f(\B)\subseteq \A$ and
$u_\B(\B)\subseteq \tilde\B$, we will write $E|_\B=\iota|_\B$ instead
of $E\circ f=\iota\circ u_\B$.   When these identifications are made, we will
simply say $E$ is a pseudo-expectation for $(\A,\B)$.
\end{definition}

\begin{remark}{Notation} \label{PsNot} For an inclusion $(\A,\B,f)$, $\PsExp(\A,\B,f)$
  will denote the collection of all pseudo-expectations for $(\A,\B,f)$ (relative to a
  fixed injective envelope $(I(\B),\iota)$ for $\B$).   As usual, when
               $f(\B)$ is identified with $\B$, we  write $\PsExp(\A,\B)$ instead
               of $\PsExp(\A,\B,f)$.
\end{remark}

In general, there are
  many pseudo-expectations.  In some cases however, there is a unique
  pseudo expectation, and this property plays an essential role in
  our study.
\begin{definition} \label{f!pseDef} We will say that 
  the inclusion $(\A,\B,f)$ has the \textit{unique
    pseudo-expectation property}\index{Pseudo-expectation!unique pseudo-expectation property} if $\PsExp(\A,\B,f)$ is a singleton
  set.  When $(\A,\B,f)$ has the unique pseudo-expectation property, and
  the (unique) pseudo-expectation is faithful, we say $(\A,\B,f)$ has
  the \textit{faithful unique pseudo-expectation
    property.}\index{Pseudo-expectation!faithful unique
    pseudo-expectation property}
\end{definition}

Note that by Remark~\ref{abcat}, if $\A$ is abelian, and $(\A,\B, f)$
has the unique pseudo-expectation property, the pseudo-expectation is
multiplicative.

That $\PsExp(\A,\B)$ is not empty
follows from the following general fact.
\begin{flexstate}{Lemma}{}\label{injext}   Suppose $\fI$ is an injective
  \cstaralg\ and $(\A,\B)$ is an inclusion.   If $\phi:\B\rightarrow \fI$
  is a contractive and completely positive  map, then $\phi$ extends to a
  contractive and completely positive map $\Phi: \A\rightarrow \fI$.
\end{flexstate}
\begin{proof}
  We may assume $\A\subseteq \bh$ satisfies $\overline{\A\H}=\H$ for
  some Hilbert space $\H$.  Thus $\tilde\A=\A+\bbC I_\H$.  Let
  $(u_\lambda)$ be an approximate unit for $\B$ and put
  $Q:=\sot\lim u_\lambda$.  Then $(Q\tilde\A Q, \B+\bbC Q)$ is a
  unital inclusion.  As $\tilde\B\simeq \B+\bbC Q$, we may
  apply~\cite[Proposition~2.2.1]{BrownOzawaC*AlFiDiAp} to obtain a
  unital and completely positive map
  $\tilde \phi: \B+\bbC Q\rightarrow \fI$ which extends $\phi$.  Injectivity
  of $\fI$ shows that we may extend $\tilde \phi$ to a unital, completely
  positive map $\Delta: Q\tilde \A Q\rightarrow I(\B)$.  Now observe
  that the map $\Phi: \A\rightarrow \fI$ given by $\Phi(a)=\Delta(QaQ)$
  is a contractive and completely positive map extending $\phi$.
\end{proof}

The next several results explore properties of pseudo-expectations.  We
begin with the relationship between pseudo-expectations for $(\A,\B,f)$
and pseudo-expectations for $(\tilde\A,\tilde\B,\tilde f)$. 
\begin{lemma}\label{prepse}  Let $(\A,\B,f)$ be an inclusion, and let
  $(I(\B),\iota)$ be an injective envelope for $\B$. The following
  statements hold.
\begin{enumerate}
\item\label{prepse0}  Let $(h_\lambda)$ be an approximate unit for
  $\B$.  Suppose $x\in \tilde \A$ satisfies $0\leq x$,
  $\norm{x}\leq 1$, and $u_\A(f(h_\lambda))\leq x$ for every
  $\lambda$.   If $E$ is a
pseudo-expectation for $(\A,\B,f)$, then $\tilde E(x)=I_{I(\B)}$.   
\item \label{prepse1}  If 
  $E$ is a pseudo-expectation for $(\A,\B, f)$, then \[\tilde
  E(I_{\tilde \A})=I_{I(\B)}=\tilde E(\tilde f(I_{\tilde\B}))\] and  $\tilde E$ is a
  pseudo-expectation for $(\tilde\A,\tilde\B,\tilde f)$.
\item \label{prepse2} Suppose $\Phi: \tilde\A\rightarrow I(\B)$ is a
  pseudo-expectation for $(\tilde\A,\tilde\B,\tilde f)$, and let
  $E:=\Phi\circ u_\A$.  Then:
  \begin{enumerate}
    \item $E$ is a
      pseudo-expectation for $(\A,\B,f)$,
    \item $\Phi(I_{\tilde \A})=\Phi(\tilde f(I_{\tilde\B}))=I_{I(\B)}$; and
      \item $\tilde E=\Phi$.
      \end{enumerate}
    \end{enumerate}
\end{lemma}
\begin{proof}
  \eqref{prepse0}   We have already noted in Observation~\ref{ccp}
  that $\tilde E$ is
  contractive and completely positive.  Let $z=\tilde E(x)$.  Then
  \[
   \iota(u_\B(h_\lambda)) =E(f(h_\lambda))=\tilde E(u_\A(f(h_\lambda)))\leq
    \tilde E(x)=z.\]   Since
  $(I(\B),\B,\iota\circ u_\B)$ is
  \wnd\
  (Lemma~\ref{ndegIB}), Proposition~\ref{sameunit}\eqref{sameunit3}
  shows $z=I_{I(\B)}$.
  
  \eqref{prepse1}   Taking $x=I_{\tilde\A}$ in part~\eqref{prepse0}
  shows  $\tilde E$ is a unital completely
  positive map and hence is contractive.   It remains to show
  $\tilde E\circ\tilde f=\iota$.  To do this, note that
  part~\eqref{prepse0} applied with $x=I_{\tilde\B}$ gives 
  $I_{I(\B)}=\tilde E(\tilde f(I_{\tilde\B}))$.   Since $\tilde\B
  =\bbC I_{\tilde\B} +\tilde f(u_\B(\B))$ and 
  $\tilde E\circ\tilde f\circ u_\B=E\circ f=\iota\circ u_\B$,  it follows that
  $\tilde E\circ \tilde f=\iota$, so part~\eqref{prepse1} holds.

\eqref{prepse2}   Since $\Phi$ is contractive and completely positive, so is $E$.
As $E\circ f=\Phi\circ u_\A\circ f=\Phi\circ \tilde f\circ u_\B=\iota\circ u_\B$,  $E$ is a
pseudo-expectation for $(\A,\B,f)$.   

Part~\eqref{prepse0} shows  $\Phi$ is a unital map.   Since
$\tilde\A=\bbC I_{\tilde\A} + u_\A(\A)$, it follows that $\Phi=\tilde
E$.  Part~\eqref{prepse1} now implies item (ii), which completes  the
proof of  part
\eqref{prepse2}.
\end{proof}
It is worth noting a bijection between $\PsExp(\tilde\A,\tilde\B)$ and
$\PsExp(\A,\B)$.
\begin{corollary}\label{PsNot1}
  The map,
  \begin{align*} \PsExp(\tilde\A,\tilde\B)\ni\Phi&\mapsto
                                                   \Phi|_\A\in \PsExp(\A,\B)\\
    \intertext{is a bijection with inverse}
    \PsExp(\A,\B)\ni E&\mapsto \tilde E\in \PsExp(\tilde\A,\tilde\B).
  \end{align*}
\end{corollary}
\begin{proof}
  Apply~Lemma~\ref{prepse}.
\end{proof}

Combining Proposition~\ref{wndbi} and Lemma~\ref{ndegIB} we obtain the
following.
\begin{flexstate}{Proposition}{} \label{PsExBi} Suppose $(\A,\B)$ is a \wnd\ inclusion and
  let $(I(\B),\iota)$ be an injective envelope for $\B$.
  If $E:\A\rightarrow I(\B)$ is a pseudo-expectation, then for every
  $x\in \A$ and $h, k\in \B$,
  \[E(hxk)=\iota(h)E(x)\iota(k).\]
\end{flexstate}

Our next goal to show that an inclusion $(\A,\B)$ has the faithful
unique pseudo-expectation property if and only if
$(\tilde\A,\tilde\B)$ does; this is Theorem~\ref{f!ps4rin} below.
The first step is to show that when $(\A,\B)$ has a faithful
pseudo-expectation, then $(\A,\B)$ is \wnd.

Our original statement of the following lemma had the stronger
hypothesis that $\tilde E$ is faithful.  We are grateful to \referee\ 
for suggesting the proof below, which shows the lemma holds with
the relaxed hypothesis that 
$E$ is faithful.  This 
yields an efficient route to the version of 
Theorem~\ref{f!ps4rin} presented here; the present version does not
require $\B$ to be abelian, unlike our previous version.

\begin{lemma} \label{pseprop}Let  $E$ be a pseudo-expectation for
  $(\A,\B)$.  If $E$ is faithful, then $(\A,\B)$ is
      \wnd. 
    \end{lemma}
    \begin{proof}   
     Observe that  $\ann(\A,\B)$ is a \cstar-subalgebra of $\A$.  Define a
map $\psi: \B\oplus\ann(\A,\B)\rightarrow \A$ by
$\psi(b\oplus x)=b+x$.   Since $\B\ann(\A,\B)=\ann(\A,\B)\B=\{0\}$,
$\psi$
is a $*$-homomorphism. If $\psi(b\oplus
x)=0$, 
\[x^*x=x^*(b+x)=0=b^*(b+x)=b^*b,\] whence $\psi$ is one-to-one.   By
definition, $\B\cap \ann(\A,\B)=\{0\}$.  Hence if
$z\in \B+\ann(\A,\B)$, there exists unique $b\in \B$ and
$x\in\ann(\A,\B)$ so that $z=b+x$.  Therefore $\ran(\psi)=\B+\ann(\A,\B)$,
whence $\psi$ is a
$*$-isomorphism of $\B\oplus\ann(\A,\B)$ onto the \cstaralg\ $\B+\ann(\A,\B)$.
Thus for $b+x\in \B+\ann(\A,\B)$,
\[\norm{b+x}=\max\{\norm{b},\norm{x}\}.\]

Now let $(I(\B),\iota)$ be an injective envelope for $\B$ and
let $(e_\lambda)$ be an approximate unit for $\B$.   Let   $x\in
\ann(\A,\B)$ satisfy $\norm{x}\leq 1$ and $x\geq 0$.  Since $0\leq
e_\lambda+x$ and $\norm{e_\lambda+x}\leq 1$, $e_\lambda +x \leq
I_{\tilde\A}$.   Then 
\[\iota(e_\lambda)+E(x)=E(e_\lambda+x)\leq I_{I(\B)}.\] Taking the
supremum in $I(\B)_{s.a.}$ gives $I_{I(\B)}+E(x)\leq I_{I(B)}$, so
$E(x)=0$.  By the faithfulness of $E$, $x=0$.  It follows that 
$\ann(\A,\B)=\{0\}$, as desired.
\end{proof}

\begin{proposition}\label{faith+ndg}   Let $(\A,\B)$ be an inclusion
  and suppose $E$ is a pseudo-expectation for $(\A,\B)$.  Then 
 $\tilde E$ is a faithful pseudo-expectation for
 $(\tilde \A,\tilde \B)$ if and only if  
 $E$ is faithful.
\end{proposition}
\begin{proof}
($\Rightarrow$) This is obvious.

($\Leftarrow$) Suppose $E$ is faithful.  By Lemma~\ref{pseprop}
$(\A,\B)$ is \wnd.  Since $\tilde E=E$ when $\A$ is unital, we may as
well assume $\A$ is not unital.

Suppose $0\leq p\in \tilde\A$ satisfies $\tilde E(p)=0$.   Write
$p=(x,\lambda)^2$, where $(x,\lambda)\in \tilde\A_{s.a.}$.  Using~Lemma~\ref{prepse}\eqref{prepse1}, 
\begin{align*} 0&=E(x^2)+2\lambda E(x)
   + \lambda^2 I_{I(\B)} \\
  &\geq E(x)^2 +2\lambda E(x)
+ \lambda^2 I_{I(\B)}\\ &= (E(x)+\lambda
I_{I(\B)})^2,
\end{align*}
whence \[E(x)=-\lambda I_{I(\B)}.\]    We claim that $\lambda=0$. 

We argue by contradiction.  Suppose $\lambda\neq 0$. The product of 
$(x,\lambda)$ with $-1/\lambda$ is still self-adjoint, so we  
may assume $\lambda=-1$.  With this assumption, $E(x)=I_{I(\B)}$.  
Then \begin{equation*} E(x^2)=I_{I(\B)}
\end{equation*}
because $0=\tilde E((x,-1)^2)$.

Given $b\in\B$,
Proposition~\ref{PsExBi} shows
\[E((xb-b)^*(xb-b))=\iota(b)^*E(x^2)\iota(b)-\iota(b)^*E(x)\iota(b)-\iota(b)^*E(x)\iota(b)
  +\iota(b^*b)=0.\]   Faithfulness of $E$ gives $xb=b$.   As $x$ is
self-adjoint, we find that for every $b\in \B$,
\[xb=b=bx.\]

Lemma~\ref{sameunit}\eqref{sameunit2} shows
that $x$ is the identity for $\A$, contradicting the assumption that
$\A$ is not unital.  Hence $\lambda=0$.

Thus $(x,\lambda)=(x,0)$ and using faithfulness of $E$ once again, we
find $x=0$.  Therefore $p=0$, so $\tilde E$ is faithful.
\end{proof}

\begin{theorem}\label{f!ps4rin}  Suppose $(\A,\B)$ is an
inclusion.
Then $(\A,\B)$ has
the faithful unique pseudo-expectation property if and only if
$(\tilde\A,\tilde\B)$ has the faithful unique pseudo-expectation
property.  When this occurs, $(\A,\B)$ is \wnd.
\end{theorem}
\begin{proof}
  Combine Corollary~\ref{PsNot1} with Proposition~\ref{faith+ndg}  to see that
  $(\A,\B)$ has the faithful unique 
  pseudo-expectation property if and only if $(\tilde\A,\tilde\B)$
  does.    The last statement follows from Lemma~\ref{pseprop}.
\end{proof}

\begin{corollary}\label{heriditary} Suppose $(\A,\B)$ is an
  inclusion and let $\E$ be a \cstaralg\ such that
  $\B\subseteq \E\subseteq \A$.    If $(\A,\B)$ 
has the faithful unique pseudo-expectation property, then $(\E,\B)$ has the faithful
  unique pseudo-expectation property.
\end{corollary}
\begin{proof}
By Theorem~\ref{f!ps4rin}, $(\A,\B)$ is \wnd.
Corollary~\ref{interinc} shows
$(\tilde\E, \tilde\B)$ is a
  unital inclusion with
  $\tilde\B\subseteq \tilde\E\subseteq \tilde\A$.  Thus, if
  $E:\A\rightarrow I(\B)$ is the pseudo-expectation for $(\A,\B)$,
  then $\tilde E$ is the faithful and unique pseudo-expectation for
  $(\tilde\A,\tilde\B)$.  Then
  by~\cite[Proposition~2.6]{PittsZarikianUnPsExC*In},
  $(\tilde\E,\tilde\B)$ has the faithful unique
  pseudo-expectation property and $(\tilde E)|_{\tilde\B}$ is the
  unique pseudo-expectation for $(\tilde\E,\tilde\B)$.
  Another application of Theorem~\ref{f!ps4rin} shows $(\E,\B)$ has
  the faithful unique pseudo-expectation property.
\end{proof}

For a unital inclusion $(\A,\B)$,~\cite[Corollary~3.14]{PittsZarikianUnPsExC*In}
shows that the faithful unique pseudo-expectation property implies
$\B^c$ is abelian.  
We now note that this useful structural fact  holds when $(\A,\B)$ is
  not assumed unital. 

  \begin{proposition}\label{f!pse->abelcom}
    If the inclusion $(\A,\B)$ has the faithful unique
  pseudo-expectation property, then $\B^c$ is abelian.
\end{proposition}
\begin{proof}
 We may suppose that for some Hilbert space $\H$,
  $\A\subseteq \bh$ and $\overline{\A\H}=\H$.  Then
  $\tilde\A=\A+\bbC I$.  Let $E:\A\rightarrow I(\B)$ be the faithful
  unique pseudo-expectation, and put \[\Phi:=\tilde E.\] Then $\Phi$
  is the unique pseudo-expectation for $(\tilde\A,\tilde\B)$ and by
  Proposition~\ref{faith+ndg}, $\Phi$ is faithful.  Applying
  \cite[Corollary~3.14]{PittsZarikianUnPsExC*In}, we see that
  $(\tilde\B)^c$ is abelian.  Since $\B^c\subseteq (\tilde\B)^c$,
  $\B^c$ is abelian.
\end{proof}

We now  extend~\cite[Proposition~5.5(b)]{PittsStReInII}
from the unital setting to include the non-unital case.  The utility
of the following 
comes from the fact that it is sometimes easier to establish the
faithful unique pseudo-expectation property than to show a
subalgebra is a MASA.
\begin{proposition}\label{Car+upse}   Let $(\C,\D)$ be a regular
  inclusion with $\D$ abelian. 
  The following statements are equivalent.
  \begin{enumerate}
  \item\label{Car+upse1} $(\C,\D)$ is a Cartan inclusion.
    \item \label{Car+upse2}  $(\C,\D)$ has
  the faithful unique pseudo-expectation property and there is a
  conditional expectation $E:\C\rightarrow \D$.
\end{enumerate}
\end{proposition}
\begin{proof}
 First notice that if $(\C,\D)$ satisfies either condition
\eqref{Car+upse1} or \eqref{Car+upse2}, then $(\C,\D)$ is \wnd.
  Indeed, if \eqref{Car+upse1} holds, then
  Lemma~\ref{->ndeg}\eqref{->ndeg2} implies $(\C,\D)$ is \wnd;
    when \eqref{Car+upse2} holds, apply Lemma~\ref{pseprop}. Thus in both cases,
      $(\tilde\C,\tilde\D)$ is a unital inclusion
      (Lemma~\ref{sameunit}\eqref{sameunit5}). 

      Our arguments establishing
      \eqref{Car+upse1}$\Leftrightarrow$\eqref{Car+upse2} differ
        depending on whether $\C$ is unital or non-unital, so we
        consider those cases separately.

 \textsc{Case 1}: Assume $\C$ is unital. \par
\eqref{Car+upse1}$\Rightarrow$\eqref{Car+upse2}   
      By hypothesis, $(\C,\D)$ is a MASA inclusion, hence $(\C,\D)$ is
      a unital inclusion.   Then
      \cite[Proposition~5.5(b)]{PittsStReInII} gives
      \eqref{Car+upse2}.

 \eqref{Car+upse2}$\Rightarrow$\eqref{Car+upse1}  For any $d\in\D$, we
   have
   \[E(I_\C)d=E(d)=d=dE(I_\C),\] so $E(I_\C)$
is the unit for $\D$.  As $(\C,\D)$ is \wnd, it 
   is a unital inclusion.   Now 
   \cite[Proposition~5.5(b)]{PittsStReInII} shows $(\C,\D)$ is a
   Cartan inclusion.

   \textsc{Case 2}: Assume $\C$ is not unital.   Then $\D$ is not
   unital by Lemma~\ref{sameunit}.\par

   \eqref{Car+upse1}$\Rightarrow$\eqref{Car+upse2}
     By~\cite[Proposition~3.2]{PittsNoApUnInC*Al},
     $(\tilde\C,\tilde\D)$ is a Cartan inclusion.  Case 1 shows: i)
     $(\C,\D)$ has the faithful unique pseudo-expectation property
     (see Corollary~\ref{PsNot1}); and ii) there exists a conditional
     expectation $\tilde E:\tilde\C\rightarrow\tilde\D$.  We wish to
     show $\tilde E|_\C$ is a conditional expectation of $\C$ onto
     $\D$.  Let
     $(u_\lambda)$ be an approximate unit for $\D$.  Then
     $(u_\lambda)$ is an approximate unit for $\C$ by
     \cite[Theorem~2.5]{PittsNoApUnInC*Al}.  So for $x\in \C$, the
     fact that $\D$ is an ideal in $\tilde\D$ gives
\[\tilde E(x)=\lim_\lambda \tilde E(xu_\lambda)=\lim_\lambda \tilde E(x)
  u_\lambda\in \D.\]  Thus $\tilde E|_\C$ is a conditional expectation of
$\C$ onto $\D$.  This establishes~\eqref{Car+upse2}.

  \eqref{Car+upse2}$\Rightarrow$\eqref{Car+upse1} By
    Proposition~\ref{faith+ndg} and Corollary~\ref{PsNot1},
      $(\tilde\C,\tilde\D)$ has the faithful unique pseudo-expectation
      property.  Noting that $\tilde E$ is a conditional
      expectation of $\tilde\C$ onto $\tilde\D$, we conclude from Case
      1 that $(\tilde\C,\tilde\D)$ is a Cartan inclusion.  Applying
      \cite[Proposition~3.2]{PittsNoApUnInC*Al} again, we find
      $(\C,\D)$ is a Cartan inclusion.  This completes the proof.
\end{proof}

\subsection{The Ideal Intersection Property}

The purpose of this section is to show that the inclusions $(\A,\B,f)$ and
$(\tilde\A, \tilde\B,\tilde f)$ both have the \iip\ or both do not,
and also to explore some
relationships between the \iip\ and the faithful unique
pseudo-expectation property.
\begin{lemma}\label{iipiffuiip}
  Let $(\A,\B,f)$ be an inclusion.
  Then $(\A,\B,f)$ has the \iip\ if and only if
  $(\tilde\A,\tilde\B,\tilde f)$ has the \iip.  
\end{lemma}
\begin{proof}
  Suppose $(\A,\B,f)$ has the \iip\ and let $J\idealin \tilde\A$
satisfy $J\cap \tilde f(\tilde\B)=\{0\}$. As 
\[(J\cap u_\A(\A))\cap u_\A(f(\B)) =(J\cap u_\A(\A))\cap \tilde
  f(u_\B(\B))\subseteq J\cap\tilde f(\tilde\B)=\{0\},\] the \iip\ for
$(\A,\B,f)$ gives $J\cap u_\A(\A)=\{0\}$.  Since $u_\A(\A)$ is an
essential ideal in $\tilde\A$, we obtain $J=0$. Thus $(\tilde\A,\tilde\B,\tilde f)$ has
the \iip.

  Now suppose $(\tilde\A,\tilde\B,\tilde f)$ has the \iip\ and let $J\idealin\A$
  satisfy  $J\cap f(\B)=\{0\}$.   Then $u_\A(J)\idealin \tilde\A$, and we claim
  \begin{equation}\label{iipiffuiip1}
    u_\A(J)\cap \tilde f(\tilde\B)=\{0\}.
  \end{equation}
  When $\B$ is unital, $\tilde\B=\B$ and $\tilde f=u_\A\circ f$, so
  \eqref{iipiffuiip1} holds because $u_\A$ is one-to-one.

Suppose then that $\B$ is not unital.  Recall that $\{0\}$
is a unital algebra, so in particular,  $\B\neq \{0\}$.
If $x\in u_\A(J)\cap \tilde f(\tilde\B)$, then
$x$ has the form 
\[x=\tilde f(b,\lambda)= u_\A(f(b))+\lambda
I_{\tilde\A}\] 
for some $b\in \B$ and
$\lambda\in \bbC$.  Since $J\cap f(\B)=\{0\}$,~\eqref{iipiffuiip1} will
follow once we show $\lambda=0$.

Arguing by contradiction, suppose $\lambda\neq 0$.  By scaling, we may
assume $\lambda=-1$, so that $x=u_\A(f(b))-I_{\tilde\A}$, for some
$b\in \B$.  For any $h\in \B$, we have $xu_\A(f(h))\in u_\A(J)\cap u_\A(f(\B))=\{0\}$,
so $xu_\A(f(h))=0$.    This gives $bh=h$.  Similarly, for $h\in\B$,
$hb=h$, so $b$ is the unit for $\B$, contrary to hypothesis.  Hence
$\lambda=0$, 
completing the proof of~\eqref{iipiffuiip1}.
  
  Since $(\tilde\A,\tilde\B,\tilde f)$ has the \iip, ~\eqref{iipiffuiip1} gives
  $J=\{0\}$.  Thus $(\A,\B,f)$ has the \iip.
\end{proof}

Our next result concerns 
the \iip\ for intermediate inclusions in the abelian case.
It is the same
as~\cite[Lemma~5.4]{PittsStReInII}  except the hypothesis that the
\cstaralg s involved have a common unit is dropped and we explicitly
consider the inclusion mappings.
\begin{lemma}\label{interiip}  
For $1\leq i < j\leq 3$, Let $(\D_j,\D_i,f_{ji})$ be inclusions such
that $f_{31}=f_{32}\circ f_{21}$ and $\D_3$ is abelian.
          Then
$(\D_3,\D_1,f_{31})$ has the \iip\ if and only if both $(\D_3,\D_2,f_{32})$ and
$(\D_2,\D_1,f_{21})$ have the \iip.
\end{lemma}
\begin{proof}
  The implication $(\Leftarrow)$ is 
left to the reader.

$(\Rightarrow)$ Suppose $(\D_3,\D_1,f_{31})$ has the \iip.  
Lemma~\ref{->ndeg}\eqref{->ndeg1} shows
$(\D_3,\D_1,f_{31})$ is \wnd.  Therefore, both
$(\D_2,\D_1,f_{21})$ and $(\D_3,\D_2,f_{32})$ are \wnd.
Then
 for $1\leq i<j\leq 3$, each of  $(\tilde\D_j, \tilde\D_i, \tilde
f_{ji})$ is a unital inclusion, hence  $\tilde f_{31}=\tilde f_{32}\circ
\tilde f_{21}$ (Corollary~\ref{interinc}).  By
Lemma~\ref{sameunit}\eqref{sameunit1}, these
inclusions have a common unit in the sense that  
\[I_{\tilde\D_3}=\tilde f_{32}(I_{\tilde\D_2})\dstext{and}
I_{\tilde\D_2}=\tilde f_{21}(I_{\tilde\D_1}).\]  By \cite[Lemma~5.4]{PittsStReInII}, 
$(\tilde \D_2,\tilde\D_1,\tilde f_{21})$ and
$(\tilde\D_3,\tilde\D_2,\tilde f_{32})$ have the \iip.  An application
of 
Lemma~\ref{iipiffuiip} completes the proof.
\end{proof}

Our final result of this section extends a portion
of~\cite[Corollary~3.22]{PittsZarikianUnPsExC*In} to include 
 non-unital settings.
\begin{proposition}\label{PZC22}  Suppose $(\D_2,\D_1)$ is an inclusion
  with $\D_2$ abelian and let $(I(\D_1),\iota_1)$ be an injective envelope
  for $\D_1$.  Then $(\D_2,\D_1)$ has the \iip\ if and only
  if $(\D_2,\D_1)$ has the faithful unique pseudo-expectation
  property.

  When this occurs, the unique pseudo-expectation $\iota_2:
  \D_2\rightarrow I(\D_1)$ is a $*$-monomorphism and
  $(I(\D_1),\iota_2)$ is an injective envelope for $\D_2$.
\end{proposition}
\begin{proof}
  ($\Rightarrow$) Suppose $(\D_2,\D_1)$ has the \iip.
 An application of  Lemma~\ref{->ndeg}\eqref{->ndeg3} shows $(\D_2,\D_1)$ is \wnd. 
By Proposition~\ref{sameunit}\eqref{sameunit5} and 
Lemma~\ref{iipiffuiip}, $(\tilde\D_2,\tilde\D_1)$ is a unital
inclusion having the \iip, so 
\cite[Corollary~3.22]{PittsZarikianUnPsExC*In} shows 
 $(\tilde\D_2,\tilde\D_1)$ has the faithful unique pseudo-expectation
 property.  Corollary~\ref{PsNot1} shows $(\D_2,\D_1)$ has the faithful
 unique pseudo-expectation property.

($\Leftarrow$)   Suppose $(\D_2,\D_1)$ has the faithful unique
pseudo-expectation property. By Lemma~\ref{pseprop},
$(\D_2,\D_1)$ is \wnd, whence 
$(\tilde\D_2,\tilde\D_1)$ is a unital inclusion.
Apply \cite[Corollary~3.22]{PittsZarikianUnPsExC*In} to see
 $(\tilde\D_2,\tilde\D_1)$ has the \iip. Lemma~\ref{iipiffuiip} shows
 $(\D_2,\D_1)$ has the \iip.

 Turning to the last statement, suppose $(\D_2,\D_1)$ has the faithful
 unique pseudo-expectation property with pseudo-expectation $\iota_2$.
 Applying Lemma~\ref{pseprop} and Lemma~\ref{sameunit}\eqref{sameunit5},
 we see that  $(\tilde\D_2,\tilde\D_1)$ is a unital inclusion.
 Corollary~\ref{PsNot1} shows
 $\tilde\iota_2:\tilde\D_2\rightarrow I(\D_1)$ is the unique
 pseudo-expectation for $(\tilde\D_2,\tilde\D_1)$ and by
 Theorem~\ref{f!ps4rin}, $\tilde\iota_2$ is faithful.
 By~\cite[Corollary~3.22]{PittsZarikianUnPsExC*In},
 $\tilde\iota_2$ is a $*$-monomorphism, hence so is $\iota_2$.

 Finally, since $\iota_1(\D_1)\subseteq \iota_2(\D_2)\subseteq
 I(\D_1)$, the minimality of injective envelopes shows
 $(I(\D_1),\iota_2)$ is an injective envelope for $\D_2$. 
 \end{proof}



\mysec[Cartan Envelopes]{Cartan Envelopes}\label{mainr}
\numberwithin{equation}{section}

Up to this point, we have mostly considered general inclusions.  We now restrict
attention to inclusions where the subalgebra is abelian. 
\begin{remark}{Standing Assumption}\label{SA1}
  \it Unless explicitly stated otherwise, for the remainder of this work, whenever 
  $(\C,\D)$ is an inclusion,  we shall always assume $\D$ is  \sc abelian.\rm\ 
\end{remark}

In~\cite{PittsStReInII}, we defined the notion of a Cartan envelope
for unital and regular inclusions and characterized when such
inclusions have a Cartan envelope, see
\cite[Theorem~5.2]{PittsStReInII}.  The purpose of this \LCsection\  is to
extend the characterization of regular inclusions having a Cartan
envelope from unital regular inclusions to all regular inclusions.
This is accomplished in Theorem~\ref{!pschar}; it is among our main
results.  While the statement of Theorem~\ref{!pschar} parallels that
of \cite[Theorem~5.2]{PittsStReInII}, discarding the hypothesis that
$(\C,\D)$ is unital presents challenges which cannot be overcome by
simply adjoining a unit.

 We begin by recording a few useful facts about
               inclusions satisfying Assumption~\ref{SA1}.
\begin{lemma}[{\cite[Corollary~2.2]{PittsNoApUnInC*Al}}]\label{Tv2.2}
  Let $(\C,\D)$ be an inclusion and fix a normalizer $v\in \N(\C,\D)$.  Then
  $\overline{vv^*\D}$ and $\overline{v^*v\D}$ are ideals in $\D$ and
  the map $vv^*d\mapsto v^*dv$ uniquely extends to a $*$-isomorphism
   $\theta_v:\overline{vv^*\D}\rightarrow \overline{v^*v\D}$; moreover,
   for every $h\in \overline{vv^*\D}$, $v\theta_v(h)=hv$.
\end{lemma}

               Our proofs of the statements in the next lemma
               depend on results from \cite{PittsNoApUnInC*Al}.

\begin{lemma}\label{relcom}  Suppose $(\C,\D)$ is an inclusion 
 and let
  $\D^c$ be the relative commutant of $\D$ in $\C$.  The
  following statements hold.
  \begin{enumerate}
  \item \label{relcom1} The identity mapping on $\C$ is a regular map
    from $(\C,\D)$ into $(\C,\D^c)$, that is,
    \[\N(\C,\D)\subseteq \N(\C,\D^c).\]
\item \label{relcom2} Suppose $(\C_1,\D_1)$ is an inclusion and
  $\alpha:(\C,\D)\rightarrow (\C_1,\D_1)$ is a regular
  $*$-monomorphism.  (While regularity gives $\alpha(\D)\subseteq
  \N(\C_1,\D_1)$, we do not assume $\alpha(\D)\subseteq \D_1$, that
  is, we do not assume $(\C_1,\D_1\ms\alpha)$ is an
  expansion of $(\C,\D)$.)
  Then
  \[\alpha(\D)\subseteq \D_1^c.\]  In addition, if $(\C_1,\D_1)$ has
  the AUP, then $\alpha(\D)\subseteq \D_1$. 

\item \label{relcom3} If $(\C,\D)$ has the AUP, then $\N(\C,\D)\subseteq
  \N(\tilde\C,\tilde\D)$.
\item \label{relcom4}   Suppose $(\C,\D)$ is a regular
and \wnd\ inclusion with $\C$ not unital.   If $\N(\C,\D)\subseteq
  \N(\tilde\C,\tilde\D)$, then $(\C,\D)$ has the AUP. 
\item \label{relcom5} Suppose $(\C,\D)$ is regular and has the AUP.   Then 
  $(\tilde\C,\tilde\D)$ is regular.
\end{enumerate}
\end{lemma}
\begin{proof}
  \ref{relcom1}) Let $v\in \N(\C,\D)$.
Let $d\in \D$ and
  $x\in \D^c$.  Using Lemma~\ref{Tv2.2}, we find that for every $h\in\overline{vv^*\D}$, 
  \[(vxv^*)dh=vx\theta_v(dh)v^*=v\theta_v(dh)xv^*=dh(vxv^*).\] Taking
  $h=u_\lambda$ where $(u_\lambda)$ is an approximate unit for
  $\overline{vv^*\D}$, and using the fact that $\norm{v^*u_\lambda -
  v^*}\rightarrow 0$, we obtain $vxv^*\in \D^c$.  Likewise
  $v^*xv\in \D^c$, so \eqref{relcom1} holds.

    \ref{relcom2}) Let $0\leq d\in \D$.  Since $d^{1/2}\in\D\subseteq
    \N(\C,\D)$,  we have $\alpha(d^{1/2})\in \N(\C_1,\D_1)$.  For
    every $h\in \D_1$,~\cite[Proposition~2.1]{PittsNoApUnInC*Al} gives,
    \begin{equation}\label{relcom2.1}\alpha(d)h
      =\alpha(d^{1/2})^*\alpha(d^{1/2})h=h\alpha(d) \in
      \D_1.
    \end{equation}
    Thus $\alpha(d)\in \D_1^c$.   Since $\D$
  is the linear span of its positive elements, 
  $\alpha(\D)\subseteq \D_1^c$. 

  Let $(u_\lambda)\subseteq \D_1$ be an approximate unit for $\C_1$.
  Replacing $h$ in \eqref{relcom2.1} with $u_\lambda$ and taking the limit along
  $\lambda$ gives $\alpha(d)\in\D_1$.

  \ref{relcom3}) Suppose $(\C,\D)$ has the AUP.  By Lemma~\ref{sameunit}\eqref{sameunit5},
  $(\tilde\C,\tilde\D)$ is a unital inclusion.  The conclusion is
  obvious when $\C$ is unital.  So  assume $\C$ is not unital and
  let $(u_\lambda)$ be a net in $\D$ which is an approximate unit for
  $\C$.  For $v\in \N(\C,\D)$,
  ~\cite[Proposition~2.1]{PittsNoApUnInC*Al} shows $v^*vu_\lambda$ and
  $vv^* u_\lambda$ belong to $\D$.  Taking the limit shows $v^*v$ and
  $vv^*$ are elements of $\D$.  Thus, if $(d,\lambda)\in \tilde\D$,
\begin{equation}\label{relcom6} (v,0)^*(d,\lambda)(v,0)=(v^*dv+
    \lambda v^*v,0)\dstext{and} (v,0)(d,\lambda)(v,0)^*=(vdv^*+
    \lambda vv^*,0).
  \end{equation}  
  This gives $\N(\C,\D)\subseteq \N(\tilde\C,\tilde\D)$.  

\ref{relcom4})   Suppose $(\C,\D)$ is regular, \wnd, and $\N(\C,\D)\subseteq
\N(\tilde\C,\tilde\D)$.  Since $\C$ is not unital, neither is $\D$, so
that $(\tilde\C,\tilde\D)=(\C^\dag,\D^\dag)$. 
The calculation in \eqref{relcom6} with $d=0$ and 
$\lambda=1$ shows that for every $v\in
  \N(\C,\D)$, $v^*v\in \D$.
  Then~\cite[Observation~1.3(2)]{PittsNoApUnInC*Al} shows $(\C,\D)$
  has the AUP.

  \ref{relcom5})  Suppose $(\C,\D)$ is regular and  has the AUP.
  There is nothing to do if $\C$ is unital, for then
  $(\tilde\C,\tilde\D)=(\C,\D)$.   Assume then, that $\C$ is not
  unital.  Regularity of $(\C,\D)$ gives \[\tilde\C=\{(0,\lambda):
  \lambda\in \bbC\}
  +\overline{\spn}\{(v,0): v\in\N(\C,\D)\}.\] By part
\eqref{relcom3}, $(\tilde\C,\tilde\D)$
  is regular.
\end{proof}

While the following corollary is immediate from parts~\eqref{relcom3}
and~\eqref{relcom4} of Lemma~\ref{relcom}, it is worth making
explicit.  In particular, $u_\C$ is not regular when $\D$ is an
essential and proper ideal in the abelian \cstaralg\ $\C$; this
also shows that~\eqref{uCnr3}$\not\Rightarrow$\eqref{uCnr2} in Corollary~\ref{uCnr}.     

\begin{corollary}\label{uCnr}  Let $(\C,\D)$ be a regular and \wnd\
  inclusion with $\C$ not unital.  Consider the following statements.
\begin{enumerate}
\item\label{uCnr1} $(\C,\D)$ has
  the AUP;
\item\label{uCnr2} the inclusion mapping $u_\C:(\C,\D)\rightarrow (\tilde\C,\tilde\D)$
  is a regular $*$-monomorphism;
\item\label{uCnr3} $(\tilde\C,\tilde\D)$ is a regular inclusion.
\end{enumerate}  The following implications hold:
 \eqref{uCnr1}$\Leftrightarrow$\eqref{uCnr2}$\Rightarrow$\eqref{uCnr3}.
\end{corollary}
\begin{proof}
\eqref{uCnr1}$\Leftrightarrow$\eqref{uCnr2}  Use 
\ref{relcom}\eqref{relcom3} and~\ref{relcom}\eqref{relcom4}.

\eqref{uCnr1}$\Rightarrow$\eqref{uCnr3} Apply \ref{relcom}\eqref{relcom5}. 
\end{proof}

The next corollary of Lemma~\ref{relcom} is sometimes useful.
\begin{corollary}\label{cor:rmasa}  Suppose for $i=1,2$, $(\C_i,\D_i)$ are regular
  MASA inclusions.  If $\alpha:\C_1\rightarrow \C_2$ is a $*$-isomorphism satisfying
  $\alpha(\N(\C_1,\D_1))\subseteq\N(\C_2,\D_2)$, then
  $\alpha(\N(\C_1,\D_1))=\N(\C_2,\D_2)$ and $\alpha(\D_1)=\D_2$. 
\end{corollary}
\begin{proof}
By~\cite[Theorem~2.5]{PittsNoApUnInC*Al}, every regular MASA inclusion
has the AUP. 
Lemma~\ref{relcom}\eqref{relcom2} gives $\alpha(\D_1)\subseteq
\D_2$.   As $\alpha^{-1}(\D_2)\subseteq \D_1^c=\D_1$, we
obtain $\alpha(\D_1)=\D_2$.   It follows that  $w\in \N(\C_2,\D_2) \iff
\alpha^{-1}(w)\in \N(\C_1,\D_1)$, hence 
$\alpha(\N(\C_1,\D_1))=\N(\C_2,\D_2)$.
\end{proof}

\begin{example}\label{failinc}
  Lemma~\ref{relcom}(b)  shows $\alpha(\D)\subseteq \B$ whenever
$(\A,\B)$ has the AUP and $\alpha$ is a regular map.  However, if
$(\A,\B)$ does not have the AUP, this can fail. For an elementary
example of this behavior, let
\[\C=\A=C_0(\bbR), \quad \D=\B=\{h\in C_0(\bbR): h(x)=0 \text{ for
    all } x>1\},\] and define $\alpha:\C\rightarrow \A$ by
\[\alpha(h)(x)=h(x/2).\]   Since 
$\N(\C,\D)=\N(\A,\B)=C_0(\bbR)$,  $\alpha$ is a regular map, but
$\alpha(\D)$ is not contained in $\B$.
\end{example}

The following extends parts of~\cite[Theorem~3.5]{PittsStReInI} to
include inclusions which may not be unital.
\begin{theorem}\label{regmasaps}  Suppose $(\C,\D)$ is a regular MASA
inclusion.  Then $(\C,\D)$ has the unique pseudo-expectation
property.  Furthermore, $(\C,\D)$ has the faithful unique
pseudo-expectation property if and only if $(\C,\D)$ has the \iip.
\end{theorem}
\begin{proof}
  Since $(\C,\D)$ is a regular MASA inclusion, it has the AUP
  (\cite[Theorem~2.5]{PittsNoApUnInC*Al}), whence
  $(\tilde\C,\tilde\D)$ is a unital inclusion
  (Lemma~\ref{->ndeg}\eqref{->ndeg3}), and regular by
  Lemma~\ref{relcom}\eqref{relcom5}.  

We claim that $(\tilde\C,\tilde\D)$ is a MASA inclusion.  If $\D$ is
unital, then so is $\C$ by Lemma~\ref{sameunit}\eqref{sameunit1};  thus in
this case $(\tilde\C,\tilde\D)=(\C,\D)$ and all is well.
Suppose then that $\D$ is not unital.  Since $\D$ is a MASA,  $\C$
cannot be unital.  Thus $(\tilde\C,\tilde\D)=(\C^\dag,\D^\dag)$ and a
routine argument shows $\tilde\D$ is a MASA  in $\tilde\C$.  

We have established that $(\tilde\C,\tilde\D)$ is a regular, unital, 
  MASA inclusion.  By ~\cite[Theorem~3.5]{PittsStReInI}
  $(\tilde\C,\tilde\D)$ has the unique
  pseudo-expectation property, so  Corollary~\ref{PsNot1} shows $(\C,\D)$
  also has the
  unique pseudo-expectation property.

  Turning to the second statement, let
  $E:\C\rightarrow I(\D)$ be the pseudo-expectation.  Then $\tilde E$
  is the pseudo-expectation for $(\tilde \C,\tilde\D)$. 

Suppose $(\C,\D)$ has the \iip. Lemma~\ref{iipiffuiip}
implies $(\tilde \C,\tilde\D)$ has the \iip.
By~\cite[Theorem~3.15]{PittsStReInI}, the left kernel $\L:=\{x\in
\tilde\C: \tilde E(x^*x)=0\}$ is an ideal of $\tilde\C$ having trivial
intersection with $\tilde \D$.  Therefore,  $\tilde E$ is faithful,
whence $E$ is faithful.

Now suppose $E$ is faithful.  By Proposition~\ref{faith+ndg},
$\tilde E$ is faithful. If $J\idealin \tilde\C$ has trivial
intersection with $\tilde\D$, \cite[Theorem~3.15]{PittsStReInI} gives
$J\subseteq\L=\{0\}$.  So $(\tilde\C,\tilde\D)$ has the \iip.
By Lemma~\ref{iipiffuiip},  $(\C,\D)$ has the \iip, completing the proof.
\end{proof}

\begin{definition}\label{defextinc}
  Suppose $(\C,\D)$ is an inclusion and $(\A,\B\ms \alpha)$ is an
  expansion of $(\C,\D)$. 
  \begin{enumerate}
  \item \label{defextinc2}
We call $(\A,\B\ms \alpha)$ an \textit{essential
  expansion}\index{Expansion!essential} of $(\C,\D)$ if $(\B, \D, \alpha|_{\D})$ is an inclusion
having the \iip.
  \item\label{defextinc1}  If the $*$-monomorphism $\alpha$ is 
a regular map, we will say 
  $(\A,\B\ms \alpha)$ is a \textit{regular expansion}\index{Expansion!regular} of $(\C,\D)$.
\item\label{defextinc4}
If $(\A,\B)$ is a Cartan pair,  we say
  $(\A,\B\ms \alpha)$ is a \textit{Cartan expansion}\index{Expansion!Cartan} of
  $(\C,\D)$.
 \end{enumerate}
\end{definition}

\begin{remark}{Remarks}\label{R:dextinc}
    \begin{enumerate}
      \item \label{R:dextinc1} It may seem odd that in Definition~\ref{defextinc}\eqref{defextinc2},
we only require that  $(\B, \D, \alpha|_{\D})$ has the \iip\ instead
of also placing that requirement on $(\A,\B,\alpha)$.  In the context of
most interest to us, that is, when $(\C,\D)$ and $(\A,\B)$  both have the
faithful unique pseudo-expectation property, we will see in Observation~\ref{Cext}
that this is automatic.

\item \label{R:dextinc2} In some cases, regular maps produce
  expansions. For example,
Lemma~\ref{relcom}\eqref{relcom2} shows that  
when $\alpha:(\C_1,\D_1)\rightarrow (\C_2,\D_2)$ is a regular
$*$-monomorphism, and $(\C_2,\D_2)$ is a Cartan inclusion, then
$(\C_2,\D_2\ms \alpha)$ is automatically a regular expansion of
$(\C_1,\D_1)$.  
\end{enumerate}
\end{remark}

The next two lemmas give some useful properties of essential expansions.
\begin{lemma}\label{cewnd}  Suppose $(\A,\B\ms \alpha)$ is an
essential expansion for the inclusion $(\C,\D)$ and there is a
faithful conditional expectation $\Delta:\A\rightarrow \B$.    
Then $(\A,\D, \alpha|_\D)$
and $(\A,\C,\alpha)$ are \wnd\ inclusions.  In particular, if $\C$ is
unital, then $\A$ is unital and $\alpha(I_\C)=I_\A$. 
\end{lemma}
\begin{proof}  By hypothesis, 
$(\B,\D,\alpha|_\D)$ has the \iip, so it is \wnd\ by
Lemma~\ref{->ndeg}\eqref{->ndeg1}. 
Suppose $0\leq a\in \A$ belongs to $\ann(\A,\D,\alpha|_\D)$.   Then for every $d\in \D$,
\[0=\Delta(\alpha(d)a)=\alpha(d)\Delta(a)=\Delta(a)\alpha(d)=\Delta(a\alpha(d)).\]
Thus $\Delta(a)\in \ann(\B,\D,\alpha|_\D)$. 
Faithfulness of $\Delta$ gives $a=0$.   As $\ann(\A,\D,\alpha)$ is the
span of its postive elements, we obtain $\ann(\A,\D,\alpha)=\{0\}$, so
$(\A,\D,\alpha|_\D)$ is \wnd.   That $(\A,\C,\alpha)$ is \wnd\ follows
from the fact that
$\ann(\A,\C,\alpha)\subseteq \ann(\A,\D,\alpha|_\D)$.  If $\C$ is
unital, Lemma~\ref{sameunit}\eqref{sameunit1} shows $\A$ is unital and
 $\alpha(I_\C)=I_\A$.
\end{proof}

\begin{lemma}\label{claim2} Let 
      $(\C_1,\D_1)$ be an inclusion.  Suppose $(\C_2,\D_2\ms \alpha)$ is
      an essential 
      expansion of $(\C_1, \D_1)$ and $(\C_2,\D_2)$ has the
      faithful unique pseudo-expectation property. 
        Then
          \begin{enumerate}
            \item \label{claim2:a} $(\C_1,\D_1)$ has the faithful unique pseudo-expectation
              property; and
              \item\label{claim2:b} $\alpha(\D_1^c)\subseteq \D_2^c$.
              \end{enumerate}
            \end{lemma}
        \begin{proof}  Before giving the proofs of \eqref{claim2:a}
          and \eqref{claim2:b}, we make a few remarks.
          Let $(I(\D_1),\iota_1)$ be an injective envelope for $\D_1$.
          By Proposition~\ref{PZC22},  $(\D_2,\D_1,\alpha|_{\D_1})$
          has the
          faithful unique pseudo-expectation property;   let
          $\iota_2: \D_2\rightarrow I(\D_1)$ be the pseudo-expectation
          for $(\D_2,\D_1,\alpha|_{\D_1})$
          relative to $(I(\D_1),\iota_1)$.  Then
          $(I(\D_1),\iota_2)$ is an injective envelope for $\D_2$ (Proposition~\ref{PZC22}).
          Let $$E_2: \C_2\rightarrow I(\D_1)$$ be the pseudo-expectation
          for $(\C_2, \D_2)$ relative to $(I(\D_1),\iota_2)$.

\eqref{claim2:a} Let us first show $(\C_2,\D_1, \alpha|_{\D_1})$ has the faithful
        unique pseudo-expectation property.
Suppose $F: \C_2\rightarrow I(\D_1)$ is a
      pseudo-expectation for $(\C_2,\D_1, \alpha|_{\D_1})$ relative to
      $(I(\D_1), \iota_1)$.  Then for $d\in \D_1$,
      $F(\alpha(d))=\iota_1(d)=\iota_2(\alpha(d))$.  Thus $F|_{\D_2}$
      is a pseudo-expectation for $(\D_2,\D_1, \alpha|_{\D_1})$ relative to
      $(I(\D_1), \iota_1)$.  Therefore
      $F|_{\D_2}=\iota_2$, that is, $F$ is a pseudo-expectation for
      $(\C_2,\D_2)$ relative to $(I(\D_1),\iota_2)$.  This forces
      $F=E_2$, showing $(\C_2,\D_1,\alpha|_{\D_1})$ has the faithful unique
      pseudo-expectation property.
          
Corollary~\ref{heriditary} applied to $(\C_2,\D_1,\alpha|_{\D_1})$
with $\B=\alpha(\C_1)$ implies
        $(\C_1,\D_1)$ has the faithful unique pseudo-expectation
        property (the pseudo-expectation is $E_2\circ\alpha$).

   \eqref{claim2:b}      By Proposition~\ref{f!pse->abelcom}, $\rcom(\C_2, \alpha(\D_1))$ is abelian.   As
        $\alpha(\D_1)\subseteq \D_2$, we find
        \[\rcom(\C_2,\D_2)\subseteq \rcom(\C_2, \alpha(\D_1)).\]  Since
        $(\C_2,\D_2)$ has the faithful unique pseudo-expectation
        property,  
        $\rcom(\C_2,\D_2)$ is a MASA in $\C_2$, so
        $\rcom(\C_2,\D_2)=\rcom(\C_2,\alpha(\D_1))$.   As
        \[\alpha(\D_1^c)\subseteq \rcom(\C_2,\alpha(\D_1))=\D_2^c,\]
        part \eqref{claim2:b} holds.
      \end{proof}

The
        definition of Cartan envelope given in~\cite{PittsStReInII}
         extends   to any regular inclusion $(\C,\D)$.  Here is the
          definition, along with the definitions of  other useful
          expansions.

  \begin{definition}{(cf.~\cite[Definition~5.1]{PittsStReInII})} \label{defCarEnv}
  Let $(\C,\D)$
  be a regular inclusion.
\begin{enumerate}
 \item\label{defCarEnv1}  A \textit{package}\index{Package} for $(\C,\D)$ is a regular
   expansion $(\A,\B\ms \alpha)$ of $(\C,\D)$ such that
  \begin{enumerate}
      \item there is a faithful conditional
      expectation $\Delta:\A\rightarrow\B$; 
    \item $(\C,\D)$ generates $(\A,\B)$ in the sense that
      $\B=C^*(\Delta(\alpha(\C)))$ and $\A=C^*(\alpha(\C)\cup
      \B)$.
    \end{enumerate}
    If in addition,  $(\A,\B)$ is a Cartan inclusion, we say
    $(\A,\B\ms \alpha)$ is a \textit{Cartan package}\index{Package!Cartan} for $(\C,\D)$.
\item\label{defCarEnv2} If $(\A,\B\ms \alpha)$ is a package for
        $(\C,\D)$ which is also an essential expansion for $(\C,\D)$ (see
        Definition~\ref{defextinc}\eqref{defextinc2}), then we call
      $(\A,\B\ms \alpha)$ an \textit{envelope}\index{Envelope} for $(\C,\D)$.
\item \label{defCarEnv3}  We say $(\A,\B \ms \alpha)$ is a \textit{Cartan
    envelope}\index{Envelope!Cartan} for $(\C,\D)$ if it is an envelope and a Cartan package.
        \end{enumerate}
            \end{definition}

We now show that when
$(\A,\B\ms \alpha)$ is an essential  and regular Cartan expansion of $(\C,\D)$, the
inclusion generated by $(\C,\D)$ is a Cartan envelope.  

      \begin{proposition}\label{reg+Del}  Suppose $(\C,\D)$ is a regular
  inclusion  and $(\A,\B\ms \alpha)$ is an essential, regular,
  and Cartan
  expansion for $(\C,\D)$. 
  Let $\Delta: \A\rightarrow \B$ be the
conditional expectation and put 
  \[\B_1:=C^*(\Delta(\alpha(\C))) \dstext{and} \A_1:=C^*(\alpha(\C)\cup\B_1). \]
  Then  $(\A_1,\B_1\ms \alpha)$ is a Cartan envelope for $(\C,\D)$ and
               $\Delta|_{\A_1}: \A_1\rightarrow \B_1$ is the
               conditional expectation.
\end{proposition}
\begin{proof}
By construction, $\alpha(\C)\subseteq \A_1$ and $\alpha(\D)\subseteq
\B_1$, so $(\A_1,\B_1\ms \alpha)$ is an expansion of $(\C,\D)$.   
\begin{flexstate}{\sc Claim}{} \label{clm3} $\alpha: (\C,\D)\rightarrow (\A_1,
        \B_1)$ is a regular map.
      \end{flexstate}
      \noindent \textit{Proof of Claim~\ref{clm3}.}
        For $v\in \N(\C,\D)$, we must
        show that $\alpha(v)$ normalizes $\B_1$. 
 As $(\A,\B)$ is a Cartan inclusion, it has the
faithful unique pseudo-expectation property by Proposition~\ref{Car+upse}.
Lemma~\ref{claim2} shows $(\C,\D)$ has the faithful unique
pseudo-expectation property
and 
\begin{equation}\label{clm3.c} \alpha(\D^c)\subseteq \B^c=\B.
\end{equation}

  An argument
similar to that used for establishing~\cite[(5.14)]{PittsStReInII}
shows that for $n\in\bbN$, and every collection of $n$ elements
$\{x_j\}_{j=1}^n\subseteq \C$,
\begin{equation}\label{nu5.14}
  \alpha(v)^*\left(\prod_{j=1}^n\Delta(\alpha(x_j))\right)\alpha(v)\in
  \B_1.\end{equation}
However, the proofs of~\eqref{nu5.14} and \cite[(5.14)]{PittsStReInII}
have some technical
differences, so we include the proof of~\eqref{nu5.14} here.   The
argument is by induction.  When  $n=1$,
the invariance of $\Delta$ under $\N(\A,\B)$ (see~\eqref{exinv}) gives,
\[\alpha(v)^*\Delta(\alpha(x_1))\alpha(v)=\Delta(\alpha(v^*)\alpha(x_1)\alpha(v))\in
  \Delta(\alpha(\C))\subseteq 
  \B_1.\]
Suppose now that~\eqref{nu5.14} holds for some $n\in \bbN$
and every collection of $n$ elements of $\C$.
Let $\{x_j\}_{j=1}^{n+1}\subseteq \C$ and set $y=\prod_{j=1}^{n}
\Delta(\alpha(x_j))$.
For 
$h\in\overline{vv^*\D^c}$, Lemma~\ref{relcom}\eqref{relcom2}  and
Lemma~\ref{claim2} show $\alpha(h) \in
\alpha(\overline{vv^*\D^c})\subseteq 
\overline{\alpha(vv^*)\B}$;  the induction hypothesis gives $
\alpha(v^*)y\alpha(v)\in \B_1$.
Given $h_0\in vv^*\D^c$, write
$h_0=vv^*k$ where $k\in \D^c$.   Using Lemma~\ref{Tv2.2},  
\begin{align*}
  \alpha(v^*)\>y\Delta(\alpha(x_{n+1}))
  \alpha(h_0)\>\alpha(v)&=\alpha(v^*)y\alpha(v) \>\theta_{\alpha(v)}
                      (\Delta(\alpha(x_{n+1}))\alpha(h_0))\\
&=\alpha(v^*)y\alpha(v) \>\theta_{\alpha(v)}
                                                             (\alpha(vv^*)\Delta(\alpha(x_{n+1}k)
                                                                  )) \\
  &=(\alpha(v^*)y\alpha(v) ) \>   \left(\alpha(v)^*
    \Delta(\alpha(x_{n+1}k))\alpha(v)\right) \\
    &=(\alpha(v^*)y\alpha(v) ) \>   \Delta(\alpha(v)^*
                      \alpha(x_{n+1}k)\alpha(v)) \in \B_1.
\end{align*}
Therefore, for every $h\in \overline{vv^*\D^c}$,
  \begin{equation}\label{nu5.15}
     \alpha(v^*)\>y\Delta(\alpha(x_{n+1}))
    \alpha(h)\>\alpha(v) \in \B_1.
   \end{equation}
   Note that $v=\lim_{n\rightarrow\infty} (vv^*)^{1/n}v$.  Thus for any
$d\in \D^c$, choosing $h=d(vv^*)^{1/n}\in \overline{vv^*\D^c}$ in
\eqref{nu5.15}, we find
\begin{equation}\label{nu5.16}\alpha(v)^* y
  \Delta(\alpha(x_{n+1}))  \alpha(d)\alpha(v) = \lim_n \alpha(v)^* y
  \Delta(\alpha(x_{n+1}))  \alpha(d (vv^*)^{1/n})\alpha(v)\in\B_1.
\end{equation}
As $\D^c$ is abelian (by~\eqref{clm3.c}), $(\C,\D^c)$ is a MASA
inclusion.
Lemma~\ref{relcom}\eqref{relcom1} shows that $(\C,\D^c)$ is
a regular MASA inclusion, and therefore has the AUP
by~\cite[Theorem~2.5]{PittsNoApUnInC*Al}.  Let $(u_\lambda)\subseteq
\D^c$ be an approximate unit 
for $\C$.  Taking $d=u_\lambda$ in~\eqref{nu5.16}, we conclude
\begin{align*}\alpha(v)^* y
  \Delta(\alpha(x_{n+1}))\alpha(v)&=\lim_\lambda \alpha(v)^* y
  \Delta(\alpha(x_{n+1}u_\lambda))\alpha(v)\\&=\lim_\lambda \alpha(v)^* y
  \Delta(\alpha(x_{n+1}))\alpha(u_\lambda)\alpha(v)\in\B_1.
\end{align*}
Thus~\eqref{nu5.14} holds.

Since $\B_1$ is generated by $\Delta(\alpha(\C))$,~\eqref{nu5.14}
implies $\alpha(v)^*\B_1\alpha(v)\subseteq \B_1$.  As this holds for
every $v\in \N(\C,\D)$, we may replace $v$ with $v^*$ to obtain
$\alpha(v)\B_1\alpha(v)^*\subseteq \B_1$.  Therefore, $\alpha(v)\in \N(\A_1,\B_1)$.
\hfill $\diamondsuit$
\vskip 6pt

By Claim~\ref{clm3}, $\alpha(\N(\C,\D))\cup\B_1\subseteq
\N(\A_1,\B_1)$, whence $(\A_1,\B_1)$ is a regular inclusion.
Also,
\[\alpha(\D)\subseteq \B_1\subseteq \B.\]
Applying Lemma~\ref{interiip}, $(\B,\B_1)$ and $(\B_1,\D,\alpha|_\D)$
have the \iip.  Thus,  $(\A,\B\ms  \, \subseteq)$ is an essential
expansion of $(\A_1,\B_1)$.  By Lemma~\ref{claim2}, $(\A_1,\B_1)$ has
the faithful unique pseudo-expectation property.

We next show $\Delta(\A_1)=\B_1$; since $\B_1\subseteq \B$, the
non-trivial inclusion is 
\begin{equation}\label{reg+DelA}
  \Delta(\A_1)\subseteq \B_1.
\end{equation}
Let $\P\subseteq \A_1$ be the $*$-semigroup generated by
$\alpha(\N(\C,\D))\cup\B_1$.  By Claim~\ref{clm3}, $\P\subseteq \N(\A_1,\B_1)$.
Let \[Y:=\{w\in \P: \Delta(w)\in \B_1\}.\]  For $n\in\bbN$, set
\[X_1:=\alpha(\N(\C,\D)) \dstext{ and } 
  X_{n+1}:=X_1\B_1 X_n.\]  Note that
each element of
  $X_{n+1}$ has the form,
  \begin{equation}\label{reg+DelA.1}
    \alpha(v_0)\left(\prod_{j=1}^nb_j\alpha(v_j)\right),
  \end{equation}
  where
  $\{b_j\}_{j=1}^n\subseteq\B_1$ and  $\{v\}_{j=0}^n\subseteq
  \N(\C,\D)$.

  We claim that for $n\in \bbN$,
$X_n\subseteq Y$.  The argument is again by induction. 
By definition of $\B_1$, $X_1\subseteq Y$.
Suppose that $n\in\bbN$ and $X_n\subseteq Y$.   Choose $w\in X_{n+1}$
and let $\eps>0$.  We may find $u\in X_n$, $b\in \B_1$ and $v\in
\N(\C,\D)$ so that
\[w=\alpha(v^*)bu.\]  As $v^*=\lim v^*(vv^*)^{1/m}$ we may find $m\in
\bbN$ so that
\[\norm{w-\alpha(v^*)\alpha(vv^*)^{1/m}bu}<\eps.\]  Note that
$\alpha(v^*)u\in X_n$ by definition of $X_n$ and $\B_1X_n\subseteq Y$  because $\Delta$ is a
$\B_1$-bimodule map.  Since $\alpha(vv^*)^{1/m}b\in
\overline{\alpha(v)\alpha(v)^*\B_1}$, Lemma~\ref{Tv2.2} gives 
\[\alpha(v^*)\alpha(vv^*)^{1/m}
  bu=\theta_{\alpha(v)}(\alpha(vv^*)^{1/m}b) \,\, \alpha(v^*)u\in
  \B_1X_n\subseteq Y.\]
Thus every neighborhood of $\Delta(w)$ contains an element of
$\Delta(Y)\subseteq \B_1$.  Therefore, $\Delta(w)\in \B_1$, that is,
$w\in Y$.  The claim follows.

  Now put $X:=\bigcup_{n=1}^\infty X_n$.   Note that for $n\in
  \bbN$, $X_n^*=X_n$ and  
  $X_1X_n\subseteq X_n\supseteq X_nX_1$.  Also, ~\eqref{reg+DelA.1}
  shows that for $n, m\in\bbN$, there exists $k\in\bbN$
  so that $X_nX_m\subseteq X_k$.   It follows that \[\S:=\B_1\cup X\cup \B_1X\cup
    X\B_1\cup \B_1X\B_1\] is a $*$-semigroup containing $\B_1\cup
  \alpha(\N(\C,\D))$.  As $\S\subseteq \P$, we obtain $\S=\P$.     The bimodule property of $\Delta$ gives 
  $\P=Y$. 
  Finally, since $\spn\P$ is a $*$-algebra containing $\alpha(\N(\C,\D))\cup
  \B_1$, 
  \[\overline\spn(\P) =\A_1,\] which gives~\eqref{reg+DelA}.

 It is now evident that  $\Delta|_{\A_1}:\A_1\rightarrow \B_1$ is a
 conditional expectation.  Proposition~\ref{Car+upse} shows  $(\A_1,\B_1)$ is a Cartan
  inclusion.   Therefore,  $(\A_1,\B_1\ms \alpha)$ is a Cartan package for
  $(\C,\D)$.  As $(\B_1,\D,\alpha|_\D)$ has the \iip,
  $(\A_1,\B_1\ms \alpha)$ is  also an envelope for 
  $(\C,\D)$.  Thus $(\A_1,\B_1\ms \alpha)$ is a Cartan envelope
    for $(\C,\D)$.
\end{proof}

The following will be used in the proof of Theorem~\ref{!pschar}; 
we will establish its converse in Proposition~\ref{CE4rcom}.  

\begin{proposition}\label{rcie}  Suppose $(\C,\D)$ is a regular
  inclusion such that $(\D^c,\D)$ has the \iip\ and $\D^c$ is abelian.
  If $(\A,\B\ms \alpha)$ is a Cartan envelope for $(\C,\D^c)$, then
  $(\A,\B\ms \alpha)$ is also a Cartan envelope for $(\C,\D)$.
\end{proposition}
\begin{proof}
Let $\Delta:\A\rightarrow\B$ be the conditional
  expectation.   

  Since $\D\subseteq \D^c$, $(\A,\B\ms \alpha)$ is an expansion of
  $(\C,\D)$.  Lemma~\ref{relcom}\eqref{relcom1} and the regularity of
  $\alpha:(\C,\D^c)\rightarrow (\A,\B)$ show
  $\alpha:(\C,\D)\rightarrow (\A,\B)$ is a regular mapping.  Since
  $(\A,\B\ms \alpha)$ is a Cartan envelope for $(\C,\D^c)$,
  $\B=C^*(\Delta(\alpha(\C)))$ and $\A=C^*(\B\cup\alpha(\C))$.  Thus,
  $(\A,\B\ms \alpha)$ is a Cartan package for $(\C,\D)$.
    
  We have $\alpha(\D)\subseteq \alpha(\D^c)\subseteq \B$.
  By definition of Cartan envelope, $(\B,\D^c,\alpha)$ has the \iip.
  Since $(\D^c,\D)$ has the \iip, so does $(\B,\D,\alpha)$
  (Lemma~\ref{interiip}).  Hence
  $(\A,\B\ms \alpha)$ is a Cartan envelope for $(\C,\D)$.
 \end{proof}

We turn to some technical preparations for the proof of
the uniqueness statement for Cartan envelopes found in
Theorem~\ref{!pschar} below.  We will use Lemma~\ref{Tv2.2} repeatedly, often
without comment. The purpose of Proposition~\ref{*semigp} is to 
produce useful $*$-semigroups of intertwiners from $*$-semigroups of
normalizers; these semigroups 
are key to the uniqueness statement in Theorem~\ref{!pschar}.
Lemma~\ref{cint} exhibits some intertwiners associated to a normalizer
and is the foundation for the proof of Proposition~\ref{*semigp}.
Recall from Definition~\ref{innore}\eqref{innore3.5} that $\I(\C,\D)$
denotes the collection of all intertwiners for $(\C,\D)$.

 When $(\C,\D)$ is an inclusion and  $v\in \N(\C,\D)$, it is convenient to
use the notation,
\begin{equation}\label{RnS}
  R(v):=\overline{vv^*\D} \dstext{and} S(v):=\overline{v^*v\D}.
\end{equation}

For an abelian \cstaralg\ $\D$ and an  ideal $J\idealin \D$, let
\begin{equation}\label{J_c} J_c:=\{d\in J: \text{ there exists } h\in
  J \text{ such that } dh=d\}.
\end{equation}
\begin{flexstate}{Observation}{}\label{JcDense}
$J_c=J_c^*$ and $J_c$ is an algebraic ideal of $\D$ whose
closure is $J$.
(In general,
$J_c$ is not closed.)
\end{flexstate}
\begin{proof}[Sketch of Proof]
Since $J\simeq C_0(U)$ for some open subset $U$ of $\hat\D$, we may 
identify $J$ with $C_0(U)$.  Write $C_c(U)=\{f\in C_0(U):
\overline{\supp(f)} \text{ is a compact subset of } U\}$.  The
observation will follow once we establish
\[\{d\in C_0(U): \exists \, h\in C_0(U) \text{ such that }
dh=d\}= C_c(U).\] If $d\in
C_0(U)$ and there exists $h\in C_0(U)$ such that $dh=d$, then $h(x)=1$
for 
every $x\in \overline{\supp(d)}$.
As $\{u\in U: |h(u)|\geq 1\}$ is a compact set, we find 
$\overline{\supp(d)}$ is compact, whence $d\in C_c(U)$.
Thus $\{d\in C_0(U): \exists \, k\in C_0(U) \text{ such that }
kd=d\}\subseteq C_c(U)$.   For the reverse inclusion, given $d\in
C_c(U)$,   an application of
Urysohn's lemma shows there exists $h\in C_0(U)$ such that $h(x)=1$
for every $x\in \overline{\supp(d)}$.  
\end{proof}

\begin{lemma} \label{cint}  Suppose $(\C,\D)$ is an inclusion.   Let 
 $v, w\in \N(\C,\D)$.
  \begin{enumerate}
 \item \label{cint2} Suppose $vv^*\in \D$ (resp.\ $v^*v\in \D$).
      Then 
      given $\eps>0$, there exists $k\in R(v)_c$, (resp.\ $k\in
      S(v)_c$) such that $\norm{v-kv}<\eps$ (resp.\ $\norm{v-vk}<\eps$). 
  \item\label{cint0} If $k\in R(v)_c$ (resp.\
    $k\in S(v)_c$), then there exists
    $s\in R(v)_c$ (resp.\
    $s\in S(v)_c$) such that $k= vv^*s$ (resp.\ $k=v^*vs$).
  \item \label{cint1} If $k\in R(v)_c$ (resp.\ $k\in
      S(v)_c$),
      then $kv\in \I(\C,\D)$ (resp.\ $vk\in \I(\C,\D)$).
\item\label{cint2.5}  For $h\in R(v)_c$ (resp.\ $h\in S(v)_c$),
  $\theta_v(h^*)\in R(v^*)_c$ (resp.\ $\theta_v^{-1}(h^*)\in S(v^*)_c$).
    \item \label{cint3}  Suppose $h\in R(v)_c$ and  $k\in R(w)_c$
        (resp.\ $h\in S(v)_c$ and $k\in S(w)_c)$.  Then
        $\theta_v^{-1}(\theta_v(h)k)\in R(vw)_c$ (resp.\ $\theta_w(h \theta_w^{-1}(k))\in S(vw)_c$).
      \end{enumerate}
    \end{lemma}
Before beginning the proof, we make a comment.  In general, when $v\in\N(\C,\D)$,
$vv^*$ need not belong to $\D$. However given $\sigma\in \hat\D$, we may use
\cite[Proposition~2.1]{PittsNoApUnInC*Al} to define $\sigma(vv^*)$:
choose $d\in\D$ with $\sigma(d)=1$ and then set
$\sigma(vv^*):=\sigma(vv^*d)$.  This is independent of the choice of
  $d$ by \cite[Fact~1.2]{PittsNoApUnInC*Al}.  Note also that the map
  $\hat\D\ni\sigma\mapsto \sigma(vv^*)$ is continuous on $\hat\D$;
  indeed, given
  $\sigma\in\hat\D$, consider
  $d\in\D$ chosen so that $\tau(d)=1$ for all $\tau$ in an open 
  neighborhood of $\sigma$. 

    \begin{proof}
We give the proofs  for $vv^*$; the proofs for $v^*v$
are similar.

\eqref{cint2}   Suppose $vv^*\in \D$ and let $\eps>0$.  Since $\widehat{vv^*}\in
C_0(\hat\D)$, 
there exists a compact set $C\subseteq \hat\D$ such that
$\sigma(vv^*)<\eps^2/4$ whenever $\sigma\in \hat\D\setminus C$.  Let 
\[K:= C\cap \{\sigma\in \hat\D: \sigma(vv^*)\geq \eps^2/4\}.\] Then $K$
is a compact subset of $\hat\D$.   By local
compactness of $\hat\D$, we may find  open sets $U, V\subseteq \hat\D$
such that:
\begin{itemize}
\item[]  $\overline U$ and $\overline V$ are  compact;
\item[]  $K\subseteq U \subseteq \overline U\subseteq \{\sigma\in\hat\D:
\sigma(vv^*)\geq \eps^2/9\}$; and
\item[]  $\overline U \subseteq V\subseteq \overline V \subseteq \{\sigma\in
\hat\D: \sigma(vv^*)\geq \eps^2/16\}$. 
\end{itemize}
Urysohn's lemma ensures  there exists $h, k\in \{d\in \D: d\geq 0
\text{ and } \norm{d}\leq 1\}$ such that
$\hat k\equiv 1$ on $K$, $\hat k$ vanishes off $ U$, $\hat
h\equiv 1$ on $\overline
U$, and $\hat h$ vanishes off $V$. Then $h, k\in R(v)$ and $hk=k$, so
that $k\in R(v)_c$.  

Note that for $\sigma\in\hat\D\setminus  K$, $\sigma(vv^*)<\eps^2/4$.   Also,
 $(kv-v) (kv-v)^*=k^2vv^*-2kvv^*+vv^*\in \D$,
so 
\[\norm{kv-v}^2=\sup_{\sigma\in \hat\D}\sigma(({kv-v})(kv-v)^*)
  =\sup_{\sigma\in \hat\D}  \sigma(vv^*)(\sigma(k)-1)^2\leq
\eps^2/4<\eps^2.\]

\eqref{cint0} 
Let \[\ran(v)=\{\sigma\in \hat\D: \sigma(vv^*)\neq
0\}.\] Then \[R(v)=\{d\in \D: \hat d \text{ vanishes
  off }\ran(v)\}\cong  C_0(\ran(v)).\]  As $k\in R(v)_c$, we may find $h\in
R(v)$ so that $kh=k$.    We may therefore find a
compact set $F\subseteq \ran(v)$ such that $|\sigma(h)|< 1/2$ for
$\sigma\in \hat\D\setminus F$.  As $kh=k$, $\hat k$
vanishes off $F$.  Define a function $\delta$ on $\hat\D$ by
\[\delta(\sigma)=
  \begin{cases} \ds \frac{\sigma(k)}{\sigma(vv^*)} & \text{for } \sigma\in F;\\
    0& \text{for }\sigma\in \hat\D\setminus F.\end{cases}\]   Since
$\widehat{vv^*}$ is bounded away from $0$ on $F$, $\delta\in
C_0(\hat\D)$, and we define $s\in \D$ by $\hat s=\delta$.  By
construction, $k=vv^* s$, and since $\hat s \hat h=\hat s$,  we find
$s\in R(v)_c$. 

      \eqref{cint1} Let $k\in R(v)_c$ and choose 
$h\in R(v)=\overline{vv^*\D}$ such that $kh=k$.   Set $u:=kv$.   Then for
$d\in \D$ we have
\[du=k(dhv)=kv\theta_v(dh)=u\theta_v(dh).\] Thus, $\D u\subseteq u\D$.

For the reverse inclusion, note that
$\theta_v(kh)=\theta_v(k)\theta_v(h)=\theta_v(k)$, so that
$\theta_v(k)\in S(v)_c$.   Thus, for $d\in \D$ we have
\[ud=kv\,(\theta_v(h) d)=\theta_v^{-1}(\theta_v(h) d) kv\in \D u.\]
Therefore $kv\in \I(\C,\D)$, as claimed.

\eqref{cint2.5} This follows from the fact that
$\theta_v:R(v)\rightarrow S(v)$ is a $*$-isomorphism.  If
$h\in R(v)_c$ and $k\in R(v)$ satisfies $hk=h$; taking adjoints and
applying $\theta_v$ gives $\theta_v(h^*)\theta_v(k^*)=\theta_v(h^*)$,
so $\theta_v(h^*)\in S(v)_c=R(v^*)_c$.

\eqref{cint3}  Using part~\eqref{cint0} we may choose $s\in R(v)_c$,
$t\in R(w)_c$ so that $h=(vv^*)^2s$ and $k=ww^*t$.  Then
\begin{align*}
  \theta_v^{-1}(\theta_v(h)k)&= \theta_v^{-1}( \theta_v((vv^*)^2 s) k)=
                               \theta_v^{-1}(v^*(vv^*s)vk)\\
                             &=vv^*svkv^*=(vv^*)sv(tww^*)v^*\\
  &= s (vv^*)v tww^*v^*\\
&=  s\,  (vtv^*) (v w)(vw)^*\in (vw)(vw)^*\D\subseteq R(vw).
\end{align*}
By definition, we may find $x\in R(v)$ and $y\in R(w)$ so that $h=hx$
and $k=ky$.  Since $\overline{R(v)_c}=R(v)$ and
$\overline{R(w)_c}=R(w)$ (see Observation~\ref{JcDense}), we may find sequences $(x_k)$ in $R(v)_c$
and $(y_k)$ in $R(w)_c$ which converge to $x$ and $y$ respectively.
The previous calculations applied to $x_k$ and $y_k$ show
$\theta_v^{-1}(\theta_v(x_k)y_k)\in R(vw)$.  Taking limits yields
$\theta_v^{-1}(\theta_v(x)y)\in R(vw)$.
But then
\[\theta_v^{-1}(\theta_v(h)k)\theta_v^{-1}(\theta_v(x)y)=\theta_v^{-1}(\theta_v(h)k),\]
showing that $\theta_v^{-1}(\theta_v(h)k)\in R(vw)_c$.
\end{proof}

\begin{proposition}\label{*semigp}  Let $(\C,\D)$ be an inclusion.  Let $\P\subseteq \N(\C,\D)$ be a $*$-semigroup and put 
  \[\Q:=\{vh: v\in \P\text{ and } h\in S(v)_c\}\dstext{and} \S:=\Q\cup
    \D.\]  The following statements hold.
  \begin{enumerate}
  \item \label{*semigp1} $\Q$ and $\S$ are $*$-semigroups and both are
    contained in
    $\N(\C,\D)\cap \I(\C,\D)$.
    \item \label{*semigp3} When $(\C,\D)$  has the approximate
  unit property and the $*$-algebra generated by $\P\cup\D$  
   is dense in $\C$,  $\spn\S$ is  a
   dense $*$-subalgebra of $\C$.
 \end{enumerate}
\end{proposition}
\begin{proof}
\eqref{*semigp1} That $\Q\subseteq \N(\C,\D)$ is clear and Lemma~\ref{cint}\eqref{cint1} shows
$\Q\subseteq  \I(\C,\D)$.  For $v, w\in \P$, $h\in S(v)_c$ and $k\in S(w)_c$,
  \[vhwk=vh\theta_w^{-1}(k)w= vw \theta_w(h\theta_w^{-1}(k))\in \Q,\]
  by   Lemma~\ref{cint}\eqref{cint3}.
  Next, Lemma~\ref{cint}\eqref{cint2.5} implies $\Q$ is closed under adjoints: indeed,
    $(vh)^*=h^*v^*=v^*\theta_v^{-1}(h^*)\in \Q$.   Thus $\Q$ is a
    $*$-semigroup.  

Turning our attention to $\S$, let $v\in \P$ and $h\in S(v)_c$. For $d\in \D$,  the
fact that $vh\in
\I(\C,\D)$ gives $d(vh)\in \Q$ and $(vh)d\in \Q$.  Thus, $\Q \D \cup
\D\Q \subseteq \Q$.  It follows that $\S$ is a $*$-semigroup. Since $\N(\C,\D)$ and $\I(\C,\D)$ are $*$-semigroups
containing $\D$, $\S\subseteq \N(\C,\D)\cap \I(\C,\D)$.

\eqref{*semigp3} Note that $(\C,\D)$ is a regular inclusion because
$\P\cup\D$ generates a $*$-semigroup of normalizers.  As $(\C,\D)$ has
the AUP, \cite[Observation~1.3]{PittsNoApUnInC*Al} shows
  $v^*v\in \D$ for every $v\in \N(\C,\D)$.     
  Lemma~\ref{cint}\eqref{cint2} shows $v\in \overline{\S}$ for every
  $v\in \P$.  Thus $\overline{\spn\S}$ contains  $\P\cup\D$,
 whence  $\overline{\spn \S}=\C$.  Finally, as $\S$ is a
  $*$-semigroup, $\spn\S$ is a $*$-algebra.
\end{proof}

By taking $\P=\N(\C,\D)$, we see that for a broad class of inclusions,
the supply of
intertwiners is large.
\begin{corollary}\label{deninter}  If $(\C,\D)$ is a regular inclusion
  with the AUP, then $\I(\C,\D)\cap \N(\C,\D)$ has dense span in $\C$.
\end{corollary}

We come to a main result, which is the promised characterization of those
regular inclusions which have a Cartan envelope.  This characterization significantly
extends \cite[Theorem~5.2]{PittsStReInII}, which dealt with unital
inclusions.

\begin{theorem}[c.f.~{\cite[Theorem~5.2]{PittsStReInII}}]\label{!pschar}   Let $(\C,\D)$ be a regular
inclusion.  The following are equivalent.
\begin{enumerate}
\item \label{!pschar0}   $(\C,\D)$ has a Cartan envelope.
\item \label{!pschar1} $(\C,\D)$ has the faithful unique pseudo-expectation property.
\item \label{!pschar2} $\D^c$ is abelian and both $(\C,\D^c)$ and 
  $(\D^c,\D)$ have the \iip.
\end{enumerate}

When $(\C,\D)$ satisfies any of  conditions (a)--(c), $(\C,\D)$ is
\wnd\ and the following
statements hold.
\begin{description}
\setcounter{enumi}{3}
\item[\sc Uniqueness] If for $j=1,2$, $(\A_j,\B_j\ms \alpha_j)$ are Cartan envelopes for
    $(\C,\D)$, there exists a unique regular $*$-isomorphism
    $\psi: (\A_1,\B_1)\rightarrow (\A_2,\B_2)$ such that $\psi\circ\alpha_1=\alpha_2$.

\item[\sc Minimality] If $(\A,\B\ms  \alpha)$ is a Cartan  
      package for $(\C,\D)$, there is an ideal $\fJ\subseteq \A$ such
      that $\fJ\cap \alpha(\C)=\{0\}$ and, letting
      $q:\A\rightarrow \A/\fJ$ denote the quotient map, 
      $(\A/\fJ, \B/(\fJ\cap \B)\ms  q\circ\alpha)$ is a Cartan
      envelope for $(\C,\D)$.
    \end{description}
\end{theorem}
The proof of Theorem~\ref{!pschar} is organized as follows.  We  show
\eqref{!pschar0}$\Rightarrow$\eqref{!pschar1}$\Rightarrow$\eqref{!pschar2}$\Rightarrow$\eqref{!pschar0}
and 
weak non-degeneracy for $(\C,\D)$, then the uniqueness statement.
The proof of the minimality statement may be found after the proof of Proposition~\ref{invPsE}.
During all parts of the proof of Theorem~\ref{!pschar}, assume that 
  an injective envelope $(I(\D),\iota)$ for $\D$ has been fixed.

  \begin{proof}[\eqref{!pschar0}$\Rightarrow$\eqref{!pschar1}] 
Since $(\C,\D)$ has a
               Cartan envelope, Lemma~\ref{cewnd} shows it is \wnd.   
 
  Suppose $(\A,\B\ms \alpha)$ is a Cartan envelope for $(\C,\D)$.
  By Proposition~\ref{Car+upse}, $(\A,\B)$ has the faithful unique
  pseudo-expectation property. By definition, $(\A,\B\ms \alpha)$ is an essential expansion of
  $(\C,\D)$,  so Lemma~\ref{claim2} shows $(\C,\D)$
   has the faithful unique pseudo-expectation property.
 \end{proof}

\begin{proof}[\eqref{!pschar1}$\Rightarrow$\eqref{!pschar2}]  By
Proposition~\ref{f!pse->abelcom}, $\D^c$ is abelian, so $\D\subseteq
\D^c\subseteq \C$.   Corollary~\ref{heriditary} shows $(\D^c,\D)$ has
the faithful unique pseudo-expectation property, whence $(\D^c,\D)$
has the \iip\ by Proposition~\ref{PZC22}.

Let $E:\C\rightarrow I(\D)$ be the faithful unique pseudo-expectation.
Since $\iota_2:=E|_{\D^c}$ is a pseudo-expectation for $(\D^c,\D)$,
Proposition~\ref{PZC22} shows $\iota_2:=E|_{\D^c}$ is the only
pseudo-expectation for $(\D^c,\D)$, and $(I(\D),\iota_2)$ is an
injective envelope for $\D^c$.

It follows from Lemma~\ref{relcom}\eqref{relcom1} that $(\C,\D^c)$ is
regular MASA inclusion, whence  $(\C,\D^c)$ has the  unique
pseudo-expectation property by Theorem~\ref{regmasaps}.   Letting
$F:\C\rightarrow I(\D)$ be the pseudo-expectation for $(\C,\D^c)$
relative to $(I(\D),\iota_2)$, we see that $F|_{\D^c}$ is a
pseudo-expectation for $(\D^c,\D)$, so $F|_{\D^c}=\iota_2$.   But then
$F|_\D=\iota_2|_{\D}=\iota$, hence $F$ is a pseudo-expectatation for
$(\C,\D)$.  This forces $F=E$.   As $E$ is faithful, we conclude that
$(\C,\D^c)$ has the faithful unique pseudo-expectation property.
Theorem~\ref{regmasaps} now shows $(\C,\D^c)$ has
the \iip. 
\end{proof}

\begin{proof}[\eqref{!pschar2}$\Rightarrow$\eqref{!pschar0}]
For notational convenience, let $\M:=\D^c$.   By Lemma~\ref{->ndeg}\eqref{->ndeg3}, $(\C,\D)$ is \wnd; thus
  Lemma~\ref{sameunit}\eqref{sameunit5} shows $(\tilde\C,\tilde\D)$ is a unital
  inclusion, and Corollary~\ref{interinc} gives 
  $\tilde\D\subseteq \tilde\M\subseteq \tilde\C$.
  By~\cite[Theorem~2.5]{PittsNoApUnInC*Al}, $(\C,\M)$ has the AUP, so
 \begin{equation*}\N(\C,\D)
 \stackrel{\ref{relcom}\eqref{relcom1}}{\subseteq}
 \N(\C,\M)
 \stackrel{\ref{relcom}\eqref{relcom3}}{\subseteq}
 \N(\tilde\C,\tilde\M).
\end{equation*}  Therefore, $(\C,\M)$ is a regular MASA inclusion; it
has the \iip\ by hypothesis.    It follows from Lemma~\ref{iipiffuiip} that $(\tilde\C,\tilde\M)$ is also a regular
MASA inclusion with the \iip.

By~\cite[Theorem~5.2]{PittsStReInII}, $(\tilde\C,\tilde\M)$ has a
Cartan envelope $(\A,\B\ms \alpha)$.  Since $(\B, \alpha(\tilde\M))$ and
$(\tilde\M, \M)$ have the \iip, so does $(\B, \alpha(\M))$ (Lemma~\ref{interiip}).  Thus 
$(\A,\B\ms \alpha|_{\C})$ is an essential and regular Cartan expansion  for $(\C,\M)$. 
Proposition~\ref{reg+Del} 
shows that $(\C,\M)$ has a
  Cartan envelope.   Finally, Proposition~\ref{rcie} implies $(\C,\D)$
  has a Cartan envelope.
\end{proof}

\begin{proof}[Proof of the 
 Uniqueness Statement in Theorem~\ref{!pschar}] \rm The proof follows the same
pattern as the proof of \cite[Proposition~5.24]{PittsStReInII}.  Suppose for $i=1,2$ that
$(\A_i,\B_i\ms \alpha_i)$ are Cartan envelopes for $(\C,\D)$.   The
equivalence of statements~\eqref{!pschar0} and \eqref{!pschar1} in
Theorem~\ref{!pschar} shows 
$(\C,\D)$ has the faithful unique pseudo-expectation property.   Since
$(\A_i,\B_i\ms \alpha)$ is an essential expansion of $(\C,\D)$, 
$(\B_i,\D, \alpha_i|_{\D})$ has the \iip, so it has the faithful unique
pseudo-expectation property (Proposition~\ref{PZC22}).  Thus
 there exist unique $*$-monomorphisms
$\iota_i:\B_i\rightarrow I(\D)$ such that $\iota_i\circ(\alpha_i|_\D)=\iota$.
As $\iota_i\circ\Delta_i\circ\alpha_i$ is a pseudo-expectation for
$(\C,\D)$, we obtain
\begin{equation} \label{!pscharu}
  \iota_1\circ\Delta_1\circ\alpha_1=E=\iota_2\circ\Delta_2\circ\alpha_2.
\end{equation} 
But $\B_i=C^*(\Delta_i(\alpha_i(\C)))$, so  $\iota_i(\B_i)=C^*(E(\C))$.
As $\iota_i$ is a $*$-monomorphism and is the unique pseudo-expectation
for $(\B_i,\D, \alpha_i|_{\D})$, we conclude  $\psi_\B:=\iota_2^{-1}\circ \iota_1$
is the unique 
$*$-isomorphism of $\B_1$ onto $\B_2$ satisfying
\[\psi_\B\circ \alpha_1|_{\D}=\alpha_2|_\D.\]

The following
diagram illustrates the maps just discussed.
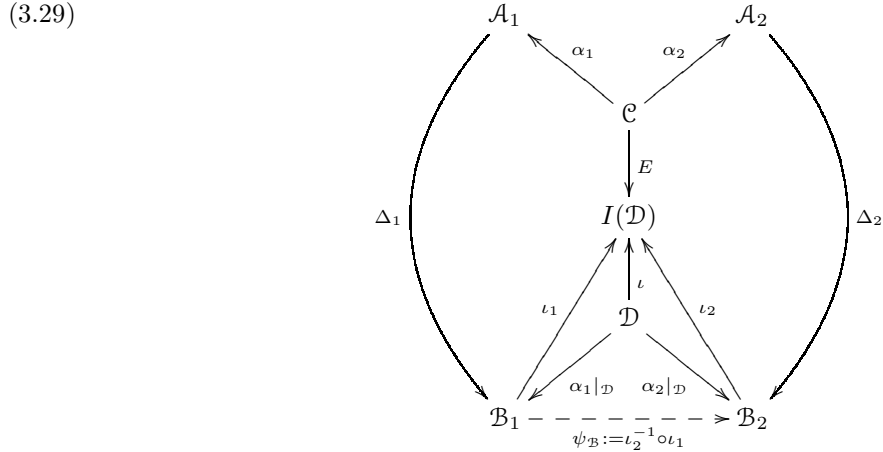
\begin{equation}\label{unCD}
  \xymatrix{
    \A_1\ar@/_3pc/[dddd]_{\Delta_1}&& \A_2\ar@/^3pc/[dddd]^{\Delta_2}\\
    & \C\ar[lu]_{\alpha_1}\ar[ru]^{\alpha_2}\ar[d]^{E}& \\
    & I(\D)&\\
    &\D\ar[u]_(.35)\iota\ar[dl]^{\alpha_1|_{\D}}\ar[dr]_{\alpha_2|_{\D}}&\\
    \B_1\ar[uur]^{\iota_1}\ar@{-->}[rr]_{\psi_\B:=\iota_2^{-1}\circ \iota_1}&&\B_2 \ar[uul]_{\iota_2}}
\end{equation}
We shall extend $\psi_\B$ to the desired isomorphism of $\A_1$ onto
$\A_2$.  We do this in stages.

By the regularity of the map $\alpha_i$, \[\P_i:=\alpha_i(\N(\C,\D))\] is a
$*$-semigroup contained in $\N(\A_i,\B_i)$.  Since Cartan inclusions
have the AUP,
\cite[Observation~1.3]{PittsNoApUnInC*Al}
shows
that for $v\in \N(\C,\D)$, $\alpha_i(v^*v)\in \B_i$.
Also, 
\eqref{!pscharu} gives
\begin{equation}\label{!pscharu1}
  \psi_\B\circ\Delta_1\circ\alpha_1=\Delta_2\circ\alpha_2.
\end{equation}

Put 
\[ \Q_i:= \{\alpha_i(v)h: v\in \N(\C,\D), h\in
  (\overline{\alpha_i(v^*v)\B_i})_c \}\dstext{and}
  \S_i:=\B_i\cup\Q_i.\] 
Since  $(\A_i,\B_i)$ has the AUP and $\A_i$ is generated by $\B_i\cup\P_i$,
Proposition~\ref{*semigp} shows that $\Q_i$ and $\S_i$ are  $*$-semigroups
contained in $\N(\A_i,\B_i)$ and  $\spn\S_i$ is a dense
$*$-subalgebra of $\A_i$.\footnote{The proof
  of~\cite[Proposition~5.24]{PittsStReInII} states without proof that the set
  $\M=\{\alpha_i(v)h: v\in \N(\C,\D), h\in
  \overline{\alpha_i(v^*v)\D_i}\}$ is a $*$-semigroup and uses $\M$
  instead of $\Q$.}    (Note that $\spn \Q_i$ is also a 
$*$-subalgebra of $\spn\S_i$.)

The first stage in extending $\psi_\B$ is to extend it to an isomorphism 
$\psi_\Q: \spn \Q_1\rightarrow \spn\Q_2$. 
Let $n\in \bbN$ and suppose for
$1\leq k\leq n$, $v_k\in \N(\C,\D)$ and
$h_k\in (\overline{\alpha_1(v_k^*v_k)\B_1})_c$.
Then
\begin{align*} 0=\sum_{k=1}^n \alpha_1(v_k)h_k&
   \iff \Delta_1\left(\sum_{k, \ell=1}^n
                    h_k^*\alpha_1(v_k^*v_\ell) h_\ell\right)=0\\
                                          &       \iff \psi_\B\left(\sum_{k, \ell=1}^n
                    h_k^*\Delta_1(\alpha_1(v_k^*v_\ell))
                                            h_\ell\right)=0\\
                  &\stackrel{\eqref{!pscharu1}}{\iff} \sum_{k, \ell=1}^n
                    \psi_\B(h_k^*)\Delta_2(\alpha_2(v_k^*v_\ell)) \psi_\B(h_\ell)=0\\
                  &\iff \Delta_2\left(\sum_{k, \ell=1}^n
                    \psi_\B(h_k^*)\alpha_2(v_k^*v_\ell)
                    \psi_\B(h_\ell)\right)=0\\
  &\iff \sum_{k=1}^n \alpha_2(v_k)\psi_\B(h_k) =0.
\end{align*}
Therefore, $\psi_\B$
extends uniquely to a $*$-isomorphism
$\psi_\Q:\spn(\Q_1)\rightarrow\spn(\Q_2)$ given by
\[\sum_{k=1}^n\alpha_1(v_k)h_k\mapsto\sum_{k=1}^n\alpha_2(v_k)\psi(h_k).\]

For $v\in \N(\C,\D)$ and $h\in(\overline{\alpha_1(v^*v)\B_1})_c$,
  \begin{align*}
    \psi_\B(\Delta_1(\alpha_1(v)h))
    &=\iota_2^{-1}(\iota_1(\Delta_1(\alpha_1(v))h))=\iota_2^{-1}(E(v)\iota_1(h))\\
     & =\iota_2^{-1}(E(v)) \psi_\B(h)
    =
       \iota_2^{-1}(\iota_2(\Delta_2(\alpha_2(v))))\psi_\B(h)\\
    &=\Delta_2(\alpha_2(v)\psi_\B(h))
      =\Delta_2(\psi_\Q(\alpha_1(v)h)).
  \end{align*}

  It follows that 
  \begin{equation}\label{psiIDelta}\psi_\B\circ\Delta_1|_{\spn\Q_1}=\Delta_2\circ\psi_\Q.
\end{equation}

Next we  extend $\psi_\Q$ to a $*$-isomorphism  $\psi_\S: \spn
\S_1\rightarrow \spn\S_2$.
For $b_1\in \B_1$ and $t_1\in \spn
\Q_1$ we have
\begin{align*}
\begin{split}  0=b_1+t_1
  \iff & \Delta_1((b_1+t_1)^*(b_1+t_1))=0\\
  \iff & b_1^*b_1+b_1^*
         \Delta_1(t_1)+\Delta_1(t_1)^*b_1+\Delta_1(t_1^*t_1)
         =0\\
  \stackrel{\eqref{psiIDelta}}{\iff} &\psi_\B(b_1^*b_1)+\psi_\B(b_1)^*\Delta_2(\psi_\Q(t_1))+
         \Delta_2(\psi_\Q(t_1^*))\psi_\B(b_1)\\ &\qquad +\Delta_2(\psi_\Q(t_1^*t_1))
         =0\\
  \iff& \Delta_2((\psi_\B(b_1)+\psi_\Q(t_1))^*(\psi_\B(b_1)+\psi_\Q(t_1)))=0\\
  \iff & \psi_\B(b_1)+\psi_Q(t_1)=0.
\end{split}
\end{align*}
Therefore, the map
$\spn \S_1=\B_1+\spn \Q_1 \ni b_1+t_1\mapsto
\psi_\B(b_1)+\psi_Q(t_1)\in \spn\S_2$ is well-defined, so we obtain a 
$*$-isomorphism $\psi_\S:\spn\S_1\rightarrow\spn\S_2$.

Now let
$\tilde \S_i:=u_{\A_i}(\S_i)\cup \bbC I_{\tilde\A_i}\subseteq
\tilde\A_i$.  Then $\tilde\S_i$ is a MASA skeleton for
$(\tilde\A_i,\tilde\B_i)$ in the sense of
\cite[Definition~3.1]{PittsStReInI} (see
also~\cite[Definition~1.7]{PittsStReInI}).  Note that $\psi_\S$ extends uniquely
to a $*$-isomorphism
$\tilde\psi_\S:\spn \tilde\S_1\rightarrow \spn\tilde\S_2$.  Since
$(\tilde\A_i,\tilde\B_i)$ is a Cartan inclusion, we may argue as in
the proof of \cite[Proposition~5.24]{PittsStReInII} to conclude that
$\tilde\psi_\S$ extends to a regular $*$-isomorphism $\tilde\psi$ of
$(\tilde\A_1,\tilde\B_1)$ onto $(\tilde\A_2,\tilde\B_2)$. Thus
$\psi:=\tilde\psi|_{\A_1}$ is an isomorphism of $\A_1$ onto $\A_2$.
The constructions show $\psi\circ \alpha_1=\alpha_2$.

Suppose $\psi':\A_1\rightarrow \A_2$ is a $*$-isomorphism such that
$\psi'\circ\alpha_1=\alpha_2$.  The uniqueness assertion for $\psi_\B$
shows $\psi'|_{\B_1}=\psi|_{\B_1}$.  Examining each stage of the
construction of $\psi|_\S$, we obtain $\psi'|_{\S_1}=\psi|_{\S_1}$.
As $\spn\S_i$ is dense in $\A_i$, this forces $\psi'=\psi$ and
completes the proof of the uniqueness assertion.
\end{proof}

Before turning to the proof of minimality, we need a bit of
terminology and a fact about abelian \cstaralg s.
Suppose $\A$ is an arbitrary \cstaralg\ and $J\idealin\A$.
Set \[J^\perp:=\{a\in \A:
a j=ja=0 \text{ for all } j\in J\}\] and recall that $J$ is called a \textit{regular
ideal} if $J^\dperp=J$, where $J^\dperp$ means $(J^\perp)^\perp$.  (Because $J$ is an ideal, so is
$J^\perp$.)

Next,
  let $X$ be a locally compact
Hausdorff space
and suppose $J\idealin C_0(X)$.    For an open set $G\subseteq X$, we write $G^\perp$ for the interior of
$X\setminus G$ and  write $\supp(J):=\{x\in X: f(x)\neq 0 \text{
  for some } f\in J\}$.  The proof of the following  is left to
the reader.
\begin{flexstate}{Fact}{}\label{risup}   Let $J\idealin C_0(X)$ and
  put $G:=\supp(J)$.  Then $\supp(J^\dperp)=G^\dperp$, $G$ is a dense
  open subset of $G^\dperp$, and  $(J^\dperp, J)$ has the \iip.
\end{flexstate}

We require the following invariance type result; it is the extension
of~\cite[Proposition~3.14]{PittsStReInI} to our context.
\begin{proposition}[c.f.~{\cite[Proposition~3.14]{PittsStReInI}}]\label{invPsE}
  Suppose $(\C,\D)$ is a regular inclusion having the unique faithful
  pseudo-expectation property, and let $E:\C\rightarrow I(\D)$ be the
  pseudo-expectation. Fix $v\in \N(\C,\D)$, and let $\overline\theta_v$
  be the result of applying Lemma~\ref{AbInEn} to
  $\theta_v:\overline{vv^*\D}\rightarrow \overline{v^*v\D}$.  Then
  for every 
  $x\in \C$, $E(vv^*x)\in\dom\overline\theta_v$ and 
  \begin{equation}\label{invPsE1}
    E(v^*xv)=\overline\theta_v(E(vv^*x)).
  \end{equation} 
\end{proposition}
\begin{proof}
  By the equivalence of statements~\eqref{!pschar1}
  and~\eqref{!pschar2} in Theorem~\ref{!pschar}, $\D^c$ is abelian and
  $(\D^c,\D)$ has the \iip.  Proposition~\ref{PZC22} shows $\iota_1:=E|_{\D^c}$ is a
  $*$-monomorphism; also
   $(I(\D), \iota_1)$ is an injective envelope
  for $\D^c$.  Recall that $\iota_1|_\D=\iota$.

  Let $Q$ and $R$ be the support projections in $I(\D)$
  for $\overline{vv^*\D}\idealin\D$ and $\overline{vv^*\D^c}\idealin\D^c$ respectively.
  Let us show $Q=R$.   
Since $\iota(\overline{vv^*\D})\subseteq
\iota_1(\overline{vv^*\D^c})$, we find $Q\leq R$. 

Observe that $vv^*\in \D^c$ and $e_n:=(vv^*)^{1/n}$ is an  approximate unit for $\overline{vv^*\D^c}$.
If $d\in \D$ and $\iota(d)=(R-Q)\iota(d)$, then as $\iota(d)=\lim_n\iota(e_nd)\in
\iota(\overline{(vv^*)\D})\subseteq QI(\D)$, we obtain $\iota(d)=Q\iota(d)$.  Thus
$(R-Q)\iota(d)=0$, so that 
\[(R-Q)I(\D)\cap \iota(\D)=\{0\}.\]  As $(R-Q)I(\D)\idealin I(\D)$
and $(I(\D), \D,
\iota)$ has the ideal intersection property, we find
$R-Q=0$ as desired.

For $x\in\C$,  Proposition~\ref{PsExBi} yields,
\[E(vv^*x)=\iota_1(vv^*)E(x)\in RI(\D)=QI(\D)=\dom\overline\theta_v.\]
Thus, the right side of \eqref{invPsE1} is defined.

To simplify notation throughout the remainder of the proof, we will
identify $\A$ with its image under $u_\A$ in $\tilde\A$, and we will
write
  \[\M:=\D^c.\]

  We wish to apply \cite[Proposition~3.14]{PittsStReInI} to
  $(\tilde\C,\tilde\M)$. To do so, we  show that
  $(\tilde\C,\tilde\M)$ is a regular and unital MASA inclusion with
  the faithful unique pseudo-expectation property, that
  $(I(\D), \tilde E|_{\tilde\M})$ is an injective envelope for
  $\tilde\M$, and that $\tilde E: \tilde\C\rightarrow I(\D)$ is the
  pseudo-expectation for $(\tilde\C,\tilde\M)$.

  By Lemma~\ref{pseprop}, $(\C,\D)$ is \wnd, so
  Lemma~\ref{sameunit}\eqref{sameunit5} shows $(\tilde\C,\tilde\D)$ is
  a unital inclusion.  Also, Corollary~\ref{interinc} gives
  $\tilde\D\subseteq \tilde\M\subseteq \tilde\C$.  In particular,
  $(\tilde\C,\tilde\M)$ is a unital inclusion.
  Proposition~\ref{f!pse->abelcom} shows $\M$ is abelian,
  Corollary~\ref{heriditary} shows $(\M,\D)$ has the faithful unique
  pseudo-expectation property, and by Proposition~\ref{PZC22},
  $(I(\D), E|_{\M})$ is an injective envelope for $\M$.   Thus,
  $(I(\D),\tilde E|_{\tilde\M})$ is an injective envelope for
    $\tilde\M$.   
  Next,  $(\C,\M)$ is a regular MASA inclusion by
  Lemma~\ref{relcom}\eqref{relcom1}.  Therefore $(\C,\M)$ has the
  AUP by ~\cite[Corollary~2.6]{PittsNoApUnInC*Al}, so
  $(\tilde\C,\tilde\M)$ is a unital, regular MASA inclusion
  (Lemma~\ref{relcom}\eqref{relcom5}).  Corollary~\ref{PsNot1}
  shows $\tilde E$ is the pseudo-expectation for
  $(\tilde\C,\tilde\M)$ and $\tilde E$ is faithful by
  Theorem~\ref{f!ps4rin}.

   Using parts \eqref{relcom1} and \eqref{relcom3} of
 Lemma~\ref{relcom},  we obtain
  \[\N(\C,\D)\subseteq \N(\C,\M)\subseteq \N(\tilde\C, \tilde{\M}).\]   
  Thus by~\cite[Proposition~3.14]{PittsStReInI}, for
  $v\in \N(\C,\D)$ and $x\in \C$,
  \[\tilde E(v^*xv)=\overline\theta_v(\tilde E(vv^*x)), \] that is, $E(v^*xv)=\overline\theta_v(E(vv^*x))$. 
\end{proof}

\begin{proof}[Proof of the Minimality Statement in
  Theorem~\ref{!pschar}]
    Suppose that  $(\A,\B\ms \alpha)$ is a Cartan
package for $(\C,\D)$ and  $\Delta:\A\rightarrow \B$ is the
conditional expectation.
By injectivity of $I(\D)$, there exists a $*$-homomorphism
$\tilde\gamma: \tilde\B\rightarrow I(\D)$ such that
\begin{equation}\label{ush}\iota=\tilde\gamma\circ\tilde\alpha|_{\tilde\D}.
\end{equation}
To show $\tilde\gamma$ is the unique such $*$-homomorphism, suppose
$\tilde\tau: \tilde\B\rightarrow I(\D)$ is another $*$-homomorphism
such that $\iota=\tilde\tau\circ\tilde\alpha|_{\tilde\D}$.   Note that
$\tilde\gamma\circ\tilde\Delta\circ \tilde\alpha|_{\tilde\D}=\iota$, so
$\tilde\gamma\circ\tilde\Delta\circ \tilde\alpha$ is a
pseudo-expectation for $(\tilde\C,\tilde\D)$.   Likewise
$\tilde\tau\circ\tilde\Delta\circ \tilde\alpha$ is a
pseudo-expectation for $(\tilde\C,\tilde\D)$.   Since $(\C,\D)$ has
the unique pseudo-expectation property, so does $(\tilde\C,\tilde\D)$
(Corollary~\ref{PsNot1}).  Therefore,
\[\tilde\gamma\circ\tilde\Delta\circ\tilde\alpha=\tilde\tau\circ\tilde\Delta\circ\tilde\alpha.\]
The definition of Cartan package shows $\B=C^*(\Delta(\alpha(\C)))$,
hence $\tilde\B=C^*(\tilde\Delta(\tilde\alpha(\tilde\C)))$.  Thus
$\tilde\tau=\tilde\gamma$, which shows $\tilde\gamma$ is the unique
$*$-homomorphism satisfying~\eqref{ush}. 

Let $E:\C\rightarrow I(\D)$ be the (unique and faithful) pseudo-expectation.
Since  $\tilde\gamma\circ \tilde\Delta\circ\tilde\alpha$ is the 
pseudo-expectation for $(\tilde\C,\tilde\D)$, Lemma~\ref{prepse} shows
$\tilde E=\tilde\gamma\circ \tilde\Delta\circ\tilde\alpha$ and  $E=\tilde\gamma\circ \tilde\Delta\circ\tilde\alpha\circ
u_\C$.  Define
\[\gamma:=\tilde\gamma\circ u_\B.\] Using~\eqref{stext1}, we find
\begin{equation}\label{ush1}
  E=\tilde\gamma\circ u_\B\circ \Delta\circ\alpha=\gamma\circ\Delta\circ\alpha.
\end{equation}

\begin{flexstate}{\sc Claim}{}\label{fJreg}
Put
 \[\fJ:=\{a\in \A: \Delta(a^*a)\in \ker \gamma\}.\]
Then $\fJ$ is a regular ideal in $\A$ such that $\fJ\cap
\B=\ker\gamma$.
\end{flexstate}
\proof[Proof of Claim~\ref{fJreg}]
To verify Claim~\ref{fJreg}, we  use
\cite[Proposition~3.19(ii)]{BrownFullerPittsReznikoffReIdIdInQu}.  To
do this, we  require the following facts:
\begin{enumerate}
\item\label{fJrega} $\ker\gamma$ is a regular ideal of $\B$; and 
\item \label{fJregb} $\ker\gamma$ is $\alpha(\N(\C,\D))$ invariant, that is,
for every $v\in \N(\C,\D)$,
\begin{equation}\label{!pchar4}
  \alpha(v^*)(\ker \gamma)\alpha(v)\subseteq \ker \gamma.
\end{equation}
\end{enumerate}

Our first step in verifying  \eqref{fJrega} 
is to show
$\ker\gamma$ is the unique ideal in $\B$ maximal with respect to
having trivial intersection with $\alpha(\D)$.
By \cite[Corollary~3.21]{PittsZarikianUnPsExC*In},
 $\ker \tilde \gamma$ is the unique maximal 
 $\tilde\alpha(\tilde\D)$-disjoint ideal of $\tilde\B$.
Suppose $J\idealin\B$ satisfies $J\cap \alpha(\D)=\{0\}$.
Then $u_\B(J)\idealin \tilde\B$.   Let us show
\begin{equation}\label{fJreg1}
  u_\B(J)\cap \tilde\alpha(\tilde\D)=\{0\}.
\end{equation}
We do this by considering two cases.

\vskip 4pt \noindent \textsc{Case 1: $\D$ is unital.}  Since $(\C,\D)$ is \wnd, $\C$ is unital.
Let us show $\B$ is unital and that $e:=\alpha(I_\C)$ is its unit.  Since
$\alpha$ is a regular $*$-monomorphism, $e\in \N(\A,\B)$,
hence $e=e^*e$ commutes with $\B$ by
\cite[Proposition~2.1]{PittsNoApUnInC*Al}.  Thus $e\in \B$ because
$\B$ is a MASA in $\A$.   Let
$n\geq 1$ and $x_1,\dots,x_n\in \C$.  
Then 
\[e\Delta(\alpha(x_1))=\Delta(e\alpha(x_1))=\Delta(\alpha(x_1)); \text{
    likewise } \Delta(\alpha(x_n)) e=\Delta(\alpha(x_n)).\]  Hence  
\[e\left(\prod_{k=1}^n\Delta(\alpha(x_k))\right)=\left(\prod_{k=1}^n\Delta(\alpha(x_k))\right)
  e= \prod_{k=1}^n\Delta(\alpha(x_k)).\]   By the definition of a
package,
the span of finite
products from $\Delta(\C)$ is dense in $\B$.  Thus $e=I_\B$.

As $(\A,\B)$ is a Cartan inclusion, it has  the AUP.  Therefore $I_\B$ is
the unit for $\A$.
The equality~\eqref{fJreg1} now follows  because 
  $u_\B$ is the identity map on $\B$, $\tilde\alpha=\alpha$ and
  $\tilde\D=\D$.  

\vskip 4pt \noindent \textsc{Case 2: $\D$ is not unital.}  Let $y\in u_\B(J)\cap
  \tilde\alpha(\tilde\D)$.  Then for some $(d,\lambda)\in \tilde\D$
  and
  $j\in J$,
  \[y=\alpha(d)+\lambda I_{\tilde\B}=u_\B(j).\]  For any $h\in
  \D$,  multiplying each side of the second equality by $\tilde\alpha(h,0)$ gives 
  \[\alpha(dh)+\lambda\alpha(h)=u_\B(j)\tilde\alpha(h,0).\]  Since the
  left side belongs to $u_\B(\alpha(\D))$ and the right
  belongs to $u_\B(J)$, the hypothesis that $\alpha(\D)\cap J=\{0\}$
  gives $\alpha(dh)+\lambda\alpha(h) =0$.
  Thus for every $h\in \D$,
  \[dh=-\lambda h.\] Let $(e_\lambda)$ be an approximate unit for
  $\D$.  Choosing $h=e_\lambda$ and taking the limit, we obtain
  $\lim
\lambda  e_\lambda$ exists in $\D$ and 
  $d=- \lim \lambda e_\lambda$.  As $\D$ is non-unital, this forces
 $y=0$.   This completes the proof of~\eqref{fJreg1}.

Equality~\ref{fJreg1} and  the fact that $\ker\tilde\gamma$ is the maximal ideal of $\tilde\B$
 disjoint from $\tilde\alpha(\tilde\D)$ gives $u_\B(J)\subseteq \ker
 \tilde \gamma$.  Therefore
 $u_\B(J)=u_\B(J)\cap u_\B(\B)\subseteq \ker\tilde \gamma \cap
 u_\B(\B)=u_\B(\ker \gamma)$, whence 
 $J\subseteq \ker \gamma$.  Thus $\ker\gamma$ is the unique
 maximal $\alpha(\D)$-disjoint ideal of $\B$.

 We are now prepared to show $\ker\gamma$ is a regular ideal of $\B$.
 Suppose $d\in \D$ satisfies $\alpha(d)\in (\ker\gamma)^\dperp$.  Let
 $\I_d$ be the closed ideal of $\B$ generated by $\alpha(d)$.  Then
 $\I_d$ is contained in $(\ker\gamma)^\dperp$, so
 $\alpha(d)\in \I_d\cap (\ker\gamma)^\dperp$.  But
 $\I_d\cap \ker\gamma\subseteq \alpha(\D)\cap \ker\gamma=\{0\}$. 
 Fact~\ref{risup} shows $((\ker\gamma)^\dperp, \ker\gamma)$ has the \iip,
 so $\I_d=\{0\}$.  That $d=0$ follows.  Therefore
 $(\ker\gamma)^\dperp$ is an ideal of $\B$ having trivial intersection
 with $\alpha(\D)$.  As $\ker\gamma\subseteq (\ker\gamma)^\dperp$, the
 maximality property of $\ker\gamma$ gives
 $\ker\gamma=(\ker\gamma)^\dperp$, so $\ker\gamma$ is a regular ideal
 of $\B$.

Turning now to \eqref{fJregb},   we will use the notation
from~\eqref{RnS}.    We claim that for $v\in \N(\C,\D)$,
 \begin{equation}\label{fJreg2a}
   \gamma\circ\theta_{\alpha(v)}=\overline\theta_v\circ \gamma|_{R(\alpha(v))}.
 \end{equation}  The next few paragraphs are devoted to establishing this.

 Fix $v\in \N(\C,\D)$.
 Regularity of the map $\alpha$ and the fact that Cartan inclusions
 have the AUP, gives $\alpha(vv^*)\in \B$.
  Let $\bbP_1\idealin \bbC[t]$ be the  ideal
generated by $t$ in the
ring $\bbC[t]$ of all polynomials with complex coefficients.

Given $p\in\bbP_1$, factor
 $p(t)=tq(t)$, where $q$ is another polynomial.  Let $x\in\C$ and set  
 \[y:=\alpha(p(vv^*))\Delta(\alpha(x))=\alpha(vv^*)\alpha(q(vv^*))\Delta(\alpha(x)).\]
 Using Lemma~\ref{Tv2.2},~\eqref{exinv}, and the fact that $\Delta$ is
 a conditional expectation,
\begin{align*}
   \gamma(\theta_{\alpha(v)}(y))
   &=\gamma(\alpha(v)^*\Delta(\alpha(q(vv^*)x))\alpha(v)) =
\gamma(\Delta(\alpha(v^*q(vv^*)xv))) \\
&\stackrel{\eqref{ush1}}{=}E(v^*q(vv^*)xv) \stackrel{\eqref{invPsE}}{=}
\overline\theta_v(E(vv^*q(vv^*)x))=\overline\theta_v(\gamma(\Delta(\alpha(vv^*q(vv^*)
                                        x))))\notag\\
&=\overline\theta_v(\gamma(y)).\notag
\end{align*}
Let $\bbP(\alpha(vv^*))=\{\alpha(p(vv^*)): p\in \bbP_1\}$.  Then $\bbP(\alpha(vv^*))$ is a
$*$-subalgebra of $\B$.  Keeping $x$ fixed, continuity and the calculation just done
 shows that for every
$k\in \overline{\bbP(\alpha(vv^*))}$,
\[ \gamma(\theta_{\alpha(v)}(k\Delta(\alpha(x)))) = \overline\theta_v
  (\gamma(k\Delta(\alpha(x)))). \]

Now let $x_1, \dots, x_N$ be elements of $\C$ and let
$k\in\overline{\bbP(\alpha(vv^*))}$.   Factor $k=\prod_{j=1}^N k_j$ where $k_j\in
\overline{\bbP(\alpha(vv^*))}$ (when $N>1$, apply
\cite[Proposition~1.4.5]{PedersenC*AlAuGp2018} to obtain $k=u|k|^{1/4}=u(|k|^{1/4(N-1)})^{N-1}$).   Then
\begin{align*} \gamma\left(\theta_{\alpha(v)}\left( k\prod_{j=1}^N
  \Delta(\alpha(x_j))\right)\right)
  &=\gamma\left(\theta_{\alpha(v)}\left(\prod_{j=1}^N
    k_j\Delta(\alpha(x_j))\right)\right) =\overline\theta_v\left(\gamma\left(\prod_{j=1}^N  k_j\Delta(\alpha(x_j)) \right)\right)\\
  & =
    \overline\theta_v\left(\gamma\left(k\prod_{j=1}^N \Delta(\alpha(x_j))\right)\right).
\end{align*}
Since $\B=C^*(\Delta(\alpha(\C)))$, it follows that for $b\in\B$ and
$k\in \overline{\bbP(\alpha(vv^*))}$,
\[\gamma(\theta_{\alpha(v)}(kb))=\overline\theta_v(\gamma(kb)).\]
As $\bbP(\alpha(vv^*))\B$ is dense in $R(\alpha(v))$, we obtain~\eqref{fJreg2a}.

With~\eqref{fJreg2a} in hand, we now complete the proof of the
invariance of $\ker\gamma$.   Suppose $t\in\ker\gamma$, and let $b\in\B$.   Then
\[\gamma(\alpha(v)^*tb\alpha(v))=\gamma(\theta_{\alpha(v)}(tb\alpha(vv^*)))
  \stackrel{\eqref{fJreg2a}}{=}
   \overline\theta_v(\gamma(tb\alpha(vv^*))) =0.\]  By replacing $b$
 with elements from an approximate unit for $\B$, we find
 $\alpha(v^*)t\alpha(v)\in \ker\gamma$, so \eqref{!pchar4} holds.
 This 
 completes the proof of~\eqref{fJregb}.

 Let $N$ be the $*$-semigroup generated by
 $\B\cup\alpha(\N(\C,\D))$. For every $w\in N$, using \eqref{fJregb},
 we find that $w(\ker\gamma)w^*\subseteq \ker\gamma$, that is,
 $\ker\gamma$ is an $N$-invariant regular ideal in $\B$.   While the
 statement of
 \cite[Proposition~3.19(ii)]{BrownFullerPittsReznikoffReIdIdInQu},
 concerns a regular invariant ideals, the proof applies
 for $N$-invariant regular ideals (in fact, the statement
 should have been for $N$-invariant regular ideals).  Thus,
 \cite[Proposition~3.19(ii)]{BrownFullerPittsReznikoffReIdIdInQu} shows
  $\fJ$ is a regular ideal in $\A$ such that $\fJ\cap \B=\ker\gamma$.  
 The proof of Claim~\ref{fJreg} is now complete.
\hfill $\diamondsuit$ \par \vskip 4pt

Let $q: \A\rightarrow \A/\fJ$ be the quotient map. We must show
$(\A/\fJ, \B/(\fJ\cap \B)\ms q\circ\alpha)$ is a Cartan envelope for
$(\C,\D)$.  

First we show 
$q\circ \alpha:(\C,\D)\rightarrow (\A/\fJ,\B/(\J\cap\B))$ is a regular
$*$-monomorphism. 
To show $q\circ\alpha$ is one-to-one, suppose  $x\in \fJ\cap \alpha(\C)$.  Then there exists $c\in \C$ so
that $x=\alpha(c)$.  Since $x^*x\in \fJ\cap\alpha(\C)$,~\eqref{ush1} shows  \[0=\gamma(\Delta(\alpha(c^*c))) =
E(c^*c).\]  Since the pseudo-expectation for $(\C,\D)$ is
faithful, we find $c=0$,  hence $x=0$.  Therefore $q\circ\alpha$ is faithful.
For
$v\in \N(\C,\D)$ and $h\in \B$, \eqref{!pchar4} gives
\[\alpha(v)^* (h+\ker \gamma)\alpha(v)\subseteq
\alpha(v^*)h\alpha(v)+\ker \gamma\in \B/(\B\cap\fJ),\]
whence $q\circ \alpha$ is a regular $*$-monomorphism.

By
\cite[Theorem~4.8]{BrownFullerPittsReznikoffReIdIdInQu},
$(\A/\fJ, \B/(\B\cap\fJ))$ is a Cartan inclusion.
 Clearly  $(\C,\D)$ generates
$(\A/\fJ,\B/(\fJ\cap\B))$ in the sense of part (ii) of
Definition~\ref{defCarEnv}\eqref{defCarEnv1}.

It remains to establish that $(\B/(\fJ\cap \B), \D, q\circ\alpha)$ has
the \iip.  Suppose $K\idealin \B/(\B\cap\fJ)$ satisfies
$K\cap (q(\alpha(\D)))=\{0\}$.  Since $q\circ\alpha$ is faithful, 
 $\B\cap q^{-1}(K)$ is an ideal of
$\B$ having trivial intersection with $\alpha(\D)$.  By the maximality
property of $\ker\gamma$, $\B\cap q^{-1}(K)\subseteq \ker\gamma$, so
$K=\{0\}$.  
Therefore $(\B/(\B\cap \fJ), \D, q\circ\alpha)$ has the \iip.  Thus
$(\A/\fJ, \B/(\B\cap\fJ)\ms  q\circ\alpha)$ is a Cartan envelope for
$(\C,\D)$.  This finishes the proof of the minimality statement and
also concludes the  proof of Theorem~\ref{!pschar}.
\end{proof}



\mysec[Envelopes for Pseudo-Cartan Inclusions]{\Pd s and their Cartan Envelopes}\label{Sec: pd}

\numberwithin{equation}{subsection}

In this \LCsection\  we give the formal definition of a \pd.  The key
results of this \LCsection\  are: Theorem~\ref{props}, which shows that in the
presence of the AUP, the
unitization process commutes with taking Cartan envelopes, and
Theorem~\ref{simple} which gives 
 some properties shared by a \pd\ and its Cartan envelope.
 Proposition~\ref{pdcons} gives a construction which can be used to
 produce a family of \pd s starting with a Cartan
inclusion. 

\subsection[Definitions \& Examples]{Definition of \Pd s and Examples}\label{defPD+E}

We are now prepared to define the notion of  \pd.  Standing
Assumption~\ref{SA1} remains in force throughout, that is, we assume
that $\D$ is abelian for every inclusion $(\C,\D)$.
      
\begin{definition}\label{defPsDi}  A regular inclusion $(\C,\D)$ is a
  \textit{\pd}\index{Inclusion!Pseudo-Cartan} if it satisfies the following:
  \begin{enumerate}
  \item $\D^c$ is abelian; and 
  \item both $(\C,\D^c)$ and $(\D^c,\D)$ have the \iip. 
  \end{enumerate}
\end{definition}

We will frequently use the alternate descriptions of a \pd\ found in
Theorem~\ref{!pschar}, often without comment.

The following description of abelian \pd s is immediate from Definition~\ref{defPsDi}.
\begin{flexstate}{Observation}{}\label{abelpd} Let $(\C,\D)$ be an
  inclusion with $\C$ abelian.  Then $\D^c=\C$,
  so $(\C,\D)$
  is a \pd\ if and only if $(\C,\D)$ is regular and has the \iip.
\end{flexstate}

Let us compare the notion of \pd\ with  other related classes of inclusions.
\begin{description}
\item[Cartan Inclusions]  
Proposition~\ref{Car+upse} shows that every Cartan inclusion is a \pd.
Furthermore, Proposition~\ref{Car+upse} may be interpreted as stating
that a \pd\ $(\C,\D)$ is a Cartan inclusion if and only if the
pseudo-expectation $E:\C\rightarrow I(\D)$ ``is'' a conditional
expectation in the sense that $E(\C)=\iota(\D)$.
 This characterization is the reason we  chose  the term
``\pd'' in Definition~\ref{defPsDi}.

\item[Virtual Cartan Inclusions]
  We will use the term
  \textit{virtual Cartan inclusion}\index{Inclusion!virtual Cartan} for a regular MASA inclusion
  $(\C,\D)$ having the \iip.  
Theorem~\ref{regmasaps} shows that every virtual Cartan inclusion is a \pd;
  furthermore, the \pd\ $(\C,\D)$ is a virtual Cartan inclusion if and
  only if $\D$ is a MASA in $\C$.
By~\cite[Theorem~2.6]{PittsNoApUnInC*Al}, every virtual Cartan
inclusion  has the AUP. 
\item[Weak Cartan Inclusions]
    Weak Cartan inclusions were defined
  in~\cite[Definition~2.11.5]{ExelPittsChGrC*AlNoHaEtGr}.  We shall
  show in Proposition~\ref{wC->pd} below that 
  that every weak Cartan inclusion is a \pd.  Since the axioms for weak
  Cartan inclusions include the AUP, the class of \pd s properly contains
  the class of weak Cartan inclusions.

  The example given in~\cite[Remark~2.11.6]{ExelPittsChGrC*AlNoHaEtGr}
  shows that there exists a Cartan inclusion $(\C,\D)$ (as in
  Definition~\ref{innore}\eqref{innore6} above) which is not a weak
  Cartan inclusion, so the classes of \pd s with the AUP and weak
  Cartan inclusions are also distinct.  Also, an example of a unital
  \pd\ $(\C,\D)$ with $\C$ abelian which is not a weak Cartan
  inclusion is given in the paragraph just prior
  to~\cite[2.11.16]{ExelPittsChGrC*AlNoHaEtGr}.  However, in both
  these examples, the algebra $\C$ is not separable.  We do not know
  whether every \pd\ $(\C,\D)$ with the AUP and such that $\C$ is
  separable is a weak Cartan inclusion.

 \item[Pseudo-Diagonals and Abelian Cores]  These are classes of
   regular inclusions which are not assumed to be regular.
   Nagy and Reznikoff define the notion of a pseudo-diagonal,
    see~\cite[p.~268]{NagyReznikoffPsDiUnTh}.  One of the
    requirements for the inclusion $(\C,\D)$ to be a pseudo-diagonal
    is that there exists a faithful conditional expectation
    $E:\C\rightarrow \D$, however regularity of $(\C,\D)$ is
    \text{not} required.   Nagy and Reznikoff observe that when
    $(\C,\D)$ is a pseudo-diagonal, $\D$ is necessarily a MASA in
    $\C$,~\cite[Corollary~3.2]{NagyReznikoffPsDiUnTh}.  Thus,  when
    $(\C,\D)$ is a \textit{regular} pseudo-diagonal, $(\C,\D)$ is a
    Cartan inclusion and so is a \pd.   Likewise, every
    \textit{regular} abelian core (see~\cite[p.\
    272]{NagyReznikoffPsDiUnTh}) is a Cartan inclusion and hence a \pd.

\end{description}

For a unital \pd\ $(\C,\D)$, the relative commutant $\D^c$ is 
a MASA in $\C$ and~\cite[Lemma~2.10]{PittsStReInI} shows $(\C,\D^c)$
is a regular inclusion.   As  $(\C,\D^c)$ has the \iip\
(Definition~\ref{defPsDi}(b)), it  is a virtual Cartan inclusion.   We now observe
  the same is true for possibly non-unital \pd s.
\begin{flexstate}{Observation}{}\label{vcfrompc}   If $(\C,\D)$ is a \pd,
  then $(\C,\D^c)$ is a virtual Cartan inclusion.
\end{flexstate}
\begin{proof}
Lemma~\ref{relcom}\eqref{relcom1} shows that $(\C,\D^c)$ is a regular
inclusion, so it is a regular MASA inclusion with the \iip.
\end{proof}

Observation~\ref{vcfrompc} gives the following.
\begin{flexstate}{Fact}{}\label{nopsCar0} Suppose $\C$ is a \cstaralg\ 
  such that no MASA in  $\C$ is
  regular.   Then  $\C$ cannot contain a pseudo-Cartan subalgebra.
\end{flexstate}

A viewpoint found in~\cite{PittsZarikianUnPsExC*In} is that the unique
pseudo-expectation property for an inclusion $(\C,\D)$ is a
``largeness property''.  From this perspective, when no MASA in $\C$
is regular,
no MASA in $\C$ is large enough to have a unique pseudo-expectation.

\begin{example}\label{nopsCar}
Let $\H$ be a separable and infinite
  dimensional Hilbert space.
  By~\cite[Theorem~3.8]{PittsNoApUnInC*Al}, $\bh$ (regarded as a \cstaralg) does
  not contain a regular MASA, so $\bh$ does not contain a pseudo-Cartan
  subalgebra.  
\end{example}

Our next goal is to show that every weak Cartan inclusion is a \pd.           
 Before proceeding,
we recall  some terminology
from~\cite{ExelPittsChGrC*AlNoHaEtGr}.

\begin{definition}[{\cite[Definitions 2.7.2, 2.8.1 and
    2.11.5]{ExelPittsChGrC*AlNoHaEtGr}}] \label{fdef}
  Let $(\C,\D)$ be an inclusion and $\sigma\in \hat\D$.
  \begin{enumerate}
  \item \label{fdef1}For $x\in\C$, we say $\sigma$ is
  \textit{free relative to $x$}\index{Free point!relative to $x$} if whenever $f_1, f_2$ are states on
  $\C$ such that $f_1|_\D=\sigma=f_2|_\D$, we have $f_1(x)=f_2(x)$.
  \item \label{fdef2} We say
  $\sigma\in \hat\D$ is \textit{free for $(\C,\D)$,}\index{Free point} or 
 a  \textit{free point} when the context is clear, if $\sigma$ uniquely
  extends to a state $f$ on $\C$.
\item\label{fdef3} The inclusion $(\C,\D)$ is called
  \textit{topologically free}\index{Inclusion!topologically free} if
  the set of free points for $(\C,\D)$ is dense in $\hat\D$.
  Other terminology for topologically free inclusions may be found in
  the literature, for example,  the \textit{almost extension
    property}  is used in
  \cite{KwasniewskiMeyerApAlExPrUnPsEx,NagyReznikoffPsDiUnTh,ZarikianPuExPrDiCrPr}
  instead of `topologically free
  inclusion'.
\end{enumerate}
\end{definition}

\begin{proposition}\label{wC->pd} Let $(\C,\D)$ be a weak Cartan
  inclusion.  Then $(\C,\D)$ is a \pd.
\end{proposition}
\begin{proof}
  First note that the weak Cartan inclusions considered in
  \cite{ExelPittsChGrC*AlNoHaEtGr} have the AUP and are regular
  (see~\cite[Definition~2.3.1, Hypothesis~2.3.2, and
  Definition~2.11.5]{ExelPittsChGrC*AlNoHaEtGr}).  Next,
  \cite[Proposition~2.11.7]{ExelPittsChGrC*AlNoHaEtGr} shows that
  $(\C,\D)$ is topologically free.   Combining
  \cite[Theorems~5.5 and~3.6]{KwasniewskiMeyerApAlExPrUnPsEx}, we find
  that $(\C,\D)$ has the unique pseudo-expectation property; let
  $E:\C\rightarrow I(\D)$ be its pseudo-expectation.   We must show
  $E$ is faithful.

  Since
  $(\C,\D)$ has the AUP, $(\tilde\C,\tilde\D)$ is a regular inclusion
  by Lemma~\ref{relcom}\eqref{relcom5}.
  Corollary~\ref{PsNot1} shows $\tilde E$ is the unique
  pseudo-expectation for $(\tilde\C,\tilde\D)$.  Let $\L(\tilde\C,\tilde\D)$ be the
  left kernel of $\tilde E$.  It
  follows from \cite[Theorem~6.5]{PittsStReInII} that
  $\L(\tilde\C,\tilde\D)$ is the largest ideal of $\tilde\C$ having
  trivial intersection with $\tilde\D$.  

  If $J\idealin \C$ has trivial intersection with $\D$, then $J$ is
  also an ideal of $\tilde\C$ which has trivial intersection with
  $\tilde\D$.  Thus $J\subseteq \L(\tilde\C,\tilde\D)$.  Hence $\L_1:=\C\cap \L(\tilde\C,\tilde\D)$ is the largest ideal of
  $\C$ having trivial intersection with $\D$.  Therefore,~\cite[Proposition~2.11.3]{ExelPittsChGrC*AlNoHaEtGr} shows
  that $\L_1$ is the grey ideal for $(\C,\D)$.  The grey ideal in a
weak Cartan inclusion is $\{0\}$ (see \cite[Definition~2.11.5]{ExelPittsChGrC*AlNoHaEtGr}), so $\L_1=\{0\}$.
But $\L_1$ is the left
  kernel of $E$, whence $(\C,\D)$ has the faithful unique
  pseudo-expectation property.  Theorem~\ref{!pschar} then shows
  $(\C,\D)$ is a \pd.
\end{proof}

\subsection[Unitization of Cartan Envelopes]{The AUP and Unitization of Cartan Envelopes}
Given a \wnd\ regular inclusion $(\C,\D)$, $(\tilde\C,\tilde\D)$ is a
unital inclusion by Lemma~\ref{sameunit}\eqref{sameunit5}, but we have
been unable to establish whether $(\tilde\C,\tilde\D)$ is regular.  Two conditions which ensure that regularity of  $(\C,\D)$
implies $(\tilde\C,\tilde\D)$ is also regular are: (i) $(\C,\D)$ has the AUP
(Corollary~\ref{uCnr}); and (ii) $\C$ is abelian 
(because every unitary in $\tilde\C$ normalizes $\tilde\D$ and
$\tilde\C$ is spanned by its unitaries).   On
the other hand, \cite[Example~3.1]{PittsNoApUnInC*Al} shows that
regularity of $(\tilde\C,\tilde\D)$ does not imply regularity of
$(\C,\D)$. 
Theorem~\ref{f!ps4rin} (and Theorem~\ref{!pschar}) give the following.
\begin{flexstate}{Observation}{}\label{unpdiffnupd}   Suppose both
  $(\C,\D)$ and $(\tilde\C,\tilde\D)$ are regular
  inclusions.   Then $(\C,\D)$ is a \pd\ if and only if
  $(\tilde\C,\tilde\D)$ is a \pd.
\end{flexstate}

A \pd\ $(\C,\D)$ may not have the AUP.  Our next two results show that
it is possible to enlarge $\D$ to obtain a \pd\ $(\C,\E)$ having 
the AUP such that the Cartan envelopes of $(\C,\D)$ and $(\C,\E)$
are the same.    
The significance of these results is that it is often possible to
replace $(\C,\D)$ with $(\C,\E)$ when the AUP is required.
The first of these results shows that $\E$ may be taken to be $\D^c$;
it 
extends~\cite[Proposition~5.29(b)]{PittsStReInII} to include
not-necessarily unital inclusions. 

\begin{proposition}\label{CE4rcom}   Suppose $(\C,\D)$ is a \pd\ and
  consider the virtual Cartan inclusion $(\C,\D^c)$.
  Then $(\A,\B\ms \alpha)$ is a Cartan
  envelope for $(\C,\D)$ if and only if  $(\A,\B\ms  \alpha)$ is Cartan envelope
  for $(\C,\D^c)$.
\end{proposition}
\begin{proof}
Proposition~\ref{rcie} shows that if $(\A,\B\ms \alpha)$ is a Cartan
envelope for $(\C,\D^c)$, then $(\A,\B\ms \alpha)$ is a Cartan envelope
for $(\C,\D)$.

Suppose $(\A,\B\ms \alpha)$ is a Cartan envelope for $(\C,\D)$.  As
$(\C,\D^c)$ is a virtual Cartan inclusion, 
we may find a Cartan envelope $(\A',\B'\ms \alpha')$ for $(\C,\D^c)$,
which by Proposition~\ref{rcie}, is  also a Cartan envelope for $(\C,\D)$. By
uniqueness of Cartan envelopes for $(\C,\D)$, there exists a unique regular
$*$-isomorphism $\psi: (\A',\B')\rightarrow (\A,\B)$ such that
$\psi\circ\alpha'=\alpha$.   Since
$\alpha':(\C,\D^c) \rightarrow (\A',\B')$ is regular, we conclude 
$\alpha: (\C,\D^c)\rightarrow (\A,\B)$ is regular.  Noting that
$\alpha(\D)\subseteq \alpha(\D^c)\subseteq \B$,
Lemma~\ref{interiip} shows $(\B, \alpha(\D^c))$ has the \iip.  Thus 
$(\A,\B\ms \alpha)$ is a Cartan envelope for $(\C,\D^c)$.  
\end{proof}

Given the  \pd\ $(\C,\D)$, 
examples show that $\D^c$ need not be
contained in $\N(\C,\D)$.
Nevertheless,
by suitably enlarging $\D$, we can use $(\C,\D)$ to produce \pd s
$(\C,\E)$ having the AUP such that $\E\subseteq \N(\C,\D)$ and which
have the same Cartan envelope as $(\C,\D)$.  However, unlike $\D^c$, 
$\E$ may not be a  MASA in $\C$.
\begin{proposition}\label{pdAUP} Let $(\C,\D)$ be a \pd.   Suppose 
  $\P\subseteq\N(\C,\D)$ is a $*$-semigroup such that $\D\subseteq
  \P$,  and $\spn\P$ is dense in $\C$.  Let
  \[\E:=C^*(\{p\in\P: p\geq 0\}).   \]  Then $(\C,\E)$ is a \pd\ having the
  AUP.  Furthermore,
  \begin{enumerate}
  \item \label{pdAUP1} $\D\subseteq \E\subseteq\D^c$;
  \item \label{pdAUP2}$\E\subseteq \N(\C,\D)$; and
    \item \label{pdAUP3} $(\A,\B\ms \alpha)$ is a Cartan envelope for $(\C,\D)$ if and
      only if $(\A,\B\ms  \alpha)$ is  a Cartan envelope for $(\C,\E)$.
    \end{enumerate}
  \end{proposition}
  \begin{proof} 
We first establish parts \eqref{pdAUP1} and \eqref{pdAUP2}, then show
$(\C,\E)$ is a \pd\ with the AUP, and finally show
part~\eqref{pdAUP3}.

The first inclusion in \eqref{pdAUP1} follows from the fact that
$\D\subseteq \P$ and $\D$ is the span of its positive elements.  For the
second,~\cite[Proposition~2.1]{PittsNoApUnInC*Al} shows that for
$0\leq p\in
\P$, $p^2\in\D^c$, so $p\in \D^c$; thus $\E\subseteq \D^c$.  This
gives \eqref{pdAUP1}.

Before moving to part \eqref{pdAUP2}, we introduce some notation.  Let us write,
    \[\P^+:=\{p\in\P: p\geq 0\}.\] Since $\D^c$ is
    abelian, $\P^+$ is a $*$-semigroup, whence $\E=\overline\spn\, \P^+$.

We now turn to part~\eqref{pdAUP2}.  For $p\in \P^+$ and 
and $d\in \D$, 
\cite[Proposition~2.1]{PittsNoApUnInC*Al} gives $p^2d^*d=d^*dp^2\in
\D$.  Taking square roots, we obtain $p|d|=|d|p\in \D$.  Since $\D$ is
the span of its positive elements, we see that for every $d\in \D$ and
$p\in \P^+$,
\[pd=dp\in \D.\]
Thus,
if $x\in \spn\P^+$, $xd\in \D$, so $x^*dx=xdx^*=x^*xd\in \D$ also.  Therefore, $\spn \P^+\subseteq \N(\C,\D)$.  As
$\N(\C,\D)$ is closed, part~\eqref{pdAUP2} holds.

We next show $(\C,\E)$ is a regular inclusion.  For any $v\in \P$,
$v^*\P^+ v\cup v\P^+ v^*\subseteq \P^+$, whence
$v\E v^*\cup v^*\E v\subseteq \E$.  Therefore,
$\P\subseteq \N(\C,\E)$, and because $\C=\overline\spn\, \P$,
$(\C,\E)$ is a regular inclusion. 

Note that $\E^c=\D^c$.   Thus, part~\eqref{pdAUP1} and
Lemma~\ref{interiip} show
$(\E^c,\E)$ has the \iip.   As $(\C,\D)$ is a \pd, $(\C,\E^c)$ has the \iip,
so  $(\C,\E)$
is a \pd.

To show $(\C,\E)$ has the AUP, let $(u_\lambda)$ be an approximate
unit for $\E$.  Let $\eps>0$.  For $v\in \P$, we may choose
$n\in \bbN$ so that $\norm{v(v^*v)^{1/n} -v}<\eps/3$.  Since
$(v^*v)^{1/n}\in \E$, there exists $\lambda_0$ so that for
$\lambda\geq \lambda_0$,
$\norm{v}\norm{(v^*v)^{1/n}u_\lambda-(v^*v)^{1/n}}<\eps/3.$  For
$\lambda\geq\lambda_0$, 
\begin{align*}\norm{vu_\lambda -v}&\leq
                                    \norm{vu_\lambda-v(v^*v)^{1/n}u_\lambda}+
                                    \norm{v(v^*v)^{1/n}u_\lambda -
                                    v(vv^*)^{1/n}}+\norm{v(v^*v)^{1/n}-v}
  \\
  &\leq
    2\norm{v-v(v^*v)^{1/n}}+\norm{v}\norm{(v^*v)^{1/n}u_\lambda-(v^*v)^{1/n}} <\eps.
\end{align*} Therefore $vu_\lambda \rightarrow v$; similarly
$u_\lambda v\rightarrow v$.  Since $\spn\P$ is dense in $\C$, it
follows that $(\C,\E)$ has the AUP.

Finally, using Proposition~\ref{CE4rcom} and the equality $\E^c=\D^c$, 
\begin{align*}
  (\A,\B\ms \alpha)\text{ is a Cartan envelope for } (\C,\D) \iff &
                                                              (\A,\B\ms \alpha) \text{ is a Cartan envelope for }\\
  &(\C,\D^c)=(\C,\E^c)
  \\
  \iff&
        (\A,\B\ms \alpha)\text{ is a Cartan envelope for }\\
  &(\C,\E). 
\end{align*} 
\end{proof}

As noted in Corollary~\ref{uCnr}, for the \wnd\ and non-unital
inclusion $(\C,\D)$, the inclusion mapping $u_\C$ of $(\C,\D)$ into
$(\tilde\C,\tilde\D)$ is regular if and only if $(\C,\D)$ has the AUP.
We now study the relationships between the Cartan envelopes of
$(\C,\D)$ and $(\tilde\C,\tilde\D)$ when $(\C,\D)$ has the AUP.  We
begin with a lemma concerning the image under $\alpha$ of an approximate unit for
$\D$ in a Cartan package $(\A,\B\ms \alpha)$.
\begin{lemma}\label{auinC1}  Suppose $(\C,\D)$ is an inclusion with
  the AUP  and
  $(\A,\B\ms \alpha)$ is a Cartan package for $(\C,\D)$.     If
  $(u_\lambda)\subseteq \D$ is an approximate unit for $\C$, then
 $(\alpha(u_\lambda))\subseteq\B$ is an approximate unit for
 $\A$.
\end{lemma}
\begin{proof}
   Let
$\Delta:\A\rightarrow \B$ be the conditional expectation.    For
$x\in \C$, we have
\[\lim_\lambda\Delta(\alpha(x))\alpha(u_\lambda)
  =\lim_\lambda\Delta(\alpha(xu_\lambda))=\Delta(\alpha(x)).\]  Since
$\B$ is generated by $\Delta(\alpha(\C))$, it follows that $(\alpha(u_\lambda))$
is an approximate unit for $\B$.  As $(\A,\B)$ is a regular
MASA inclusion, an application
of~\cite[Theorem~2.6]{PittsNoApUnInC*Al} shows $(\alpha(u_\lambda))$
is an approximate unit for $\A$.
\end{proof}

\begin{lemma} \label{RCaPk} Let $(\C,\D)$ be a regular inclusion with
  the AUP such that $\C$ is not unital.  Suppose
   $(\A,\B\ms \alpha)$ is a Cartan package (resp.\ Cartan
  envelope) for $(\C,\D)$.
  Then $(\tilde\A,\tilde\B\ms \tilde\alpha)$ is a Cartan package
  (resp.\ Cartan envelope) for
  $(\tilde\C,\tilde\D)$.    
\end{lemma}
\begin{proof}
  We will use $\Delta: \A\rightarrow \B$ for the conditional
  expectation.  Since $(\A,\B)$ is a Cartan pair, so is
  $(\tilde\A,\tilde\B)$; also,
  $\tilde\Delta: \tilde\A\rightarrow \tilde\B$ is the conditional
  expectation of $\tilde\A$ onto $\tilde\B$.
  
  That $\tilde\A$ is generated by
  $\tilde\alpha(\tilde\C)\cup\tilde\Delta(\tilde\alpha(\tilde\C))$
  follows from the fact that $\A$ is generated by
  $\alpha(\C)\cup\Delta(\alpha(\C))$.  Similarly, $\tilde\B$ is generated
  by $\tilde\Delta(\pi(\tilde\C))$.  It remains to show that $\tilde\alpha$
  is a regular $*$-homomorphism, that is,
  $\tilde\alpha(\N(\tilde\C,\tilde\D))\subseteq \N(\tilde\A,\tilde\B)$.  For this, we use
  the relative strict topology on $\tilde\A$ , which we now describe.

For each $x\in \A$, the map $\tilde\A\ni
a\mapsto\norm{xa} +\norm{ax}$ is a seminorm on $\tilde\A$, and the
\textit{relative strict} topology on $\tilde\A$
is the smallest topology on $\tilde\A$ making each of these seminorms continuous.
         Thus, a net $(a_\lambda)$ in $\tilde\A$ converges to $a\in \tilde\A$ in the
        relative strict topology if and only if for every $x\in \A$,
        \[\lim_\lambda
          (\norm{(a_\lambda-a)x}+\norm{(x(a_\lambda-a)})=0.\] We write
        $a=^{rs}\lim_\lambda a_\lambda$ when the net $(a_\lambda)$ in
        $\tilde\A$ converges in the relative strict topology to $a\in \tilde\A$.
        The relative strict topology on $\tilde\A$ is locally convex, and
        since $\A$ is an essential ideal in $\tilde\A$, it is Hausdorff.
        Routine arguments show that when $\tilde\A$ is equipped with the
        relative strict topology, the adjoint operation is continuous
        and multiplication is jointly
        continuous on norm-bounded subsets of $\tilde\A$.

We claim that if $(w_\lambda)$ is a bounded net in $\N(\tilde\A,\tilde\B)$ which
converges relative strictly to $w\in \tilde\A$, then $w\in \N(\tilde\A,\tilde\B)$.
Indeed, for $h, k\in \tilde\B$, as $w_\lambda^*hw_\lambda\in \tilde\B$, 
\[ w^*hwk={}^{rs}\lim_\lambda w_\lambda^* hw_\lambda k
  ={}^{rs}\lim_\lambda kw_\lambda^* h w_\lambda =k w^*hw.\] As $\tilde\B$
  is a MASA in $\tilde\A$, $w^*hw\in \tilde\B$.  Likewise $whw^*\in \tilde\B$, so our
  claim holds.  

  Now let $(u_\lambda)$ be a net in $\D$ which is an approximate unit for $\C$.
  Choose $v\in \N(\tilde\C,\tilde\D)$ and write $v=(x,\xi)$
for some $x\in \C$ and $\xi\in \bbC$.  For each $\lambda$, $u_\lambda\in
\N(\C,\D)$, so (since $\N(\tilde\C,\tilde\D)$ is closed under
multiplication) $vu_\lambda=(xu_\lambda +\xi u_\lambda,0)\in \N(\C,\D)$.  The regularity of the map
$\alpha$ and the fact that $\N(\A,\B)\subseteq \N(\tilde\A,\tilde\B)$
(Lemma~\ref{relcom}(c))  gives
\[\N(\A,\B)\ni \alpha(vu_\lambda)=\tilde\alpha(vu_\lambda)\in \N(\tilde\A,\tilde\B).\]
By Lemma~\ref{auinC1}, $\alpha(u_\lambda)$ is an approximate unit for
$\A$.   Therefore, $\alpha(u_\lambda)$ converges in the relative
strict topology to $I_{\tilde\A}$.   As $\tilde\alpha(vu_\lambda)$ is a bounded net
in $\N(\tilde\A,\tilde\B)$, the claim gives 
\[\tilde\alpha(v)=\tilde\alpha(v)\,
  {}^{rs}\lim_\lambda\alpha(u_\lambda)
  ={}^{rs}\lim_\lambda\tilde\alpha(vu_\lambda)\in \N(\tilde\A,\tilde\B).\]  Thus $\tilde\alpha$ is a
regular $*$-monomorphism.  This completes the proof that $(\tilde\A,\tilde\B\ms \tilde\alpha)$
is a Cartan package for $(\tilde\C,\tilde\D)$.

Now suppose that $(\A,\B\ms \alpha)$ is a Cartan envelope for
$(\C,\D)$.  It remains to establish part~\eqref{defCarEnv2} of
Definition~\ref{defCarEnv}.  By hypothesis,  $(\B,\alpha(\D))$ has the
\iip, so $(\tilde\B,\tilde\alpha(\tilde\D))=(\tilde\B, (\alpha(\D))^\sim)$ has the
\iip\ by Lemma~\ref{iipiffuiip}.   
\end{proof}

We will use the following in the proof of Lemma~\ref{nunfrun}, however it is
interesting in its own right.
\begin{proposition}\label{simple3.0} Let $(\C,\D)$ be a \pd\ with Cartan
      envelope $(\A,\B\ms \alpha)$.   Then $\C$ is unital if and only if $\A$ is unital.
    \end{proposition}
    \begin{proof}
Suppose $\C$ is unital.  Lemma~\ref{cewnd} shows $\A$ is unital and
  $\alpha(I_\C)=I_\A$.

Now suppose $\A$ is unital.     Proposition~\ref{CE4rcom} shows that $(\A,\B\ms \alpha)$ is also a
  Cartan envelope for $(\C,\D^c)$.  Replacing $\D$ with
  $\D^c$ if necessary, we may therefore assume
  $(\C,\D)$ has the AUP.
  
 Let
$(u_\lambda)\subseteq \D$ be an approximate unit for $\C$.
Choose $x_1,\dots, x_n\in \C$ and let 
\[z=\prod_{k=1}^n \Delta(\alpha(x_k)).\] Since
$\Delta(\alpha(x_n))\, \alpha(u_\lambda)= \Delta(\alpha(x_nu_\lambda))$,
$z\alpha(u_\lambda)\rightarrow z$; likewise,
$\alpha(u_\lambda) z\rightarrow z$.  As the collection of
all finite products of $\Delta(\alpha(\C))$ has  dense span in $\B$, it
follows that $\alpha(u_\lambda)$ is an approximate unit for $\B$, and
hence also for $\A$ (because $(\A,\B)$ has the AUP).
Therefore, 
\[\norm{\alpha(u_\lambda) I_\A-I_\A}=\norm{\alpha(u_\lambda)- I_\A}\rightarrow 0.\]  Hence $(u_\lambda)$
converges to an element $e\in \C$.  Since $\alpha(e)=I_\A$, we find
$e=I_\C$.
\end{proof}

The following result shows that in some settings,  a Cartan envelope for
a non-unital inclusion may be constructed from a Cartan envelope for its 
unitization.
\begin{lemma}\label{nunfrun}  Let $(\C,\D)$ be a regular inclusion
  having the AUP with $\C$ not unital.  Suppose $(\A,\B\ms \pi)$ is
  a Cartan envelope for $(\tilde\C,\tilde\D)$ and write
  $\Delta$ for the conditional expectation of $\A$ onto
    $\B$.
 Put
  \[
    \alpha:=\pi|_\C,\quad
    \D_1:=C^*(\Delta(\alpha(\C))), \dstext{and}
    \C_1:=C^*(\alpha(\C)\cup \D_1).
    \]
    Then $(\C_1,\D_1\ms \alpha)$ is a Cartan envelope
    for $(\C,\D)$; further, $\C_1\idealin \A$ and $\D_1\idealin
    \B$ are ideals having codimension one. 
  \end{lemma}
  \begin{proof}
Note that since $(\B,\tilde\D,\pi)$ and $(\tilde\D,\D)$ have the \iip,
 $(\B,\D,\pi)$ also has the \iip.
Proposition~\ref{reg+Del} 
shows $(\C_1,\D_1\ms \alpha)$
  is a Cartan envelope for $(\C,\D)$.  Since $\C$ was assumed
  non-unital, $\C_1$ is also non-unital by Proposition~\ref{simple3.0}.

  Let $I_\A$ denote the unit of $\A$ and consider the \cstaralg s
    \[\D_1^+= \D_1+\bbC I_\A 
      \dstext{and} \C_1^+=\C_1+\bbC I_\A.\]
    Obviously, $\D_1^+\subseteq \B$, and $\C_1^+\subseteq \A$;
Lemma~\ref{cewnd} gives
    $\pi(0,1)=I_\A$.  Since $(\A,\B\ms \alpha)$ is a Cartan envelope for
$(\tilde\C,\tilde\D)$, 
    $\B$ is
    generated by $\pi(\Delta(\tilde\C))=\pi(\Delta(\C))+\bbC I_\A$ and $\A$ is generated by
    $\B\cup \pi(\tilde\C)=\B\cup (\pi(\C)+\bbC I_\A)$. 
    $\C_1^+=\A$.  Thus,
    \begin{equation}\label{ABunt}
      (\tilde\C_1,\tilde\D_1)\simeq (\C_1^+,\D_1^+)=(\A,\B).
    \end{equation}
That
$\C_1\idealin \A$ and $\D_1\idealin \B$ are ideals of codimension one
follows from~\eqref{ABunt}.
\end{proof}

\begin{theorem} \label{props} Let  $(\C,\D)$ be a regular inclusion with the
  AUP such that $\C$ is not unital.   Then $(\C,\D)$ is a \pd\ if and
  only if $(\tilde\C,\tilde\D)$ is a \pd. 
  Moreover, the following statements hold.
\begin{enumerate}
  \item \label{p1} If $(\A,\B\ms \alpha)$ is a Cartan envelope for
    $(\C,\D)$, then
    $(\tilde\A, \tilde\B\ms  \tilde\alpha)$ is
  a
    Cartan envelope for $(\tilde\C,\tilde\D)$.
  \item \label{p2} Suppose $(\A,\B\ms \pi)$ is a Cartan envelope for
      $(\tilde\C,\tilde\D)$ and $\Delta:\A\rightarrow\B$ is the
      conditional expectation.  Let \[\alpha:=\pi|_{\C},\quad
        \D_1=C^*(\Delta(\pi(\C))), \dstext{and}
      \C_1=C^*(\pi(\C)\cup
      \D_1).\]
      Then $\D_1\idealin \B$ and $\C_1\idealin \A$ are ideals with codimension 1
      and $(\C_1,\D_1\ms \alpha)$ is a Cartan envelope for $(\C,\D)$.
\end{enumerate}
\end{theorem}
\begin{proof}
That $(\C,\D)$ is a \pd\ if and only if $(\tilde\C,\tilde\D)$ is a
\pd\ follows by combining Lemma~\ref{relcom}\eqref{relcom5} and
Observation~\ref{unpdiffnupd}.
Lemma~\ref{RCaPk} gives part \eqref{p1}, and
Lemma~\ref{nunfrun} gives part \eqref{p2}. 
\end{proof}

\subsection[Shared Properties]{Some Properties Shared by a Pseudo-Cartan Inclusion and
  its Cartan Envelope}
This subsection is devoted to establishing Theorem~\ref{simple}, which
shows that certain 
desirable properties are common to both a
 \pd\ and its Cartan envelope.  
 
Recall that
Definition~\ref{defextinc}\eqref{defextinc2} states that
$(\C_1,\D_1\ms \alpha)$ is an essential expansion of $(\C,\D)$ if the inclusion
$(\D_1,\D,\alpha|_{\D})$ has the \iip.
We next observe that more can be said for essential expansions 
in the context of \pd s.
\begin{flexstate}{Observation}{} \label{Cext} For $i=1,2$, suppose
  $(\C_i,\D_i)$ are \pd s, and $(\C_2,\D_2\ms \alpha)$ is an essential expansion of
  $(\C_1,\D_1)$.  Then the inclusion $(\C_2,\C_1,\alpha)$ has the \iip.
\end{flexstate}
\begin{proof}
  Suppose $J\idealin \C_2$ and $J\cap \alpha(\C_1)=\{0\}$.  Then
  $J\cap\alpha(\D_1)=\{0\}$.  By hypothesis, $(\D_2,\alpha(\D_1))$ has
  the \iip, so $J\cap \D_2=\{0\}$.  Since $(\C_2,\D_2)$ is a \pd, both
  $(\D_2^c,\D_2)$ and $(\C_2,\D_2^c)$ have the \iip.  As $(J\cap
  \D_2^c) \cap \D_2\subseteq J\cap \D_2=\{0\}$, we find $J\cap
  \D_2^c=\{0\}$, whence $J=\{0\}$.
\end{proof}

\begin{lemma}\label{idbim}
Let $(\C,\D)$ be a \pd\ with Cartan
  envelope $(\A,\B\ms \alpha)$.    Let $J\idealin \C$, and define
  $J_1\subseteq \A$ to be the
  norm-closed $\B$-bimodule generated  by $\alpha(J)$. Then $J_1\idealin\A$.
\end{lemma}
\begin{proof}
Suppose $x\in J$,  $h\in \B$ and $v\in \N(\C,\D)$.  We first show
that
\begin{equation}\label{idbim1}
  \alpha(v)h\alpha(x)\in J_1.
\end{equation}
Since
$\alpha$ is a regular map, $\alpha(v)\in \N(\A,\B)$.
By Lemma~\ref{Tv2.2}, given $u\in
\overline{\alpha(v^*v)\B}$, 
\[\alpha(v)uh\alpha(x)=\theta_{\alpha(v)}^{-1}(uh)\alpha(v)\alpha(x).\]
As
$\theta_{\alpha(v)}^{-1}(uh)\in
\overline{\alpha(vv^*)\B}$ and $\alpha(vx)\in \alpha(J)$, we see
that 
\begin{equation}\label{idbim1.1}\alpha(v)uh\alpha(x)\in J_1.
\end{equation}
Since $(\A,\B)$ is Cartan
and $\alpha$ is a regular map, Lemma~\ref{relcom}\eqref{relcom2} gives
$\alpha(v^*v)\in \B^c=\B$.   By considering an approximate unit for
$\B$, we find  $\alpha(v^*v)\in
\overline{\alpha(v^*v)\B}$. Since $\alpha(v)\alpha(v^*v)^{1/n}\rightarrow
\alpha(v)$, taking $u=\alpha(v^*v)^{1/n}$ in~\eqref{idbim1.1},
we obtain
\eqref{idbim1}.

Next using \eqref{idbim1}, we see that if $n\in \bbN$ and 
$x_i\in J$,  $h_i, k_i\in \B$  for $1\leq
i\leq n$  and $v\in \N(\C,\D)$,
\begin{equation*}
  \alpha(v)\left( \sum_{i=1}^n h_i\alpha(x_i)k_i\right)\in J_1,
\end{equation*}
from which it follows that
\begin{equation}\label{idbim2} \alpha(v)J_1\subseteq J_1.
\end{equation}

Let $M$ be the set of all finite products of elements of
$\B\cup \alpha(\N(\C,\D))$.  An induction argument
using~\eqref{idbim2} shows that for $w\in M$ and $x_1\in J_1$,
$wx_1\in J_1$.  As $\spn M$ is dense in $\A$, we conclude that $J_1$
is a left ideal of $\A$.  Similar arguments show $J_1$ is a right
ideal, so $J_1\idealin\B$.
\end{proof}

\begin{theorem}\label{simple}  Let $(\C,\D)$ be a \pd\ with Cartan
  envelope $(\A,\B\ms \alpha)$.   Then
  \begin{enumerate} \item\label{simple1} $\C$ is simple if and only if
    $\A$ is simple; and 
\item\label{simple3}  $\C$ is unital if and only if $\A$ is unital.
  \item\label{simple2} $\C$ is separable if and only if $\A$ is separable.
  \end{enumerate}
\end{theorem}
\begin{proof}
  Proposition~\ref{CE4rcom} shows that $(\A,\B\ms \alpha)$ is also a
  Cartan envelope for $(\C,\D^c)$.  Thus, by replacing $\D$ with
  $\D^c$ if necessary, without loss of generality we may assume
  $(\C,\D)$ is a \pd\ with the AUP.

  Let $\Delta: \A\rightarrow \B$ be the conditional expectation and
let $E:\C\rightarrow \B$ be the map \[E:=\Delta\circ\alpha.\]

\eqref{simple1} Suppose $\C$ is simple and let $J\idealin \A$.  Then
$J\cap \alpha(\C)\in \{\{0\},\alpha(\C)\}$.  If $J\cap\alpha(\C)=\{0\}$,
Observation~\ref{Cext} shows $J=\{0\}$.   On the other hand,
suppose $J\cap \alpha(\C)=\alpha(\C)$.  Since we are assuming $(\C, \D)$
has the AUP,
we may choose a net $(e_\lambda)$ in
$\D$ which is an approximate unit for $\C$.  As
$\B=C^*(E(\C))$, we find $(\alpha(e_\lambda))$ is an
approximate unit for $\B$.  But $(\A,\B)$ is a Cartan inclusion,
so by~\cite[Theorem~2.6]{PittsNoApUnInC*Al}, every approximate unit
for $\B$ is an approximate unit for $\A$.  As
$\alpha(\D)\subseteq \alpha(\C)\subseteq J$, we conclude that $J$ contains an
approximate unit for $\A$, whence $J=\A$.  Thus $\A$ is simple.

For the converse, suppose $\A$ is simple.  Let $J\idealin \C$ be a
non-zero ideal, let $K$ be the ideal in $\B$ generated by
$\alpha(J\cap\D)$  and put
\[J':=\{x\in \C: E(x^*x)\in K\}.\]

We aim to show $J'$ is a non-zero ideal in $\C$.  By construction, $J'$
is a closed subset of $\C$.

Next we show $J'$ is a non-zero subspace
of $\C$.
For $\tau\in \hat\B$, $\tau\circ E$ is a state on $\C$, so the map,
\[\C\ni x\mapsto (\tau(E(x^*x)))^{1/2}\] is a
semi-norm on $\C$.  Thus for $x,y\in\C$, $\tau(E((x+y)^*(x+y)))^{1/2} \leq
\tau(E(x^*x))^{1/2} +\tau(E(y^*y))^{1/2}$.   Allowing $\tau$ to vary throughout
$\hat\B$,
we conclude
\[E((x+y)^*(x+y))\leq (E(x^*x)^{1/2}+E(y^*y)^{1/2})^2.\]  Thus when $x, y\in J'$
we obtain  $x+y\in J'$.   That
$J'$ is invariant under scalar multiplication is obvious, so $J'$ is a
closed linear subspace of $\C$.    
Clearly $J\cap \D\subseteq J'$.   Since $(\C,\D^c)$ and $(\D^c,
\D)$ both have the ideal intersection property,  $J\cap \D\neq
\{0\}$, and hence $J'\neq \{0\}$.

We are now ready to show $J'$ is an ideal.  For $v\in \N(\C,\D)$, $v^*(J\cap \D) v\subseteq J\cap \D$, and 
regularity of $\alpha$ gives  $\alpha(v)\in
\N(\A,\B)$.  So for $b\in \B$ and $d\in J\cap D$, we have
\begin{align*}\alpha(v^*)\alpha(d)b\alpha(v)&=\lim
\alpha(v^*)\alpha(d)b\alpha(vv^*)^{1/n} \alpha(v) \\
&=\lim (\alpha(v)^*\alpha(d) \alpha(v))\theta_{\alpha(v)}(\alpha(vv^*)^{1/n}
b)\in K.
\end{align*}
It follows that
\[\alpha(v)^*K\alpha(v)\subseteq K, \quad v\in \N(\C,\D).\]
Therefore, when $v\in \N(\C,\D)$ and $x\in J'$, 
$xv\in J'$ because
\begin{align*}
  E(v^*x^*xv)&=\Delta(\alpha(v)^*\alpha(x^*x)\alpha(v))=\alpha(v)^*\Delta(\alpha(x^*x))\alpha(v)\\
             & =\alpha(v)^*E(x^*x)\alpha(v)\in K.
\end{align*}
Since $\spn\N(\C,\D)$ is
dense in $\C$, we conclude that $J'$ is a right ideal in $\C$.
For $x\in J'$ and $y\in \C$,
$E(x^*y^*yx)\leq\norm{y}^2E(x^*x)$, whence $J'$ is a left ideal in
$\C$.     Thus
$\{0\}\neq J'\idealin\C$.

Next we show that $J\cap \D=\D$.   For this, let $\sigma\in \hat\D$
and 
let $\tau\in \hat\B$ satisfy $\tau\circ\alpha=\sigma$.
Let $M$ be the norm-closed $\B$-bimodule
generated by $\alpha(J')$.   By Lemma~\ref{idbim} and
the hypothesis that $\A$ is simple, \[M=\A.\] 
Let $0\leq b\in \B$ be chosen so that $\norm{b}=1$ and
$\tau(b)=1$.   Since $b\in  M$ there is $n\in \bbN$, $x_i\in
J'$, and $h_i, k_i\in \B$ ($1\leq i\leq n$) so that
\[\norm{b-\sum_{i=1}^n h_i \alpha(x_i) k_i}<1/2.\]  Then 
\[\left|\tau\left(\Delta\left(b-\sum_{i=1}^n h_i \alpha(x_i) k_i\right)\right)\right|=\left|1-\sum_{i=1}^n
    \tau(h_i)\tau(\Delta(\alpha(x_i)))\tau(k_i)\right|<1/2.\]
It follows that there exists $x\in J'$ such that
$\tau(\Delta(\alpha(x)))\neq 0$.  
Hence
\[0\neq |\tau(\Delta(\alpha(x)))|^2=\tau(E(x^*)E(x))\leq
  \tau(E(x^*x)). \]  Since $x\in
J'$, $E(x^*x)\in K$, and thus  $\tau$ does not annihilate
$K$.  But $K$ is the ideal of $\B$ generated by $\alpha(J\cap
\D)$, so $\tau$ does not annihilate $\alpha(J\cap \D)$.
Therefore, $\sigma|_{J\cap \D}\neq 0$.   Since this holds for all
$\sigma\in \hat\D$, we conclude that $J\cap \D=\D$. 

Finally, the fact that $\D$ contains an approximate unit for $\C$
implies that $J=\C$.  Thus $\C$ is simple.

\smallskip

\eqref{simple3} This is Proposition~\ref{simple3.0}.

\smallskip
\eqref{simple2}  Since every non-empty subset of a separable metric space is a
separable metric space, separability of $\A$ implies $\alpha(\C)$ is
separable.   So $\C$ is separable.

Now suppose $\C$ is separable and let $Q\subseteq \C$ be countable and
dense.   Then $\A$ is generated by $\alpha(Q) \cup
\Delta(\alpha(Q))$, so $\A$ is separable.
\end{proof}

It would be desirable to include nuclearity in the list of properties
shared by $\C$ and $\A$.   
\begin{conjecture}\label{nuclear}  Let $(\C,\D)$ be a \pd\ 
and let $(\A,\B\ms \alpha)$ be a Cartan envelope for $(\C,\D)$.  Then
$\C$ is nuclear if and only if $\A$ is nuclear.
\end{conjecture}

\subsection[Pseudo-Cartan Inclusions from Cartan Inclusions]{Constructing Pseudo-Cartan Inclusions from Cartan
  Inclusions}

We have seen that every \pd\ has a Cartan envelope, and we now consider the reverse
process, that of  constructing \pd s from a given Cartan inclusion.
We will work more generally, starting instead with a given \pd\ and
constructing new \pd s from it.
Our next result, Proposition~\ref{pdcons},
extends
parts of~\cite[Lemma~5.26 and Proposition~5.31]{PittsStReInII} to the
setting of~\pd s.

\begin{flexstate}{Proposition}{}\label{pdcons}
  Suppose $(\C_1,\D_1)$
  is a \pd\ and $\D$ is a \cstar-subalgebra of $\D_1$ such that
  $\D\subseteq \D_1$ has the \iip.
Let $\M\subseteq \N(\C_1,\D)$ be a $*$-semigroup such
  that $\D\subseteq \overline{\spn}\, \M$ and set
  $\C=\overline\spn\,\M$.  The following statements hold.
  \begin{enumerate}
  \item\label{pdcons1} $(\C,\D)$ is a \pd.
  \item\label{pdcons3} Let $(\A_1,\B_1\ms \alpha)$  
    be a Cartan envelope for 
  $(\C_1,\D_1)$ and let $\Delta_1:\A_1\rightarrow \B_1$ be the
  conditional expectation.   Set
  \[\B:=C^*(\Delta_1(\alpha(\C)))\dstext{and} \A:=C^*(\alpha(\C)\cup
    \B).
  \]
  Then $(\A,\B\ms \alpha|_\C)$ is a Cartan envelope for $(\C,\D)$.
  \end{enumerate}
\end{flexstate}
  \begin{proof}
       \eqref{pdcons1}
  Since $\M$ is a $*$-semigroup, $\C$ is a \cstaralg.  Therefore, $(\C,\D)$ is
  a regular inclusion.   Since $(\C_1,\D_1\ms \subseteq)$ is an
  essential expansion of $(\C,\D)$, Lemma~\ref{claim2}\eqref{claim2:a}
  shows $(\C,\D)$ has the faithful unique pseudo-expectation
  property.  Hence $(\C,\D)$ is a \pd.

\eqref{pdcons3}
As $(\D_1,\D)$ and $(\B_1,\alpha(\D_1))$ have the \iip, so does
$(\B_1, \alpha(\D))$.  Therefore $(\A_1,\B_1\ms \alpha|_{\D})$ is an
essential and Cartan expansion of $(\C,\D)$.   Applying
Proposition~\ref{reg+Del}, we obtain~\eqref{pdcons3}. 
\end{proof}

While Proposition~\ref{pdcons} applies when
$(\C_1,\D_1)$ is a Cartan inclusion,  
in general, $(\C_1,\D_1)$
  need not be the Cartan envelope of $(\C,\D)$.  In fact, $(\C,\D)$
  can itself be Cartan.   A trivial example of this occurs when $\M$
  is taken to be $\D$, in which case, $(\C,\D)=(\D,\D)$ is a Cartan
  inclusion. Here is a more interesting example where $(\C,\D)$ is a
  Cartan inclusion. 

\begin{remark}{Example}  
  Let the amenable discrete group $\Gamma$ act topologically freely on the
  compact Hausdorff space $X$.  Suppose $Z\subseteq X$ is a non-empty
  closed invariant set which is nowhere dense in $X$.  Then
  \[\D:=\{f\in C(X): f|_Z \text{ is constant}\}\] has the ideal
  intersection property in $C(X)$.  Let $Y=\hat\D$ and fix $x_0\in Z$.   Since every
  $f\in\D$ is constant on $Z$, for every $f\in \D$, $f-f(x_0)\in
  C_0(X\setminus Z)$, and it follows that $Y$ is homeomorphic to the
  one point compactification of $X\setminus Z$.   Let $\pi: C(Y)\rightarrow \D$
  be the inverse of the Gelfand transformation and let $\{U_s:s\in\Gamma\}$
  be the canonical copy of $\Gamma$ in $C(X)\rtimes_r\Gamma$.
  The restriction of the action of $\Gamma$ to $X\setminus Z$ is
  topologically free, and hence determines a topologically free action
  of $\Gamma$ on $Y$.

  Let $(\C_1,\D_1)$ be the Cartan inclusion
  $(C(X)\rtimes_r \Gamma, C(X))$ and let $\M$ be the $*$-semigroup  generated by $\D\cup
  \{U_s:s\in\Gamma\}$.  Apply Proposition~\ref{pdcons} to
  $(\D_1,\D)$ 
  to produce a \pd\ $(\C,\D)$.  Since
  $Z$ is a $\Gamma$ invariant set, $\D$ is invariant under
  $\{U_s\}_{s\in\Gamma}$.  Therefore, the pair $(\pi, U)$ is a
  covariant representation for the action of $\Gamma$ on $Y$.  Since
  $\C$ is generated by $\D$ and $\{U_s\}_{s\in\Gamma}$, we obtain a
  $*$-epimorphism
  $(\pi\times U): C(Y)\rtimes_f \Gamma=C(Y)\rtimes_r\Gamma\rightarrow
  \C$.  The topological freeness of $\Gamma$ acting on $Y$ and
  Proposition~\ref{rcp} implies that $(C(Y)\rtimes_r\Gamma, C(Y))$ is
  a Cartan pair, and in particular, has the \iip.  As $\pi$ is
  faithful, we conclude that $(C(Y)\rtimes \Gamma, C(Y))$ is
  isomorphic to $(\C,\D)$.  Thus $(\C,\D)$ is already Cartan, so is
  its own Cartan envelope.  This shows $(\C_1,\D_1)$ is not the Cartan
  envelope for $(\C,\D)$.
\end{remark}

It is natural to wonder when a reduced crossed product is a \pd.
This question is answered by the following result. 
\begin{flexstate}{Proposition}{} \label{rcp} Let $X$ be a locally
  compact Hausdorff space and suppose $\Gamma$ is a discrete group acting
  as homeomorphisms on $X$.   With the corresponding action of
  $\Gamma$ on $C_0(X)$, 
  the following statements are equivalent:
  \begin{enumerate}
  \item $(C_0(X)\rtimes_r\Gamma, C_0(X))$ is a \pd;
  \item $(C_0(X)\rtimes_r\Gamma, C_0(X))$ is a Cartan inclusion; and
    \item the action of $\Gamma$ on
$X$ is
topologically free.
\end{enumerate}
\end{flexstate}
\begin{proof}
For each of $(a)\Leftrightarrow (b)$ and $(b)\Leftrightarrow (c)$, we
use the fact that 
$(C_0(X)\rtimes_r\Gamma, C_0(X))$ is a regular inclusion and there is a
faithful conditional expectation $\Delta: C_0(X)\rtimes \Gamma
\rightarrow C_0(X)$. 

  $(a)\Leftrightarrow(b)$.
Apply Proposition~\ref{Car+upse}.

$(b)\Leftrightarrow(c)$. By~\cite[Proposition~4.14]{Zeller-MeierPrCrCstAlGrAu},
$(C_0(X)\rtimes_r\Gamma, C_0(X))$ is a MASA inclusion if and only if the action
of $\Gamma$ on $X$ is topologically free.   
\end{proof}

\begin{remark}{Remark} \label{pd_rcp} 
Despite the equivalence of parts (a) and (b) in  Proposition~\ref{rcp}, reduced crossed products can be used to
construct examples of \pd s which are not Cartan inclusions.   Indeed,~\cite[Theorem~6.15]{PittsStReInI} shows that reduced
  crossed products can be used to construct unital virtual Cartan
  inclusions.  Thus by combining Proposition~\ref{pdcons} with
  \cite[Theorem~6.15]{PittsStReInI}, we obtain a wide variety of \pd s.
\end{remark}



\mysec[Twists for Cartan Envelopes] {The Twisted Groupoid of the Cartan Envelope} \label{tgpC}
\numberwithin{equation}{subsection}

In~\cite[Section~7]{PittsStReInII}, we described
the twist for the Cartan envelope of a unital regular inclusion.
Combining Theorem~\ref{props} with results of \cite{PittsStReInII}
allows us to describe the twist associated to the Cartan envelope for
any \pd\ regardless of whether it is unital.  This description is
found in Theorem~\ref{gtwCE} below.

\subsection{Reprising the Unital Case}
\label{ruc}
\rm We begin with reprising some definitions and results
from~\cite[Section~7]{PittsStReInII}.   We first consider
unital regular inclusions with the unique (but not necessarily
faithful) pseudo-expectation property.  We then turn our attention to
unital 
\pd s.

Let $(\C,\D)$ be a
regular and unital inclusion with the unique pseudo-expectation
property and 
let $(I(\D),\iota)$ be an injective
envelope for $(\C,\D)$.   Let $E: \C\rightarrow I(\D)$ be the
pseudo-expectation.   
\begin{enumerate}
\item\label{ruc1} An \textit{eigenfunctional}\index{Eigenfunctional} (see
  \cite[Definition~2.1]{DonsigPittsCoSyBoIs}) on $\C$ is a non-zero
  bounded linear functional $\phi$ on $\C$ which is an eigenvector for
  the left and right actions of $\D$ on the dual space of $\C$.  When
  $\phi$ is an eigenfunctional, there exist unique $r(\phi)\in \hat\D$
  and $s(\phi)\in \hat\D$ such that for every $d_1, d_2\in \D$ and
  $x\in \C$,
  \[\phi(d_1xd_2)=r(\phi)(d_1) \, \phi(x)\, s(\phi)(d_2).\]
If for
every $v\in \N(\C,\D)$,
\[|\phi(v)|^2\in\{0, s(\phi)(v^*v)\},\] $\phi$ is a \textit{compatible
  eigenfunctional}\index{Eigenfunctional!compatible} (\cite[Definition~7.6]{PittsStReInII}).  A
\textit{compatible state}\index{Compatible state} is a state on $\C$ which is also a
compatible eigenfunctional.   The
collection of all compatible eigenfunctionals of unit norm is denoted
$\Eigone_c(\C,\D)$ and $\fS(\C,\D)$ denotes the set of compatible states.
\item \label{ruc2} For $\phi\in \Eigone_c(\C,\D)$,
\cite[Theorem~7.9]{PittsStReInII} shows that there are unique states
$\fs(\phi), \fr(\phi)$ on $\C$ with the following properties:
\begin{itemize}
  \item $s(\phi)=\fs(\phi)|_\D$,
$r(\phi)=\fr(\phi)|_\D$;   \item when  $v\in \N(\C,\D)$ satisfies
$\phi(v)\neq 0$ and $x\in \C$,
\[\phi(vx)=\phi(v)\, \fs(\phi)(x)\dstext{and}
  \phi(xv)=\fr(\phi)(x)\, \phi(v).\]
\end{itemize}

For later use, here are formulas for $\fs(\phi)$ and
$\fr(\phi)$:  for $x\in \C$ and $v\in \N(\C,\D)$ with $\phi(v)>0$,
\begin{equation}\label{rcu2.5}
  \fs(\phi)(x)=\frac{\phi(vx)}{\phi(v)}\dstext{and}
  \fr(\phi)(x)=\frac{\phi(xv)}{\phi(v)}.
\end{equation}
In addition, note that   
(still assuming $v\in \N(\C,\D)$ satisfies $\phi(v)>0$) 
$\fs(\phi)(v^*v)\neq 0$ and for every $x\in \C$,
 \[\phi(x) =\frac{\fs(\phi)(v^*x)}{\fs(\phi)(v^*v)^{1/2}}.\]

\item For
$\rho\in \fS(\C,\D)$ and $v\in \N(\C,\D)$ with $\rho(v^*v)\neq 0$, we
use the notation $[v,\rho]$ for the linear functional 
\[[v,\rho](x):=\frac{\rho(v^*x)}{\rho(v^*v)^{1/2}}.\]
By~\cite[Corollary~7.11]{PittsStReInII},
\[\Eigone_c(\C,\D)=\{[v,\rho]: \rho\in \fS(\C,\D) \text{ and }
  \rho(v^*v)\neq 0\}.\]

\item \label{ruc3}
A \textit{strongly compatible state}\index{Compatible state!strongly} on $\C$ is a state $\rho$ on $\C$
for which there exists $\sigma\in \widehat{I(\D)}$ such that
\begin{equation*}\label{scmptd}
  \rho=\sigma\circ E.
\end{equation*}
We denote the family of all strongly
compatible states on $\C$ by $\fS_s(\C,\D)$.
Since 
  \[\hat\D=\{\rho|_\D: \rho\in \fS_s(\C,\D)\},\] there is a rich supply of
  strongly compatible states.

  By~\cite[Theorem~6.9]{PittsStReInII} (whose statement is reproduced
  in Theorem~\ref{upse=>cov}) every strongly compatible state is a
  compatible state, that is,
\[\fS_s(\C,\D)\subseteq \fS(\C,\D).\]

\item \label{ruc4} A compatible eigenfunctional
$\phi\in\Eigone_c(\C,\D)$ is called a \textit{strongly compatible
  eigenfunctional}\index{Eigenfunctional!strongly compatible} if
$\fs(\phi)\in \fS_s(\C,\D)$. When this occurs, $\fr(\phi)$ also
belongs to $\fS_s(\C,\D)$.

\item \label{ruc5} Let $\Sigma(\C,\D)$ be the collection of all norm-one strongly compatible
eigenfunctionals and put \[G(\C,\D):=\{ |\phi|: \phi\in
\Sigma(\C,\D)\},\] where $|\phi|$ denotes the composition of the absolute value
function on $\bbC$ with $\phi$.
\item \label{ruc6}
The representation of $\phi\in
\Sigma(\C,\D)$ as $\phi=[v,\fs(\phi)]$, allows the sets
$\Sigma(\C,\D)$ and $G(\C,\D)$  to
be equipped with groupoid
operations~\cite[Definition~7.17]{PittsStReInII}.   Upon doing so, 
the unit space, $\unit {G(\C,\D)}$, may be identified with $\fS_s(\C,\D)$. With the topology of
pointwise convergence, $\Sigma(\C,\D)$    and $G(\C,\D)$ become Hausdorff topological
groupoids (take $F=\fS_s(\C,\D)$
in~\cite[Theorem~7.18]{PittsStReInII})
and we obtain a
twist
\begin{equation*}\label{twist(C,D)}
\unit{G(\C,\D)}\times\bbT\hookrightarrow \Sigma(\C,\D)\twoheadrightarrow
G(\C,\D)
\end{equation*}
associated with $(\C,\D)$.
\item \label{ruc7} For each $a\in \C$, let $\fg(a):\Sigma(\C,\D)\rightarrow\bbC$ be given by
\begin{equation*}\label{fgdef} \fg(a)(\phi)=\phi(a).
\end{equation*}

\item \label{ruc8} Set 
\[\L(\C,\D):=\{a\in \C: \rho(a^*a)=0 \text{ for all } \rho\in
  \fS_s(\C,\D)\}=\{a\in \C: E(a^*a)=0\}.\]
By~\cite[Theorem~6.5]{PittsStReInII}, $\L(\C,\D)$ is an ideal of $\C$
such that $\L(\C,\D)\cap \D=\{0\}$.  Also, if $J\idealin\C$ satisfies $J\cap
\D=\{0\}$, then $J\subseteq \L(\C,\D)$. 

\end{enumerate}

Now suppose $(\C,\D)$ is a  unital \pd.  Recall from
Proposition~\ref{PZC22}
that $(\D^c,\D)$ has the faithful
unique pseudo-expectation property and $E|_{\D^c}$ is a
$*$-monomorphism;  also $E|_{\D^c}$ is the pseudo-expectation for $(\D^c,\D)$.
 Simplify the notation of~\ref{ruc}\eqref{ruc5}
somewhat by writing $\Sigma$ and $G$ instead of $\Sigma(\C,\D) $ and
$G(\C,\D)$.  As in~\ref{ruc}\eqref{ruc6}, let 
\begin{equation*}\label{tunit}
  \unit G\times \bbT\hookrightarrow \Sigma\twoheadrightarrow G,
\end{equation*} be the twist for $(\C,\D)$.  We remind the reader that
$\unit G$ is identified with $\fS_s(\C,\D)$.

By \cite[Theorem~7.24]{PittsStReInII}, $\fg$ induces a regular
$*$-monomorphism $\theta:\C\rightarrow C^*_r(\Sigma,G)$ and by~\cite[Corollary~7.30]{PittsStReInII}, 
\begin{equation}\label{ruc9} (C^*_r(\Sigma,G), C(\unit{G})\ms \theta)
\end{equation}
is a package for $(\C,\D)$.
While~\cite[Corollary~7.31]{PittsStReInII} states
$(C^*_r(\Sigma,G), C(\unit{G})\ms \theta)$ is a Cartan envelope for
$(\C,\D)$, full details
were  omitted in~\cite{PittsStReInII};~\cite[Remark~1.2]{PittsCoStReInII} provides the
missing details.

\subsection{The Non-Unital Case}\label{nuc}
Let $(\C,\D)$ be a non-unital \pd.  As $(\C,\D)$ is \wnd,
Lemma~\ref{sameunit}\eqref{sameunit1} shows there are two cases: $\C$
has no unit, in which case $\D$ also has no unit; or $\C$ is unital, and
$\D$ is not.  Should $\C$ have a unit, so does $\D^c$.  In this case, 
Proposition~\ref{CE4rcom} shows that describing the Cartan envelope for
$(\C,\D)$ amounts to describing the Cartan envelope for $(\C,\D^c)$,
so the unital case applies.  Thus, throughout the following
discussion, we assume $\C$ has no unit.

Fix an injective envelope $(I(\D),\iota)$ for $\D$.   As
$(\D^c,\D)$ has the \iip, there is a unique $*$-monomorphism
$u:\D^c\rightarrow I(\D)$ such that $u|_\D=\iota|_\D$.  Thus $(I(\D),
\tilde u)$ is an injective envelope for $\D^c$.   It follows that 
$E:\C\rightarrow I(\D)$ is a pseudo-expectation for $(\C,\D)$ relative
to $(I(\D), \iota)$ if and
only if $E$ is a pseudo-expectation for $(\C,\D^c)$ relative to
$(I(\D),u)$.
Since $(\C,\D^c)$ and $(\C,\D)$ have the faithful
unique pseudo-expectation property, we may define $\fS_s(\C,\D)$ and
$\fS_s(\C,\D^c)$ 
as in \ref{ruc}\eqref{ruc3}. 
Then 
\begin{equation}\label{scpst} \fS_s(\C,\D)=\fS_s(\C,\D^c).
\end{equation}
Recalling from Proposition~\ref{CE4rcom} that $(\C,\D)$ and $(\C,\D^c)$
have the same Cartan envelope, we will give a groupoid description
for a Cartan envelope of $(\C,\D^c)$.

Let $\tilde\Sigma$ be the set of 
 norm-one strongly compatible eigenfunctionals for
 $(\tilde\C,\udc\D)$ and let $\tilde
 G=\{|\phi|: \phi\in \tilde\Sigma\}$.   Applying the discussion for
 the unital case to 
 $(\tilde\C,\udc\D)$ we obtain the twist,
 \[\unit{\tilde G}\times \bbT\rightarrow \tilde
   \Sigma\twoheadrightarrow \tilde G\] with unit space $\unit{\tilde
   G}=\fS_s(\tilde\C,\udc\D)$.  Taking  $\theta$ as
 in~\eqref{ruc9}, $(C^*_r(\tilde \Sigma, \tilde G), C(\unit{\tilde G})\ms  \theta)$ is a Cartan
 envelope for $(\tilde\C,\udc\D)$.   We use $\Delta$ to denote
 the conditional expectation of $C^*_r(\tilde\Sigma,\tilde G)$ onto
 $C(\unit{\tilde G})$; note that  $\Delta$ arises from the
 restriction of an element in $C_c(\tilde\Sigma, \tilde G)$ to
 $\unit{\tilde G}$.

Let
$q$ denote the linear functional,  $\tilde\C\ni
(c,\lambda)\mapsto \lambda$.   Since $q$ is a multiplicative linear
functional, it is a compatible state for $(\tilde\C, \udc\D)$, and we claim it is actually a
strongly compatible state.   Let $\psi$ be a multiplicative linear
functional on $I(\D^c)$
  satisfying $\psi\circ\tilde u=q|_{\udc\D}$.  Let $E: \C\rightarrow I(\D)$  be the
  pseudo-expectation and $\tilde E$  its unitization.
  Given $v\in \N(\C,\D^c)$, we
  have
  \[\tilde E(v)^*\tilde E(v)\leq \tilde E(v^*v)=
    u(v^*v).\]  Applying $\psi$ yields,
  \[|\psi(\tilde E(v))|^2\leq \psi(u(v^*v))=q(v^*v)=0.\]  It follows that
  $\psi\circ \tilde E$ annihilates $\N(\C,\D^c)$, and hence
  annihilates $\C$.   As $q$ and $\psi\circ\tilde E$
  are states on $\tilde\C$ which have the same kernel, we obtain
  \[q=\psi\circ\tilde E,\] showing that $q$ is a strongly compatible state.

 Let
\[\Sigma:=\{\phi\in \tilde\Sigma:
\phi|_\C\neq 0\} \dstext{and} G:=\{|\phi|: \phi\in \tilde\Sigma \text{
  and } 
\phi|_\C\neq 0\}\] be equipped with their relative topologies.
Note that if $\phi\in \tilde\Sigma$ satisfies $\phi|_\C=0$, then
$\ker\phi=\ker q=\C$, so $\phi$ is a scalar multiple of $q$.   It
follows that
\[\Sigma=\tilde\Sigma\setminus \bbT q\dstext{and} G=\tilde G\setminus
  \{|q|\}.\]

The definition of the topology and groupoid operations on $G$ show
that $G$ is an open 
subgroupoid of $\tilde G$. Since $\unit{\tilde G}=\fS_s(\tilde\C,\udc\D)$,
\[\unit{G}=\{\rho\in \fS_s(\tilde\C,\udc\D): \rho|_\C\neq 0\}.\]
  By~\cite[Lemma~2.7]{BrownFullerPittsReznikoffGrC*AlTwGpC*Al}.
this produces the twist,
\begin{equation}\label{tnunit}
  \unit{G}\times \bbT\hookrightarrow \Sigma\twoheadrightarrow G.
\end{equation}

Routine
modifications to the proof of~\cite[Theorem~7.24]{PittsStReInII} show
that with $\theta: \tilde\C\rightarrow C^*_r(\tilde\Sigma, \tilde G)$
arising from $\fg$ as described in the unital case above,
\[C^*(\theta(\C)\cup C_0(\unit{G}))= C^*_r(\Sigma,G).\] We now show
that $C^*(\Delta(\theta(\C)))=C_0(\unit G)$. Let
\[\N_0(\C,\D^c):=\{v\in \N(\C,\D^c): \fg(v) \text{ is compactly
    supported}\}\] and let $\C_0=\spn\N_0(\C,\D^c)$.  As in the proof of
~\cite[Theorem~7.24]{PittsStReInII}, $\N_0(\C,\D^c)$ is a $*$-semigroup
and $\C_0$ is dense in $\C$.  Now suppose $\rho_1$ and $\rho_2$ are
distinct elements of $\unit G$, in other words, $\rho_1$ and $\rho_2$
are distinct strongly compatible states on $\tilde\C$, neither of
which vanish on $\C$.  Then $\rho_1|_\C\neq \rho_2|_\C$, so we may
find an element $v\in \N_0(\C,\D^c)$ such that
$\rho_1(v)\neq \rho_2(v)$.  Thus, $\fg(v)(\rho_1)\neq
\fg(v)(\rho_2)$.  Since $\Delta$ is determined  by restriction to
$\unit{G}$, we obtain $\Delta(\fg(v))(\rho_1)\neq \Delta(\fg(v))(\rho_2)$.  
As $\theta(v)=\fg(v)$, we conclude that $\Delta(\theta(\C))$ separates
points of $\unit G$.  Also, recalling that $\D^c$ contains
an approximate unit for $\C$, we see that given $\rho\in \unit G$,
there exists
$h\in \D^c$ such that $\rho(h)=\theta(h)(\rho)\neq 0$.  By the
Stone-Wierstrau\ss\ theorem,  $C^*(\Delta(\theta(\C)))=C_0(\unit G)$.

An application of Theorem~\ref{props} shows that $(C^*_r(\Sigma,
G), C_0(\unit G)\ms  \theta|_\C)$ is a Cartan envelope for $(\C,\D^c)$,
giving a groupoid description for $(\C,\D^c)$ (and $(\C,\D)$).

The definition of compatible eigenfunctional makes sense regardless of
whether the \pd\ $(\C,\D)$ is unital.   Since $\N(\C,\D)\subseteq \N(\C,\D^c)$,
\[\Eigone_c(\C,\D^c)\subseteq \Eigone_c(\C,\D).\]    Actually,
equality holds.  To see this, let $\phi\in \Eigone_c(\C,\D)$.
Regularity of $(\C,\D)$ implies there exists $v\in \N(\C,\D)$ such
that $\phi(v)>0$.   With $\fs(\phi)$ as defined in~\eqref{rcu2.5}, we
obtain
\[\phi=[v,\fs(\phi)]\in \Eigone_c(\C,\D^c).\]  
Using~\eqref{scpst}, we conclude that the collection of norm-one strongly compatible
eigenfunctionals for $(\C,\D)$ is the same as the set of norm-one
strongly compatible eigenfunctionals for $(\C,\D^c)$.

We 
summarize our discussion with the following result, which
extends the groupoid description of the Cartan envelope given
in~\cite[Corollary~7.31]{PittsStReInII} (once again, with
$F=\fS_s(\C,\D))$) from the unital setting to include both the unital
and non-unital cases.

\begin{theorem}\label{gtwCE}   Let $(\C,\D)$ be a \pd,
 let $\Sigma$ be the collection of all norm-one strongly
  compatible eigenfunctionals for $(\C,\D)$, and let
  $G=\{|\phi|:\phi\in\Sigma\}$.  Then the following statements hold.
  \begin{enumerate}
    \item 
  $\Sigma$ and $G$ are Hausdorff topological groupoids with $G$
  effective and 
  \'etale, $\unit G=\fS_s(\C,\D)$ and
  \[\unit G\times \bbT\hookrightarrow \Sigma\twoheadrightarrow G\] is
  a twist.
  \item The map $\fg$ extends to  a regular $*$-monomorphism $\theta: \C\rightarrow
    C^*_r(\Sigma, G)$ and
    \[(C^*_r(\Sigma,G), C_0(\unit G)\ms  \theta)\] is a Cartan envelope
    for $(\C,\D)$.
  \end{enumerate}
\end{theorem}



\mysec[Constructions and Properties]{Constructions and Properties of Pseudo-Cartan Inclusions} \label{const}
\numberwithin{equation}{subsection}

In this \LCsection, we explore the behavior of \pd s and their Cartan
envelopes under mappings and also some familiar
constructions.   Theorem~\ref{cmap} shows that under a suitable
hypothesis, a regular map between
\pd s extends to the Cartan envelopes and Proposition~\ref{autoext}
shows that under a mild hypothesis, a regular automorphism of a \pd\ uniquely extends to its Cartan
envelope.    
 We examine the behavior of \pd s and their inductive limits 
for suitable connecting maps in Theorems~\ref{indlim}
and~\ref{indlimEnv}, and Theorem~\ref{gentp} describes behavior of \pd
s and their Cartan envelops under the minimal tensor product.    In
some of the results, obtaining regularity of maps under these constructions is
a delicate and technical issue.

\subsection{Mapping Results}\label{maprop}
The purposes of this subsection are to establish
a mapping property, Theorem~\ref{cmap}, for Cartan
envelopes and to show that a regular automorphism of a \pd\ extends
uniquely to its Cartan envelope,  Proposition~\ref{autoext}.
While interesting in its own right, Theorem~\ref{cmap} is
a key tool for studying inductive limits of \pd s.   

We require some preparation, beginning with a simple fact about sums
of normalizers.
The sum of normalizers is not usually a normalizer. However,  when the
normalizers are ``orthogonal,'' their sum is again a
normalizer.  Here is the precise statement. 
\begin{flexstate}{Fact}{}\label{nsum}  Suppose $(\C,\D)$ is an
  inclusion and $w, v\in \N(\C,\D)$.   If $w$ and $v$ are orthogonal
  in the sense that $v^*w=0=wv^*$, then $v+w\in \N(\C,\D)$.
\end{flexstate}
\begin{proof}
  Let $d\in \D$.  We have
  \begin{align*}
    \begin{split}
      (v+w)^*d(v+w)&=v^*dv+w^*dw+v^*dw+w^*dv\\
    &=v^*dv+w^*dw+\lim v^*d(ww^*)^{1/n}w+\lim w^*d(vv^*)^{1/n}v\\
    &=v^*dv+w^*dw+\lim v^*w\theta_w(d(ww^*)^{1/n})\\
    &\qquad +\lim
      w^*v\theta_v(d(vv^*)^{1/n})\\
                 &=v^*dv+w^*dw\in \D.
\end{split}
               \end{align*}
  Similarly $(v+w)d(v+w)^*\in\D$, so $v+w\in \N(\C,\D)$.
\end{proof}

Suppose $(\C,\D)$ is an unital inclusion and 
 $(I(\D),\iota)$ is an injective envelope for $\D$.    Given $v\in
\N(\C,\D)$, let $P, Q$ be the support projections in $I(\D)$ for the ideals
$\overline{vv^*\D}$ and $\overline{v^*v\D}$ respectively. (Given any
ideal $J\idealin\D$, the supremum in
$I(\D)_{s.a.}$ of $\{\iota(x): 0\leq x\in J \text{ and } \norm{x}\leq 1\}$ is called the support
projection for $J$.) 
Recall from~\cite[Proposition~1.11]{PittsStReInI} that
the isomorphism 
$\theta_v: \overline{vv^*\D}\rightarrow \overline{v^*v \D}$ uniquely
extends to a $*$-isomorphism $\tilde\theta_v: I(\D)P\rightarrow
I(\D)Q$.   We then obtain a Frol\' ik decomposition~$\{R_j\}_{j=0}^4$
  for $\tilde \theta_v$, see~\cite[Definition~2.12]{PittsStReInI}, and
  also the Frol\'ik ideals $\{K_i(v)\}_{i=0}^4$,~\cite[Definition~2.13]{PittsStReInI}.    Our immediate  interest will
  be with $K_0(v)$.  This is the \textit{fixed
  point ideal} for $v\in \N(\C,\D)$; its definition is 
\begin{equation}\label{FId0}
    K_0(v):=\iota^{-1}(R_0I(\D)).
  \end{equation}   This is a regular ideal in $\D$, and it has the
  following alternate description:
\begin{equation}\label{FId}
  K_0(v)=\{d\in (vv^*\D)^\dperp: vd=dv\in \D^c\}=\{d\in
  (v^*v\D)^\dperp: vd=dv\in \D^c\}.
\end{equation}
See~\cite[Definition~2.13 and Lemma~2.15]{PittsStReInI} for further
details.

To extend the previous considerations to \wnd\ inclusions, recall that
if $(\C,\D)$ is a \wnd\ inclusion, then
$(\tilde\C,\tilde\D)$ is a unital inclusion and the standard
embeddings (described
in~\eqref{unitizedef1}) satisfy $u_\C|_\D=u_\D$.  For $v\in \N(\C,\D)$, define
\begin{equation}\label{FId0nu}
  K_0(v) :=u_\D^{-1}(K_0(u_\C(v)))
  =\{d\in \D: u_\D(d)\in K_0(u_\C(v))\}. 
\end{equation}

When a normalizer $v$ for a regular MASA inclusion is not annihilated
by the pseudo-expectation,  $K_0(v)$ is non-trivial.  In fact,
a bit more is true.
  \begin{lemma}\label{fixidealelt}
 Let
$(\C,\D)$ be a regular MASA inclusion, let $(I(\D),\iota)$ be an
injective envelope for $\D$ and let $E:\C\rightarrow I(\D)$ be the
pseudo-expectation.   Suppose $v\in \N(\C,\D)$ satisfies $E(v)\neq
0$.   Then $K_0(v)\neq \{0\}$ and there exists $k\in K_0(v)\cap
\overline{v^*v\D}$ such that $vk\neq 0$.
\end{lemma}
\begin{proof}
 Assume first that in
addition to the hypotheses stated, 
$(\C,\D)$ is a unital inclusion.
By~\cite[Theorem~3.5(iii)]{PittsStReInI}, $|E(v)|^2=R_0(\iota(v^*v))$,
so $R_0\neq 0$.  Therefore, $K_0(v)\neq \{0\}$ by~\eqref{FId0} and the
fact that $(I(\D),\D,\iota)$ has the \iip.
Next,~\cite[Lemma~2.15]{PittsStReInII} shows $K_0(v)=(K_0(v)\cap
\overline{v^*v\D})^\dperp$, thus $K_0(v)\cap
\overline{v^*v\D}\neq\{0\}$.   Choose $0\neq k\in
K_0(v)\cap\overline{v^*v\D}$.  Then $vk\neq 0$.  (Otherwise,
$v^*vk=0$, whence $k\in \{v^*v\D\}^\perp\cap
\{v^*v\D\}^\dperp$, contrary to $k\neq 0$.)  This completes the unital case.

Now assume $(\C,\D)$ is not unital.  First observe
that $(\tilde\C,\tilde\D)$ is a unital regular MASA inclusion.
Indeed, \cite[Theorem~2.5]{PittsNoApUnInC*Al} shows $(\C,\D)$ has
the AUP, so Lemma~\ref{relcom}\eqref{relcom5} yields regularity of
$(\tilde\C,\tilde\D)$ and a routine argument shows $(\tilde\C,\tilde\D)$ is a MASA
inclusion.  Since $(\C,\D)$ has the AUP, the standard
embedding $u_\C:\C\rightarrow \tilde\C$ is a regular map by
Lemma~\ref{relcom}\eqref{relcom3}.  Thus,
$(v,0)=u_\C(v)\in\N(\tilde\C,\tilde\D)$.
By the unital case, there exists $(k,\lambda)\in K_0((v,0))\cap
\overline{(v^*v,0)\tilde\D}$, with $(v,0)(k,\lambda)\neq (0,0)$.
Since $(k,\lambda)\in \overline{(v^*v,0)\tilde\D}$, we must have $\lambda=0$.   In other words, there 
exists $k\in K_0(v)\cap \overline{v^*v\D}$ such that  $vk\neq 0$.
\end{proof}

Sometimes it is possible to establish regularity of a $*$-monomorphism between inclusions with only partial knowledge of the normalizers in the domain of the map.   The following useful technical result gives a setting where this can be done.

\begin{proposition}
  \label{EWork}
  Let $(\A_2,\B_2\ms \alpha)$ be an essential and Cartan expansion
  for the \pd\ $(\A_1,\B_1)$ and suppose $N\subseteq \N(\A_1,\B_1)$ is a
        $*$-semigroup such that $\B_1\subseteq N$ and $\spn N$ is dense
        in $\A_1$.  If $\alpha(w)\in \N(\A_2,\B_2)$ for every $w\in N$,
        then $\alpha: (\A_1,\B_1)\rightarrow (\A_2,\B_2)$ is a regular
        $*$-monomorphism.

\end{proposition}
  \begin{proof}
     Throughout the proof, let $E_2:\A_2\rightarrow\B_2$ be the
     conditional expectation.   The first step in the proof is to establish the following lemma.
\begin{lemma}\label{EWorkb}  For $v\in \N(\A_1,\B_1)$, let
      \[\Z_v:=\{b_1\in \B_1: \alpha(vb_1)\in \N(\A_2,\B_2)\}.\]  If
      $J\idealin \B_1$ is an essential ideal such that 
      $J\subseteq \Z_v$, then for
      every $b_2\in \B_2$, 
      $\alpha(v)^*b_2\alpha(v)\in \B_2$.
    \end{lemma}

    \begin{proof}
We first show that for any $h\in \B_2$ and  $x\in J$,
\begin{equation}\label{EWorkb1} \alpha(v)^*h\alpha(v)\alpha(x)\in
  \B_2.
\end{equation}
Indeed for $y_n:=\alpha((vx)(vx)^*))^{1/n}\in
\overline{\alpha((vx)(vx)^*)\B_2}$,
\[\alpha(v)^*h\alpha(vx) =\lim \alpha(v)^*hy_n\alpha(vx)=
  \lim\alpha(v^*vx)\theta_{\alpha(vx)}(hy_n)\in\B_2;\] thus
\eqref{EWorkb1} holds. 

Fixing $b_2\in \B_2$, put
\[m:=\alpha(v^*)b_2\alpha(v)-E_2(\alpha(v)^*b_2\alpha(v)).\]  We shall show
$m=0$.  To do this, let
\[M:=\{h\in \B_2: E_2(m^*m)h=0\}.\]  Then $M\idealin \B_2$, and
because $M=\{E_2(m^*m)\}^\perp$,  $M=M^\dperp$.
By \eqref{EWorkb1},
\[\overline{\alpha(J)\B_2}\subseteq M.\]

We claim $M$ is an essential
ideal of $\B_2$.  Indeed, suppose $K\idealin \B_2$ satisfies
$K\cap M=\{0\}$.  Put \[L:=\alpha^{-1}(K)\idealin \B_1\] and note that
$L\cap J=\{0\}$ because
\[\alpha(L\cap J)\subseteq K\cap\alpha(J)\subseteq K\cap \overline{\alpha(J)\B_2}\subseteq K\cap M=\{0\}.\]
But $J$ is an essential ideal in $\B_1$, so $L=\{0\}$.  This gives
$K\cap \alpha(\B_1)=\{0\}$.  Then $K=\{0\}$ because $(\B_2,\B_1,
\alpha|_{\B_1})$ has the \iip.
Thus, $M$ is an essential ideal.

Since $M$ is an essential ideal, 
$M^\perp=\{0\}$.  As $M=M^\dperp$, $M=\B_2$.
Hence $E_2(m^*m)=0$, so $m=0$ by faithfulness of $E_2$.   This gives
$\alpha(v^*)b_2\alpha(v)\in \B_2$.  
\end{proof}

We now return to the proof of Proposition~\ref{EWork}.   Suppose first
that $(\A_1,\B_1)$ is a virtual Cartan inclusion. 
Given $v\in \N(\A_1,\B_1)$, we shall show that $\Z_v$
contains an essential ideal.  We require two claims.
\dstate{Claim}{Recall that $\{v^*v\}^\perp=\{b_1\in \B_1: v^*vb_1=0\}$.  Then  
  \[\{v^*v\}^\perp=\{b_1\in \B_1: vb_1=0\}\subseteq \Z_v.\]}
\dproof
To see $\{v^*v\}^\perp\subseteq \{b_1\in \B_1: vb_1=0\}$, note that
for any polynomial $p$ with $p(0)=0$, $0=vp(v^*v)b_1=p(vv^*)vb_1$;
then use
\[0=\inf\{\norm{p(vv^*)v-v}: p \text{ is a polynomial with }
p(0)=0\}.\]   The opposite inclusion is trivial, so $\{v^*v\}^\perp = \{b_1\in \B_1: vb_1=0\}$; this description gives $\{v^*v\}^\perp\subseteq \Z_v$.
\enddproof

\dstate{Claim}{Suppose $\Lambda$ is a non-empty set and $\{J_\lambda:
  \lambda\in \Lambda\}$ is a collection of ideals in $\B_1$ such that:
 for every $\lambda\in \Lambda$,  $J_\lambda\subseteq \Z_v$; and for distinct
  $\lambda, \mu\in \Lambda$, $J_\lambda\cap J_\mu=\{0\}$.   Write
  $\bigvee_{\lambda\in \Lambda} J_\lambda:=\overline\spn
  \bigcup_{\lambda\in \Lambda} J_\lambda$.  Then $\bigvee_{\lambda\in
    \Lambda} J_\lambda \idealin \B_1$ and $\bigvee_{\lambda\in
    \Lambda} J_\lambda \subseteq \Z_v$.}
\dproof
That $\bigvee_{\lambda\in
    \Lambda} J_\lambda $ is an ideal in $\B_1$ is clear.  
Let $h\in \spn\bigcup_{\lambda\in \Lambda} J_\lambda$.  Then we may
find a finite set $F\subseteq \Lambda$ and for each $\lambda\in F$, an element
$h_\lambda\in J_\lambda$
so that
\[h=\sum_{\lambda\in F} h_\lambda.\]
For distinct $\lambda,\mu\in F$, $h_\lambda h_\mu\in J_\lambda\cap
J_\mu=\{0\}$,  so $h_\lambda h_\mu=0$.   A calculation then shows that 
$\{vh_\lambda: \lambda
\in F\}$ is a collection of pairwise orthogonal normalizers in
$\N(\A_1,\B_1)$.    But $h_\lambda\in \Z_v$ for every $\lambda$,
whence   $\{\alpha(vh_\lambda)\}_{\lambda\in F}$ is a  family of
pairwise orthogonal normalizers in $\N(\A_2,\B_2)$.   Fact~\ref{nsum}
shows that $\sum_{\lambda\in F} \alpha(vh_\lambda)\in \N(\A_2,\B_2)$.  Since
\[\alpha(v)\left(\sum_{\lambda\in F} \alpha(h_\lambda)\right)\in
  \N(\A_2,\B_2),\] $h\in \Z_v$.  The proof of Claim~\thecct\ is
completed by observing that $\Z_v$ is a closed set.  (Indeed, if
$z_i\in \Z_v$ converges to $z\in \B_1$, then
\[\alpha(vz)=\lim \alpha(vz_i) \in
\overline{\N(\A_2,\B_2)}=\N(\A_2,\B_2),\] so $z\in \Z_v$.) 
\enddproof
 We are now prepared to show $\Z_v$ contains an essential
ideal.   This fact is trivial when $v=0$, so assume $v\neq
0$. 

Let $\Q$ be a collection of ideals in $\B_1$.  We shall say
$\Q$ is a \textit{nice
  collection of ideals} if:
\begin{itemize}
  \item for each
    $J\in \Q$, $J\subseteq \Z_v$;
  \item If $J, K\in \Q$ and $J\neq K$, then $J\cap K=\{0\}$; and
  \item $\{v^*v\}^\perp \in \Q$.
  \end{itemize}
By Claim~1, the singleton set, $\{\{v^*v\}^\perp\}$, is a nice
collection, so the family $\fQ$ consisting of all nice collections of ideals is non-empty.
  Partially order $\fQ$ by
   inclusion.   Zorn's lemma provides a maximal nice collection of
   ideals; call it $\fM$.  Consider the ideal,
   \[\fJ:=\bigvee_{J\in \fM} J.\]  Claim~2 shows $\fJ\subseteq \Z_v$.
   We will show $\fJ$ is an essential ideal in $\B_1$ by showing
   $\fJ^\perp=\{0\}$.

   Suppose to the contrary that $\fJ^\perp\neq
   \{0\}$.    Since $\{v^*v\}^\perp\subseteq \fJ$,
   \[\fJ^\perp \subseteq \{v^*v\}^\dperp.\]  Note that
   $\{v^*v\}^\dperp=\overline{v^*v\B_1}^\dperp$, and, since
   $\overline{v^*v\B_1}\subseteq \{v^*v\}^\dperp$
   has the \iip, we find
   \[\fJ^\perp\cap \overline{v^*v\B_1}\neq \{0\}.\]

   Let $0\neq b\in \fJ^\perp \cap \overline{v^*v\B_1}$.  If $v^*vb=0$,
   then $b\in \{v^*v\}^\perp\cap \{v^*v\}^\dperp=\{0\}$, which is
   impossible as $b\neq 0$.  Therefore, $vb\neq 0$.  Since
   $(\A_1,\B_1)$ is a \pd, it has the faithful unique
   pseudo-expectation property.  Thus, letting
   $E:\A_1\rightarrow I(\B_1)$ be the pseudo-expectation,
   $E(b^*v^*vb)\neq 0$.  Since $\spn N$ is dense in $\A_1$, there
   exists a sequence $t_\ell\in \spn N$ such that
   $t_\ell\rightarrow vb$.  Then $E(t_\ell^*vb)\neq 0$ for
   sufficiently large $\ell$, so there exists $w\in N$ such
   that \[E(w^*vb)\neq 0.\] Because we have assumed 
   $(\A_1,\B_1)$ is a virtual Cartan inclusion, we may apply
   Lemma~\ref{fixidealelt} to obtain
   $k\in K_0(w^*vb)\cap \overline{(w^*vb)^*(w^*vb)\B_1}$ such that
   \[0\neq w^*vbk\in \B_1.\]
   If $(e_\lambda)$ is an approximate unit for $\fJ^\perp\cap\overline{v^*v\B_1}$, then
   \[w^*vbk=\lim w^*vk(be_\lambda), \dstext{whence } w^*vbk\in
   \fJ^\perp\cap\overline{v^*v\B_1}.\] 
   Therefore $v(w^*vbk)^*\neq 0$, for otherwise $(w^*vbk)^*\in
   \overline{v^*v\B_1}^\perp\cap \overline{v^*v\B_1}$.    Then
   \[0\neq v(w^*vbk)^*=(vk^*b^*v^*)w\in N.\]  Since
   $\alpha(N)\subseteq \N(\A_2,\B_2)$, we conclude that 
\[h:=(w^*vbk)^*\] is a non-zero element of $\Z_v$.

 Let
   $L:=\overline{h\B_1}$ be the ideal generated by $h$.  Then
   $L\subseteq \Z_v$, and since we noted above that  $h^*\in \fJ^\perp\cap \overline{v^*v\B_1}$, we see that for every
   $J\in \fM$, $L\cap J=\{0\}$.  Thus $\fM\cup \{L\}$ is a nice collection of
   ideals properly containing $\fM$.  As this  contradicts maximality of
   $\fM$, we obtain $\fJ^\perp=\{0\}$.  Therefore, $\fJ$ is an essential ideal.

 Apply Lemma~\ref{EWorkb} to $\Z_v$ and  $\Z_{v^*}$ to 
   conclude that $\alpha(v)\in \N(\A_2,\B_2)$.  If follows that  $\alpha$ is a
   regular $*$-monomorphism.  This completes the proof of the
   proposition when $(\A_1,\B_1)$ is a virtual Cartan inclusion.

   Now suppose $(\A_1,\B_1)$ is a general \pd.   By
   Observation~\ref{vcfrompc}, $(\A_1,\B_1^c)$ is a virtual Cartan inclusion.
   Lemma~\ref{relcom}\eqref{relcom1} shows
   $\N(\A_1,\B_1)\subseteq \N(\A_1,\B_1^c)$ and 
   Lemma~\ref{claim2}\eqref{claim2:b} gives $\alpha(\B_1^c)\subseteq \B_2$.
   Let $N'$ be the $*$-semigroup generated by $N\cup\B_1^c$.  Then
   $\B_1^c\subseteq N'\subseteq \N(\A_1,\B_1^c)$.   Writing $w'\in N'$
   as a finite product with factors belonging to $N\cup \B_1^c$, we find
   $\alpha(w')\in \N(\A_2,\B_2)$.   
   Therefore, 
   \[\alpha(\N(\A_1,\B_1))\subseteq \alpha(\N(\A_1,\B_1^c))\subseteq
     \N(\A_2,\B_2),\] with the second inclusion obtained from the virtual
   Cartan case applied to $(\A_1,\B_1^c)$ and $N'$.  This completes the proof.
 \end{proof}

 We need an intertwining property for conditional expectations before
 stating and proving Theorem~\ref{cmap}.  We shall discuss this
 property further in the introduction to Section~\ref{sec:indlim}.
 
\begin{lemma} \label{EWorka}    Let $(\A_2,\B_2\ms \alpha)$ be an essential and Cartan expansion
  for the Cartan inclusion $(\A_1,\B_1)$, and for $i=1,2$, let $E_i:\A_i\rightarrow\B_i$ be the conditional expectations.   Then 
      $\alpha\circ E_1= E_2\circ \alpha$.
    \end{lemma}
    
\begin{proof} 
Let $(I(\B_1),\iota_1)$ be an injective envelope for $\B_1$.  Since
$(\A_2,\B_2\ms\alpha)$ is an essential expansion of $(\A_1,\B_1)$,
$(\B_2,\B_1,\alpha|_{\B_1})$ is an inclusion with the \iip.
Proposition~\ref{PZC22} shows there is a unique
$*$-monomorphism $\iota_2:\B_2\rightarrow I(\B_1)$ such
that \begin{equation}\label{EWork1}
  \iota_1=\iota_2\circ(\alpha|_{\B_1}).
\end{equation}
For $b_1\in \B_1$, we
have
\[(\iota_2\circ E_2\circ\alpha)(b_1)=\iota_2(\alpha(b_1))=\iota_1(b_1),\] so
$\iota_2\circ E_2\circ\alpha$ is a pseudo-expectation for
$(\A_1,\B_1)$.    Since $\iota_1\circ E_1$ is the pseudo-expectation
for $(\A_1,\B_1)$, we have
\[\iota_2\circ\alpha\circ E_1\stackrel{\eqref{EWork1}}{=}\iota_1\circ E_1=\iota_2\circ
  E_2\circ\alpha.\]  Since $\iota_2$ is one-to-one,  the lemma
follows.
\end{proof}

We are now ready to establish a key mapping property of Cartan envelopes.
\begin{theorem}\label{cmap}
  Suppose $(\C_1,\D_1)$ and $(\C_2,\D_2)$ are \pd s and 
$(\C_2,\D_2\ms \alpha)$ is a regular and essential expansion of $(\C_1,\D_1)$.
For $i=1,2$, let $(\A_i,\B_i\ms \tau_i)$ be Cartan envelopes for
$(\C_i,\D_i)$ with conditional expectations $\Delta_i:\A_i\rightarrow
\B_i$.

There is a unique $*$-homomorphism $\mac\alpha:\A_1\rightarrow\A_2$ 
 satisfying:
\begin{align}
\mac\alpha\circ \tau_1&=\tau_2\circ \alpha  \label{cmap1}\\
\mac\alpha\circ\Delta_1&=\Delta_2\circ\mac\alpha \label{cmap1.05}\\
  \intertext{and}
 \mac\alpha\circ\Delta_1\circ\tau_1&=\Delta_2\circ\tau_2\circ\alpha. \label{cmap1.1}
\end{align}
Furthermore, $\mac\alpha$ is a $*$-monomorphism and $(\A_2,\B_2\ms\mac\alpha)$
is a regular and essential expansion of $(\A_1,\B_1)$.
\end{theorem}
The following commutative diagram
illustrates the maps
involved; the unlabelled 
vertical arrows represent inclusion maps.
   \begin{equation} \label{myCD1}
          \xymatrix{
            \A_1 \ar@{-->}[rrr]^{\exists !\, \mac\alpha}\ar @/_2pc/_{\Delta_1}[ddd]& & & \A_2\ar @/^2pc/^{\Delta_2}[ddd]\\
            & \C_1\ar[ul]_{\tau_1}\ar[r]^\alpha & \C_2\ar[ur]^{\tau_2} &\\
            &\D_1\ar[dl]^{\tau_1|_{\D_1}}\ar[r]_{\alpha|_{\D_1}}\ar[u]& \D_2\ar[dr]_{\tau_2|_{\D_2}}\ar[u]& \\
            \B_1\ar[uuu] \ar@{-->}[rrr]_{\mac\alpha|_{\B_1}}& & & \B_2\ar[uuu]
            }
        \end{equation}

        \begin{proof}
  \textit{Existence:}  We first show that there exists a
  $*$-homomorphism $\mac\alpha:\A_1\rightarrow\A_2$
  satisfying~\eqref{cmap1}. 
Note that $(\A_2,\B_2\ms \tau_2\circ\alpha)$ is a Cartan expansion of
$(\C_1,\D_1)$.  Also, the inclusions $(\B_2,\D_2, \tau_2|_{\D_2})$ and $(\D_2,\D_1,
\alpha|_{\D_1})$ have the \iip.  Therefore $(\B_2,\D_1,
\tau_2\circ\alpha|_{\D_1})$ also has the \iip.  Since
$\tau_2\circ\alpha$ is a regular map,  $(\A_2,\B_2\ms
\tau_2\circ\alpha)$ is an essential and regular Cartan expansion of
$(\C_1,\D_1)$.  

 Let
        \begin{equation}\label{mpr0}\B_1':=C^*((\Delta_2\circ\tau_2\circ\alpha)(\C_1))
          \dstext{and} \A_1':=C^*(\B_1'\cup \tau_2(\alpha(\C_1))).
        \end{equation}
Proposition~\ref{reg+Del} shows $(\A_1',\B_1'\ms \tau_2\circ\alpha)$ is a
Cartan envelope for $(\C_1,\D_1)$.

By the uniqueness of the Cartan envelope, there exists
a unique $*$-isomorphism $\psi:\A_1\rightarrow \A_1'$ such that $\psi$
is regular and 
$\psi\circ\tau_1= \tau_2\circ\alpha$.  Let $f:\A_1'\hookrightarrow \A_2$
be the inclusion map.  Taking
\[\mac\alpha:=f\circ\psi,\] we obtain~\eqref{cmap1}.  Since $f$ and
$\psi$ are one-to-one, so is $\mac\alpha$.

We claim $(\A_2,\B_2\ms\mac\alpha)$ is a regular and essential
expansion of $(\A_1,\B_1)$.
By construction, 
$\mac\alpha(\B_1)=\psi(\B_1)\subseteq \B_2$.  Let us  show $(\B_2, \mac\alpha(\B_1))$ has
the \iip.  As
$(A_2,\B_2\ms\tau_2\circ \alpha)$ is a regular and essential Cartan expansion
of $(\C_1,\D_1)$, 
$(\B_2,\tau_2(\alpha(\D_1)))$ has the \iip.  Since
\[\tau_2(\alpha(\D_1))\subseteq \psi(\B_1)=\mac\alpha(\B_1)\subseteq
\B_2,\] $(\B_2, \mac\alpha(\B_1))$ has the \iip, so
$(\A_2,\B_2\ms\mac\alpha)$ is an essential expansion of $(\A_1,\B_1)$.

Let $N$ be the $*$-semigroup generated by
$\B_1\cup \tau_1(\N(\C_1,\D_1))$.  Since $(\A_1,\B_1\ms\tau_1)$ is a
Cartan envelope for $(\C_1,\D_1)$, $N\subseteq \N(\A_1,\B_1)$, and
$\spn N$ is dense in $\A_1$.  Proposition~\ref{EWork}
shows $\mac\alpha$ is a regular $*$-monomorphism.
This completes the  proof of the claim. 

Lemma~\ref{EWorka} (applied to
$\mac\alpha$) gives
~\eqref{cmap1.05}.
Finally,  \eqref{cmap1.05} gives \[\mac\alpha\circ\Delta_1\circ\tau_1=
\Delta_2\circ\mac\alpha\circ\tau_1=\Delta_2\circ\tau_2\circ\alpha,\] with
the last equality following from \eqref{cmap1}.  Thus 
\eqref{cmap1.1} holds, and the proof is existence is complete.

\medskip \noindent\textit{Uniqueness:} Suppose $\mac\alpha':\A_1\rightarrow \A_2$ is a
$*$-homomorphism 
satisfying~\eqref{cmap1}, \eqref{cmap1.05}, and \eqref{cmap1.1}.

Since
\[\mac\alpha\circ\Delta_1\circ\tau_1=\Delta_2\circ \tau_2\circ \alpha= \mac\alpha'\circ\Delta_1\circ\tau_1,\]
and $\B_1$ is generated by $\Delta_1(\tau_1(\C_1))$, we obtain
\[\mac\alpha'|_{\B_1}=\mac\alpha|_{\B_1}.\]  But $\A_1$ is generated by
$\B_1\cup\tau_1(\C_1)$ and $\mac\alpha'\circ
\tau_1=\tau_2\circ\alpha=\mac\alpha\circ\tau_1$, so
\[\mac\alpha'=\mac\alpha.\]  This completes the proof.
\end{proof}

\begin{remark}{Example} \label{nonessEx}
The purpose of this example is to show the uniqueness conclusion of Theorem~\ref{cmap} need not
  hold without the hypothesis that $\alpha$ is an essential map.  
   Section~2 of~\cite{PittsCoStReInII} gives an example of a
  virtual Cartan inclusion $(\C_1,\D_1)$ and describes a Cartan
  package $(\C_2,\D_2\ms \alpha)$ for $(\C_1,\D_1)$ which is not the
  Cartan envelope for $(\C_1,\D_1)$.  A Cartan envelope
  $(\A_1,\B_1\ms \tau_1)$ for $(\C_1,\D_1)$ is also described
  in~\cite[Section~2]{PittsCoStReInII}.  For convenience, we summarize
  a particular case of the discussion in~\cite[Section~2]{PittsCoStReInII}.
  Let $S_-=\begin{bsmallmatrix}0&1\\1&0\end{bsmallmatrix}$,
  $S_+=\begin{bsmallmatrix} 1&0\\0&-1\end{bsmallmatrix}$, and put
  $S_0:=iS_-S_+$.  Then $S_-, S_+$ and $S_0$ are self-adjoint unitary
  matrices.  Put
  \begin{align*}
    \C_1&:=C([-1,1], M_2(\bbC))\\\intertext{and}\D_1&:=\{f\in \C_1: f(t)\in
  C^*(S_-) \text{ if } t<0 \text{ and } f(t)\in C^*(S_+) \text{ if }
                                                      t>0\};
  \end{align*}
  continuity gives $f(0)\in\bbC I$ for $f\in \D_1$.   A
Cartan package $(\C_2,\D_2\ms \alpha)$ for $(\C_1,\D_1)$ is
\begin{align*}\C_2&:=C([-1,0],M_2(\bbC))\oplus M_2(\bbC)\oplus C([0,1],M_2(\bbC)),\\
   \D_2&:=C([-1,0],C^*(S_-)) \oplus C^*(S_0)\oplus C([0,1], C^*(S_+)), \dstext{and}
          \\ 
  \alpha(f)&:=f|_{[-1,0]} \oplus f(0)\oplus f|_{[0,1]}, \quad f\in \C_1.
\end{align*}
A Cartan envelope for $(\C_1,\D_1)$ is $(\A_1,\B_1\ms \tau_1)$, where
\begin{align*}
\A_1&:=C([-1,0],M_2(\bbC))\oplus C([0,1],M_2(\bbC)),\\
   \B_1&:=C([-1,0],C^*(S_-)) \oplus C([0,1], C^*(S_+)), \dstext{and}\\ 
  \tau_1(f)&:=f|_{[-1,0]} \oplus f|_{[0,1]}, \quad f\in \C_1.
\end{align*}

Since $(\C_2,\D_2)$ is already a Cartan inclusion, a Cartan envelope
for $(\C_2,\D_2)$ is just $(\C_2,\D_2\ms \Id|_{\C_2})$.  However, the
uniqueness of $\mac\alpha$ fails.  Indeed,
define  $\mac\alpha_-, \mac\alpha_+:\A_1\rightarrow \C_2$ by
\[\mac\alpha_-(f_1\oplus f_2):=f_1\oplus f_1(0)\oplus f_2\dstext{and}
  \mac\alpha_+(f_1\oplus f_2):=f_1\oplus f_2(0)\oplus f_2.\] 
Thus  the uniqueness conclusion
of Theorem~\ref{cmap} does not hold.
\end{remark}

The uniqueness and minimality statements from Theorem~\ref{!pschar}
give the following rigidity result for the Cartan envelopes of a \pd.
\begin{proposition}\label{autoext}   Let $(\C,\D)$ be a \pd, with Cartan
  envelope $(\A,\B\ms\tau)$.   If $\theta:\C\rightarrow \C$ is a
  $*$-automorphism such that \[\theta(\N(\C,\D))\subseteq
    \N(\C,\D)\dstext{and}\theta(\D)\subseteq \D,\] then $\theta$ uniquely extends to a regular 
  $*$-automorphism $\mac\theta:\A\rightarrow\A$ such that
  \begin{equation}\label{autoex1}
    \mac\theta\circ\tau=\tau\circ\theta.
  \end{equation}
\end{proposition}

The hypothesis that $\D$ is invariant under $\theta$ is automatically satisfied when
$(\C,\D)$ has the AUP, see Lemma~\ref{relcom}\eqref{relcom2}.    

\begin{proof}
We first show that $(\A,\B\ms \tau\circ\theta)$ is a Cartan package
for $(\C,\D)$.   By definition of the Cartan envelope, $\tau$ is a regular map, and as $\theta$ is assumed regular, we find $\tau\circ\theta$ is a regular map.   Let $\Delta:\A\rightarrow \B$ be the conditional
expectation.   Since $(\A,\B\ms\tau)$ is a Cartan envelope for $(\C,\D)$,
\begin{align*}\B&=C^*(\Delta(\tau(\C)))=C^*(\Delta(\tau(\theta(\C))))\\
  \intertext{and}
  \A&=C^*(\tau(\C)\cup \B)=C^*(\tau(\theta(\C))\cup\B).
\end{align*}
Since $\theta(\D)\subseteq \D$, $(\A,\B\ms\tau\circ\theta)$ is a regular
expansion of $(\C,\D)$.  An examination of Definition~\ref{defCarEnv}\eqref{defCarEnv1} shows $(\A,\B\ms\tau\circ\theta)$ is a Cartan package for $(\C,\D)$.

The ideal $\fJ\idealin\A$ obtained from the minimality statement of Theorem~\ref{!pschar} satisfies
\[\{0\}=\fJ\cap \tau(\theta(\C))=\fJ\cap\tau(\C).\]  By Observation~\ref{Cext}, $(\A,\C,\tau)$ has the \iip, so $\fJ=\{0\}$.    Therefore, $(\A,\B\ms\tau\circ\theta)$ is a Cartan envelope for $(\C,\D)$.

Apply the uniqueness statement of Theorem~\ref{!pschar} to obtain a unique regular $*$-automorphism $\mac\theta$ of $\A$ satisfying~\eqref{autoex1}.
\end{proof}
\begin{remark}{Remark}\label{autoextrem}
Let $(\A,\B)$ be a Cartan inclusion, or more generally, a regular MASA  inclusion,  and
suppose $\mac\theta: \A\rightarrow \A$ is a regular $*$-automorphism.
By Corollary~\ref{cor:rmasa}, $\mac\theta(\B)=\B$ and $\mac\theta^{-1}$ is  regular.
In particular, the family of regular
automorphisms of a Cartan inclusion or of a regular MASA inclusion
is a group.

Automorphisms of certain Cartan inclusions and virtual Cartan
inclusions satisfying $\mac\theta(\B)=\B$ have been studied,
see~\cite[Corollary~2.2]{Komura*HoBeGrC*Al} and
\cite[Theorem~6.10]{TaylorEsCaSuC*Al}.  We expect it is possible to
obtain a description of those regular $*$-automorphisms of a \pd\
$(\C,\D)$ for which $\D$ is an invariant subspace by combining these
results with Proposition~\ref{autoext}.
\end{remark}

If the MASA hypothesis in Remark~\ref{autoextrem} is dropped, the inverse of a regular automorphism need not be regular.    Here is an example of a \pd\ where the family of regular $*$-automorphisms is not a group.

\begin{example}\label{needraut}
  Let \begin{align*}
        \C&=\{h\in C(\bbR): \lim_{t\rightarrow -\infty} h(t)=0\text{ and }
            \lim_{t\rightarrow +\infty} h(t) \text{ exists}\}\\
                                        \intertext{and}
                                        \D&=\{h\in\C:   h|_\bbN \text{ is constant}\}.
                                      \end{align*}
                                      Since $\bbN\subseteq \bbR$ has empty interior,  $(\C,\D)$ has the \iip.   Also,
                                      \begin{equation}\label{needraut1}
                                        \N(\C,\D)=\{h\in \C: |h| \text{ is constant on }\bbN\},
                                      \end{equation}
                                      so
                                      \[N:=\{h\in \C: h(n)\in \bbT \text{ for all }n\in\bbN\}\] is a $*$-semigroup of normalizers.  The Stone-Weiererstrau\ss\ theorem shows $\spn N$ is a dense subalgebra of $\C$.  Thus $(\C,\D)$ is regular, and hence a \pd.

For $f\in \C$, define $\theta(f)(t)=f(t+1)$.   Using~\eqref{needraut1}
we see $\theta$ is regular and leaves $\D$ invariant.  However,
$\theta^{-1}$ is not a regular map.
                                    
\end{example}

\subsection{Inductive Limits}\label{sec:indlim}
In ~\cite{KumjianOnC*Di}, Kumjian introduced \cstar-diagonals, a class
of inclusions subsequently broadened by Renault to the class of Cartan
inclusions.  An inclusion $(\C,\D)$ is a \textit{\cstar-diagonal} if 
it is a Cartan inclusion such that every pure state on $\D$ extends
uniquely to a state on $\C$.  (This definition is equivalent to
Kumjian's original definition,
see~\cite[Proposition~2.10]{PittsNoApUnInC*Al}.)  The class of unital
\cstar-diagonals is closed under inductive limits provided the
connecting maps are unital, regular, and one-to-one
(\cite[Theorem~4.23]{DonsigPittsCoSyBoIs}), and it can be shown that
the same is true in the non-unital setting (with connecting maps again
being
regular $*$-monomorphisms).  The unique extension property implies
that the connecting mappings behave well with respect to the
conditional expectations: if for $i=1,2$, $(\C_i,\D_i)$ are
\cstar-diagonals and $E_i:\C_i\rightarrow \D_i$ is the conditional
expectation, then for any regular $*$-monomorphism
$\alpha_{21}:\C_1\rightarrow \C_2$,
\begin{equation}\label{CE+emb}
  \alpha_{21}\circ E_1=E_2\circ \alpha_{21}.
\end{equation}

Turning to the context of Cartan inclusions,
when the connecting maps are regular and satisfy  \eqref{CE+emb}, Li
showed the
inductive limit of a sequence of Cartan inclusions is again a Cartan inclusion,
see~\cite[Theorem~1.10]{LiXinEvClSiC*AlCaSu}.  Li's result was
extended  for
non-commutative Cartan inclusions by Meyer, Raad and Taylor,
see~\cite[Theorem~3.9]{MeyerRaadTaylorInLiNoCaIn}.

For Cartan inclusions, it follows from
Lemma~\ref{EWorka} and Remark~\ref{R:dextinc}\eqref{R:dextinc2}
 that~\eqref{CE+emb}
holds provided the inclusion $(\D_2, \D_1, \alpha_{21})$ has the
\iip\ and $\alpha_{21}$ is a regular map.    However, when this assumption on $(\D_2,\D_1,\alpha_{21})$ is 
dropped, examples show
\eqref{CE+emb} need 
not hold.  
 Such examples suggest that when the connecting maps are regular, but 
fail to be essential, an
 inductive limit of Cartan inclusions may not be a Cartan inclusion. 

 In this subsection, we consider inductive limits of \pd s.   We show
 that when the connecting maps in
an inductive system are essential and regular $*$-monomorphisms, then
the inductive limit of \pd s is again a \pd.  Further, we show that
the Cartan envelope of the inductive limit of such a system is the
inductive limits of the Cartan envelopes.

  The following gives
the notion of a directed system of \pd s suitable for our purposes.
\begin{definition}\label{dspd}  Let $\Lambda$ be a directed set.
  Suppose
  \begin{itemize}
    \item  for each $\lambda\in \Lambda$,  $(\C_\lambda,\D_\lambda)$ is a
      \pd;
    \item for every $(\lambda,\mu)\in \Lambda\times\Lambda$ with
      $\lambda\geq \mu$,
        $(\C_\lambda,\D_\lambda\ms \alpha_{\lambda\mu})$ is an essential
        and regular expansion of $(\C_\mu,\D_\mu)$; and
        \item whenever $\lambda\geq\mu\geq\nu$,
          $\alpha_{\lambda\nu}=\alpha_{\lambda\mu}\circ\alpha_{\mu\nu}$.
        \end{itemize}
        \begin{enumerate}
          \item We will call the collection
        $\{(\C_\lambda,\D_\lambda\ms \alpha_{\lambda\mu}): \lambda\geq
        \mu\}$ a \textit{system of \pd s directed by
          $\Lambda$}.\index{System of \pd s}  When
        the directed set is clear from context or it is not necessary
        to specify it, we will simplify terminology and say
        $\{(\C_\lambda,\D_\lambda\ms \alpha_{\lambda\mu}): \lambda\geq
        \mu\}$ is a \textit{system} of \pd s.
\item Given a system
        $\{(\C_\lambda,\D_\lambda\ms \alpha_{\lambda\mu}): \lambda\geq
        \mu\}$, recall from the definition of expansion (see
        Definition~\ref{incdef}\eqref{incdef2}) that
        $\alpha_{\lambda\mu}(\D_\mu)\subseteq \D_\lambda$.  Letting
        $\indlim\C_\lambda$ and $\indlim\D_\lambda$ be the inductive
        limit \cstar-algebras, we obtain the inclusion
        $(\indlim\C_\lambda,\indlim\D_\lambda)$.  We will call this
        inclusion the \textit{inductive limit} of the system
        $\{(\C_\lambda,\D_\lambda\ms \alpha_{\lambda\mu}): \lambda\geq
        \mu\}$.

  \item Finally, for each $\lambda\in\Lambda$, we will let
  $\oar{\alpha_\lambda}:\C_\lambda\rightarrow \indlim \C_\lambda$ be the
  \textit{canonical embedding}: it satisfies
  \begin{equation}\label{dspd1}
    \oar{\alpha_\lambda}\circ\alpha_{\lambda\mu}=\oar{\alpha_\mu} \dstext{for every}
  \lambda\geq \mu.
  \end{equation}
\end{enumerate}
\end{definition}

\begin{remark}{Remark}\label{no1}    By Observation~\ref{Cext},
if 
$\{(\C_\lambda,\D_\lambda\ms
\alpha_{\lambda\mu}):\lambda\geq \mu\}$ is a system of \pd s, then
whenever $\lambda\geq \mu$,  $(\C_\lambda,\C_\mu,
\alpha_{\lambda\mu})$ has the \iip.
\end{remark}  

We now  show that the inductive limit of a system of \pd s is again a \pd, and when the system consists of Cartan
inclusions, the inductive limit is Cartan.  (In the latter case, we do not apply
\cite[Theorem~1.10]{LiXinEvClSiC*AlCaSu}, preferring instead to use Theorem~\ref{indlim}\eqref{indlim4}.)

\begin{theorem}\label{indlim}  
Suppose $\Lambda$ is a directed
        set and $\{(\C_\lambda,\D_\lambda\ms\alpha_{\lambda\mu}):
        \lambda\geq \mu\}$ is a system of \pd s.   Let
        \[\C:=\indlim \C_\lambda\dstext{and}\D:=\indlim \D_\lambda.\]
      The following statements hold.
  \begin{enumerate}
    \item\label{indlim4}  $(\C,\D)$ is
      a \pd, and for each $\lambda\in \Lambda$,
      $(\C,\D\ms\oar{\alpha_\lambda})$ is a regular and
      essential expansion for $(\C_\lambda,\D_\lambda)$.  
\item\label{indlim4a}  Suppose each $(\C_\lambda,\D_\lambda)$ is a Cartan
  inclusion and $\Delta_\lambda: \C_\lambda\rightarrow\D_\lambda$ is
  the conditional expectation.  Then $(\C,\D)$ is a Cartan inclusion
  and letting $\Delta:\C\rightarrow \D$ be the conditional expectation,
  we have $\Delta\circ\oar{\alpha_\lambda}=\oar{\alpha_\lambda} \circ
  \Delta_\lambda$.
\end{enumerate}
\end{theorem}
\begin{proof}
  \eqref{indlim4}
             Let us show $\oar{\alpha_\lambda}$ is a regular map of $(\C_\lambda,\D_\lambda)$ into
$(\C,\D)$.   Let  $v\in
  \N(\C_\lambda,\D_\lambda)$ and $d\in \D$.  For $\eps>0$, we may find
  $\mu\geq \lambda$ and $d_\mu\in \D_\mu$ such that
  $\norm{v}^2\norm{d-\oar{\alpha_\mu}(d_\mu)}<\eps$.  Since
  $\oar{\alpha_\mu}\circ\alpha_{\mu\lambda}=\oar{\alpha_\lambda}$ and
  $\alpha_{\mu\lambda}$ is a regular map,  
  \[\oar{\alpha_\lambda}(v)\oar{\alpha_\mu}(d_\mu)\oar{\alpha_\lambda}(v^*)
=\oar{\alpha_\mu}(\alpha_{\mu\lambda}(v)
    d_\mu\alpha_{\mu\lambda}(v)^*)\in \oar{\alpha_\mu}(\D_\mu)\subseteq
    \D.\]  But
  \[\norm{\oar{\alpha_\lambda}(v)d\oar{\alpha_\lambda}(v^*)
-\oar{\alpha_\lambda}(v)\oar{\alpha_\mu}(d_\mu)\oar{\alpha_\lambda}(v^*)}\leq
    \norm{v}^2\norm{d-\oar{\alpha_\mu}(d_\mu)}<\eps.\]  It
  follows that 
  $\oar{\alpha_\lambda}(v)d\oar{\alpha_\lambda}(v)^*\in \D$.  Likewise,
  $\oar{\alpha_\lambda}(v)^*d\oar{\alpha_\lambda}(v)\in \D$.  Thus 
  $\oar{\alpha_\lambda}$ is a regular map.

We now show 
  $(\C,\D)$ is a regular inclusion. 
 Since each $(\C_\lambda,\D_\lambda)$ is a regular inclusion we find
 that 
 \[\C=\overline{\bigcup_\lambda
     \oar{\alpha_\lambda}(\C_\lambda)}=\overline{\bigcup_\lambda\spn
     \oar{\alpha_\lambda}(\N(\C_\lambda,\D_\lambda))}\subseteq \overline\spn
     \N(\C,\D).\] Therefore $(\C,\D)$ is a regular inclusion.

To complete the proof that $(\C,\D)$ is a \pd, we show that $(\C,\D)$
has the faithful unique pseudo-expectation property.  Let $(I(\D),\iota)$ be an injective envelope
for $\D$ and for $\lambda\in \Lambda$, define
$\iota_\lambda:\D_\lambda \rightarrow I(\D)$ by 
\[\iota_\lambda:=\iota\circ \oar{\alpha_\lambda}|_{\D_\lambda}.\]  Note that
for
$\mu\geq \lambda$, \[\iota_\mu\circ
\alpha_{\mu\lambda}|_{\D_\lambda}=\iota_\lambda.\]

We first claim
that for every $\lambda\in\Lambda$, $(I(\D),\iota_\lambda)$ is an
injective envelope for $\D_\lambda$.  It suffices to show the
inclusion 
$(I(\D), \D_\lambda, \iota_\lambda)$ has the \iip.
Let $J\idealin I(\D)$ satisfy $
\iota_\lambda(\D_\lambda)\cap J=\{0\}$.   Choose $\mu\geq \lambda$.  Since
$\alpha_{\mu\lambda}$ is an essential map, 
$(\D_\mu, \D_\lambda, \alpha_{\mu\lambda}|_{\D_\lambda})$ has the
\iip. Therefore, 
 \[(\iota_\mu(\D_\mu), \D_\lambda,
\iota_\mu\circ\alpha_{\mu\lambda}|_{\D_\lambda})=(\iota_\mu(\D_\mu),
\D_\lambda, \iota_\lambda)\] has the \iip\ as
 well.   As $J\cap \iota_\mu(\D_\mu)\cap
\iota_\lambda(\D_\lambda)=J\cap \iota_\lambda(\D_\lambda)=\{0\}$, we see $J\cap \iota_\mu(\D_\mu) =\{0\}$ for every $\mu\geq
 \lambda$.  But $\bigcup_{\mu\geq \lambda}
 \alpha_\mu(\D_\mu)$ is dense
    in $\D$, so $\bigcup_{\mu\geq\lambda} \iota_\mu(\D_\mu)$ is dense
    in $\iota(\D)$.   By~\cite[Proposition~II.8.2.4]{BlackadarOpAl},
    $J\cap\iota(\D)=\{0\}$.   Since $(I(\D), \D,\iota)$ has the \iip,
    we conclude $J=\{0\}$.  Thus our claim holds.

    We are now ready to show $(\C,\D)$ has a unique pseudo-expectation.
    For $i=1,2$, let $E_i:\C\rightarrow I(\D)$ be a
    pseudo-expectation.  Then $E_i\circ\oar{\alpha_\lambda}$ is a
    pseudo-expectation for $(\C_\lambda,\D_\lambda)$ relative to
    $(I(\D),\iota_\lambda)$.  Since $(\C_\lambda,\D_\lambda)$ has the
    faithful unique pseudo-expectation property,
    $E_1\circ\oar{\alpha_\lambda}=E_2\circ\oar{\alpha_\lambda}$.  As this holds
    for every $\lambda$, we conclude that $E_1=E_2$.  Thus, $(\C,\D)$
    has a unique pseudo-expectation $E:\C\rightarrow I(\D)$.

    We now
    show $E$ is faithful.  Let $\L:=\{x\in \C: E(x^*x)=0\}$ be the
    left kernel of $E$.  It follows
    from~\cite[Theorem~6.5]{PittsStReInII} and Lemma~\ref{prepse} that
    $\L$ is an ideal in $\C$.  Then $\oar{\alpha_\lambda}^{-1}(\L)$ is an ideal of
    $\C_\lambda$ contained in
    $\L_\lambda:=\{x\in \C_\lambda: E(\oar{\alpha_\lambda}(x^*x))=0\}$.
    Since $(\C_\lambda,\D_\lambda)$ has the faithful unique
    pseudo-expectation property and $E_\lambda=E\circ\oar{\alpha_\lambda}$
    is the pseudo-expectation for $(\C_\lambda, \D_\lambda)$,
    $\L_\lambda=0$.  We conclude that
    $\L\cap \oar{\alpha_\lambda}(\C_\lambda)=\{0\}$.  Thus,
    using~\cite[Proposition~II.8.2.4]{BlackadarOpAl} yet again, we find
    $\L=\{0\}$.  Therefore, $E$ is faithful.  Since $(\C,\D)$ has the
    faithful unique pseudo-expectation property, it is a \pd.

As  $(I(\D),\iota_\lambda)$ is an injective envelope for $\D_\lambda$,
$(I(\D),\D_\lambda, \iota_\lambda)$ is an essential inclusion.  Since
\[I(\D)\supseteq \iota(\D)\supseteq
\iota_\lambda(\D_\lambda)=\iota(\oar{\alpha_\lambda}(\D)),\]
Lemma~\ref{interiip} shows 
$(\D,
\D_\lambda, \oar{\alpha_\lambda})$ is also an essential inclusion. 
Since we have already shown that
$\oar{\alpha_\lambda}: (\C_\lambda,\D_\lambda)\rightarrow (\C,\D)$ is a
regular map, the proof of \eqref{indlim4} is complete.

  \eqref{indlim4a} Part~\eqref{indlim4} shows $(\C,\D)$ is a \pd.
Clearly, for each $\lambda$, $(\C_\lambda, \D_\lambda\ms  \id|_{\C_\lambda})$ is
a Cartan envelope for $(\C_\lambda,\D_\lambda)$; let
$\Delta_\lambda:\C_\lambda\rightarrow \D_\lambda$ be the conditional
expectation.    For $\mu,
\lambda\in \Lambda$ with $\lambda\geq \mu$,
  Lemma~\ref{EWorka} gives 
\[\alpha_{\lambda\mu}\circ\Delta_\mu=
  \Delta_\lambda\circ\alpha_{\lambda\mu}.\]
This means that there exists $\Delta: \bigcup_{\lambda}
\oar{\alpha_\lambda}(\C_\lambda)\rightarrow \bigcup_{\lambda}
\oar{\alpha_\lambda}(\D_\lambda)$ such that
$\Delta\circ\oar{\alpha_\lambda}=\oar{\alpha_\lambda}\circ \Delta_\lambda$.   Since each $\Delta_\lambda$ is
surjective, contractive and idempotent, $\Delta$ also has these
properties. Therefore, $\Delta$ extends to a conditional expectation
$\Delta:\C\rightarrow\D$.   Proposition~\ref{Car+upse} shows $(\C,\D)$
is a Cartan inclusion.  Lemma~\ref{EWorka} gives
$\Delta\circ\oar{\alpha_\lambda}= \oar{\alpha_\lambda}\circ\Delta_\lambda$.
\end{proof}

We next show that the Cartan envelope of an inductive limit of a
system of \pd s
is the inductive limit of the Cartan envelopes.
  
\begin{theorem}\label{indlimEnv}
 Suppose $\Lambda$ is a directed
        set and $\{(\C_\lambda,\D_\lambda\ms  \alpha_{\lambda\mu}):
        \lambda\geq \mu\}$ is a system of \pd s.
        By Theorem~\ref{indlim},
        $(\indlim\C_\lambda,\indlim\D_\lambda)$ is a  \pd.

For each $\lambda$, let $(\A_\lambda,\B_\lambda\ms \tau_\lambda)$
  be a Cartan envelope for $(\C_\lambda,\D_\lambda)$.
  Theorem~\ref{cmap} provides unique regular $*$-monomorphisms 
  $\mac\alpha_{\lambda\mu} :(\A_\mu,\B_\mu)\rightarrow (\A_\lambda,\B_\lambda)$  ($\lambda\geq
  \mu$) making the diagram~\eqref{myCD1} commute; in particular
  $\mac\alpha_{\lambda\mu}$ satisfy
    \begin{equation}\label{indlim2.751}
      \mac\alpha_{\lambda\mu}\circ\tau_\mu=\tau_\lambda\circ\alpha_{\lambda\mu}.
    \end{equation}
The following statements hold.
\begin{enumerate}
\item \label{indlim1.1}  $\{(\A_\lambda,\B_\lambda\ms  \mac\alpha_{\lambda\mu}):
  \lambda\geq\mu\}$ is a system of Cartan inclusions and
  \[(\indlim\A_\lambda, \indlim\B_\lambda)\] is a Cartan inclusion.
\item\label {indlim3.14}
There is a regular $*$-monomorphism
  $\tau: \indlim\C_\lambda \rightarrow \indlim \A_\lambda$ such
  that for every $\lambda\in\Lambda$, \begin{equation}\label{indlim3.1415}\oar{\mac\alpha_\lambda}\circ
    \tau_\lambda=\tau\circ\oar{\alpha_\lambda}
  \end{equation}
  and
  $(\indlim \A_\lambda,\indlim\B_\lambda\ms \tau)$ is a Cartan envelope
  for $(\indlim\C_\lambda,\indlim\D_\lambda)$.
\end{enumerate} 
\end{theorem}
\begin{proof}
  Throughout the proof we will sometimes simplify notation and write, 
  \begin{align*}\C&:=\indlim\C_\lambda\dstext{and}\D:=\indlim\D_\lambda.\\
    \intertext{Similarly, once part~\eqref{indlim1.1} is established, we
    will sometimes write}
    \A&:=\indlim\A_\lambda\dstext{and}\B:=\indlim\B_\lambda.
  \end{align*}
  \eqref{indlim1.1}    Theorem~\ref{cmap} 
  shows that for $\lambda\geq \mu$, $(\A_\lambda, \B_\lambda\ms
  \mac\alpha_{\lambda\mu})$ is a regular and essential expansion of
    $(\A_\mu, \B_\mu)$.    Now suppose $\lambda\geq \mu\geq \nu$.
    We have 
    \[ \mac\alpha_{\lambda\nu}\circ
      \tau_\nu=\tau_\lambda\circ\alpha_{\lambda\nu}\] so
    \[(\mac\alpha_{\lambda\mu}\circ\mac\alpha_{\mu\nu})\circ\tau_{\nu}
      =\mac\alpha_{\lambda\mu}\circ\tau_\mu\circ\alpha_{\mu\nu}
      =\tau_\lambda\circ
      \alpha_{\lambda\mu}\circ\alpha_{\mu\nu}=\tau_\lambda\circ\alpha_{\lambda\nu}.\]
    Similar considerations yield,
    \[(\mac\alpha_{\lambda\mu}\circ\mac\alpha_{\mu\nu})\circ\Delta_\nu=\Delta_\lambda\circ
      (\mac\alpha_{\lambda\mu}\circ\mac\alpha_{\mu\nu})\dstext{and}
      (\mac\alpha_{\lambda\mu}\circ\mac\alpha_{\mu\nu}) \circ
      \Delta_\nu\circ\tau_\nu= \Delta_\lambda\circ\tau_\lambda\circ 
      \alpha_{\lambda\nu}. \] 
    By the uniqueness part of Theorem~\ref{cmap},
    we find
$\mac\alpha_{\lambda\nu}
=\mac\alpha_{\lambda\mu}\circ\mac\alpha_{\mu\nu}$, so
\[\{(\A_\lambda,\B_\lambda\ms\mac\alpha_{\lambda\mu}) : \lambda\geq\nu\}\]
  is a system of Cartan inclusions.   
  Theorem~\ref{indlim}\eqref{indlim4a} shows $(\indlim
  \A_\lambda,\indlim\B_\lambda)$ is a Cartan inclusion.

  \eqref{indlim3.14}
 By properties of inductive limits and \eqref{indlim2.751}, there is a
 uniquely determined  $*$-monomorphism \[\tau:
 \C\rightarrow \A\]  
satisfying
\eqref{indlim3.1415}.

For $\lambda\in \Lambda$, we have
\[(\tau\circ\oar{\alpha_\lambda})(\D_\lambda)=
 (\oar{\mac\alpha_\lambda}\circ\tau_\lambda)(\D_\lambda)\subseteq
 \oar{\mac\alpha_\lambda} (\B_\lambda)\subseteq \B,\] from which it
follows that $\tau(\D)\subseteq \B$.   Thus $(\A,\B\ms\tau)$ is
an expansion of $(\C,\D)$.    We next show that $(\A,\B\ms\tau)$ is an
essential expansion.

Let $J\idealin\B$ satisfy $J\cap \tau(\D)=\{0\}$.
Then \[J\cap(\oar{\mac\alpha_\lambda}\circ\tau_\lambda(\D_\lambda))=
  J\cap (\tau\circ\oar{\alpha_\lambda}(\D_\lambda))\subseteq
  J\cap\tau(\D)=\{0\}.\]  As $(\B_\lambda, \D_\lambda, \tau_\lambda)$
has the \iip, we obtain
$J\cap(\oar{\mac\alpha_\lambda}\circ\tau_\lambda(\B_\lambda))
=J\cap (\tau\circ\oar{\alpha_\lambda}(\B_\lambda)) =\{0\}.$
By~\cite[Proposition~II.8.2.4]{BlackadarOpAl},
    $\bigcup_{\lambda}(J\cap
    (\tau\circ\oar{\alpha_\lambda}(\B_\lambda)))$ is dense in $J$, we
    obtain $J=\{0\}$.  Thus, $(\A,\B\ms\tau)$ is an essential and Cartan
    expansion of $(\C,\D)$.

    Let \[N:=\overline{\bigcup_\lambda
    \oar{\alpha_\lambda}(\N(\C_\lambda,\D_\lambda))}.\]  Then $N$ is a
$*$-semigroup containing $\overline{\bigcup_\lambda
  \oar{\alpha_\lambda}(\D_\lambda)}=\D$.  An application of
Proposition~\ref{EWork} shows that $\tau:
\C\rightarrow \A$ is a regular map.

Next, note that for $\lambda\in \Lambda$,
\[(\Delta\circ\tau\circ\oar{\alpha_\lambda})(\C_\lambda)
  \stackrel{\eqref{indlim3.1415}}{=}
  (\Delta\circ\oar{\mac\alpha_\lambda}\circ\tau_\lambda)
  (\C_\lambda)
  \stackrel{\ref{indlim}\eqref{indlim4a}}{=}
  (\oar{\mac\alpha_\lambda}\circ\Delta_\lambda\circ\tau_\lambda)(\C_\lambda),.\]
Since $(\A_\lambda,\B_\lambda\ms\tau_\lambda)$ is a Cartan envelope
for $(\C_\lambda,\D_\lambda)$, we have
$C^*(\Delta_\lambda(\tau_\lambda(\C_\lambda)))=\B_\lambda$, so 
\[C^*((\Delta\circ\tau\circ\oar{\alpha_\lambda})(\C_\lambda))
  =\oar{\mac\alpha_\lambda}(\B_\lambda).\]   Hence
\[\B=C^*(\Delta(\tau(\C))).\]  Similar considerations yield
\[\A=C^*(\B\cup\tau(\C)).\]   By Definition~\ref{defCarEnv},
$(\A,\B\ms\tau)$ is a Cartan envelope for $(\C,\D)$.  
\end{proof}

\subsection{Minimal Tensor Products}\label{sec:tenprod}
The purpose of this subsection  is to show that the minimal tensor product of two \pd
s is again a \pd\ and that the process of taking Cartan envelopes
commutes with  minimal tensor products; see Theorem~\ref{gentp} below for the
precise statement.  We thank \referee\ for suggesting we use 
\cite[Lemma~5.1]{BarlakLiCaSuUCTPr}; doing so gives a much
more direct route to Theorem~\ref{gentp} than our original path.

   Throughout, for $i=1,2$, let $(\C_i,\D_i)$ be \pd s.  For algebras
   $\A$ and $\B$, we use the notation $\A\odot \B$ for the linear span
   of $\{a\otimes b: a\in \A, b\in \B\}$. 

   Regularity of the tensor product inclusion is straightforward, and
 recorded in the following lemma.
\begin{lemma}\label{tensorreg} Let $\S:=\{v_1\otimes v_2:
      v_i\in\N(\C_i,\D_i)$.  Then $\S\subseteq \N(\C_1\otimes_{\min} \C_2,\D_1\otimes\D_2)$ and $\spn\S$
      is dense in $\C_1\otimes_{\min} \C_2$.  In particular,
      $(\C_1\otimes_{\min} \C_2,\D_1\otimes\D_2)$ is a regular inclusion.
    \end{lemma}
    \begin{proof}
Let $\fA:=\spn\S$.

      That $\S\subseteq \N(\C_1\otimes_{\min} \C_2,\D_1\otimes\D_2)$
      is evident.   For $i=1,2$, suppose $x_i\in
    \C_i$.  We wish to show that $x_1\otimes x_2\in \overline\fA$.
    Let $\eps>0$.  By regularity of $(\C_i,\D_i)$, there exists $N\in
    \bbN$,  $\{v_j\}_{j=1}^N \subseteq \N(\C_1,\D_1)$ and
  $\{w_j\}_{j=1}^N\subseteq 
    \N(\C_2,\D_2)$ such that
  \[\norm{x_1-\sum_{j=1}^N v_j}<\eps \dstext{and}
    \norm{x_2-\sum_{j=1}^N w_j}<\eps.\]  Then
\begin{align*}\norm{x_1\otimes x_2-\left(\sum_{j=1}^N
      v_j\right)\otimes\left(\sum_{j=1}^N w_j\right)}&\leq
  \norm{\left(x_1-\sum_{j=1}^N
      v_j\right)\otimes x_2} \\& + \norm{\left(\sum_{j=1}^N
      v_j\right)\otimes \left(x_2-\sum_{j=1}^N
                                                       w_j\right)}\\
  &<\eps\norm{x_2}+(\norm{x_1}+\eps)\eps.
\end{align*} Therefore, $x_1\otimes x_2\in \overline\fA$.  As the
span of elementary tensors is dense in $\C_1\otimes_{\min}\C_2$, 
the lemma follows. 
\end{proof}

 We next observe that tensor products of Cartan inclusions are again
 Cartan.   The following is
essentially a reformulation of a result of Barlak and Li. 
\begin{lemma}[{\cite[Lemma~5.1]{BarlakLiCaSuUCTPr}}]\label{mtpCar}
For $i=1,2$, let $(\A_i,\B_i)$ be Cartan inclusions.  Then
$(\A_1\mintp \A_2,\B_1\otimes\B_2)$ is a Cartan inclusion.
Furthermore, if $\Delta_i: \A_i\rightarrow \B_i$ is the conditional
expectation, then $\Delta_1\otimes\Delta_2: \A_1\mintp\A_2\rightarrow
\B_1\otimes\B_2$ is the conditional expectation. 
\end{lemma}
\begin{proof}
  By Renault's theorem (as extended by Raad in~\cite{RaadGeReThCaSu})
  associated to every Cartan inclusion is a twist,
  \begin{equation*}
  \unit G\times \bbT\hookrightarrow \Sigma\twoheadrightarrow G,
\end{equation*} where $G$ is a
  locally compact, Hausdorff, \' etale and effective 
  groupoid; furthermore,  given such a twist,  $(C^*_r(\Sigma,G), C_0(\unit G))$
  is a  Cartan inclusion.
Thus, we may as well assume that for $i=1,2$, 
  \[(\A_i,\B_i)=(C^*_r(\Sigma_i,G_i), C_0(\unit G_i)),\] where 
$\unit G_i\times \bbT\hookrightarrow \Sigma_i\twoheadrightarrow G_i$
are twists over the
  locally compact, Hausdorff, \' etale and effective
  groupoids $G_i$.

As in~\cite[Section~5]{BarlakLiCaSuUCTPr}, let $G=G_1\times G_2$ and 
$\Sigma=(\Sigma_1\times\Sigma_2)/\sim$, where for every $z\in\bbT$ and
$\sigma_j\in \Sigma_j$, 
$(z\sigma_1,\sigma_2)\sim (\sigma_1, z\sigma_2)$.   Then $G$ is a  Hausdorff, \'etale,
locally compact, and 
topologically principal groupoid and 
\begin{equation*}
  \unit G\times \bbT\hookrightarrow \Sigma\twoheadrightarrow G
\end{equation*} is a  twist; thus 
$(C^*_r(\Sigma, G), C_0(\unit G))$ is a Cartan inclusion.  The proof of 
\cite[Lemma~5.1]{BarlakLiCaSuUCTPr}  does not require that $G_i$ be
second countable.   As
\cite[Lemma~5.1]{BarlakLiCaSuUCTPr} gives a formula for an isomorphism
of \[(C^*_r(\Sigma_1, G_1)\mintp C^*_r(\Sigma_2, G_2), C_0(\unit
G_1)\otimes C_0(\unit G_2))\dstext{onto} (C^*_r(\Sigma), C_0(\unit
G)),\] $(\A_1\mintp\A_2, \D_1\otimes\D_2)$ is a Cartan inclusion.

One can obtain the last statement using  groupoid models, but
we shall instead obtain it as follows.  By~\cite[IV.4.23(i)]{TakesakiThOpAlI}, the  map $x_1\otimes x_2 \mapsto \Delta_1(x_1)\otimes \Delta_2(x_2)$ uniquely
extends to a completely positive map $\Delta_1\otimes \Delta_2:
\A_1\mintp\A_2\rightarrow \B_1\otimes \B_2$.   Noting that $\Delta_1\otimes
\Delta_2$ fixes $\B_1\otimes\B_2$, we
see that $\Delta_1\otimes \Delta_2$ is the conditional expectation
for $(\A_1\mintp\A_2, \B_1\otimes\B_2)$. 
\end{proof}

Dual to the notion of an essential embedding of $C(X)$ into $C(Y)$ is
the notion of essential surjection of $Y$ onto $X$:  
a continuous surjection $\pi: Y\twoheadrightarrow X$ between compact
Hausdorff spaces is called an \textit{essential surjection} if whenever
$F\subseteq Y$ is a closed set and $\pi(F)=X$, then $F=Y$.  Essential
surjections are also called \textit{irreducible surjections}.

             \begin{lemma}\label{esscovpr}  Suppose for $i=1, 2$, $Y_i$ and $X_i$
  are compact Hausdorff spaces, and $\pi_i: Y_i\twoheadrightarrow X_i$
  are continuous and essential surjections.   Then $\pi_1\times\pi_2:
  Y_1\times Y_2\twoheadrightarrow X_1\times X_2$, given by $(y_1,y_2)\mapsto
  (\pi_1(y_1),\pi_2(y_2))$,  is an essential
  surjection.
\end{lemma}
\begin{proof}  
We start by proving a special case:  assume that $Y_2=X_2$ and $\pi_2$
is the identity map.

For typographical ease, let $\pi:=\pi_1\times \id_{X_2}$.
Let $F\subseteq Y_1\times X_2$ be a closed set such that
$\pi(F)=X_1\times X_2$.   For $t\in X_2$,  let
\[F_t:=\{y_1\in Y_1: (y_1,t)\in F\}.\] Then $F_t$ is a closed subset
of $Y_1$.  Given $s\in X_1$, since $\pi(F)=X_1\times X_2$, there is
$(y_1,t)\in F$ such that $\pi(y_1,t)=(s,t)$.  Thus $\pi_1(F_t)=X_1$.
Since $\pi_1$ is essential, we obtain $F_t=Y_1$.  Hence for $t\in X_2$,
$Y_1\times\{t\}\subseteq F$. We conclude that $F=Y_1\times X_2$,
showing that the lemma holds in this special case.

For the general case, consider the composition, 
\[Y_1\times Y_2\stackrel{\pi_1\times \id_{Y_2}}{\longrightarrow}
    X_1\times Y_2 \stackrel{\id_{X_1}\times \pi_2}{\longrightarrow}
    X_1\times  X_2.\]  By the special case, each of the maps in the
  composition are
  essential surjections.  As the composition of essential surjections
  is an essential surjection, the proof is complete.
\end{proof}

\begin{proposition} \label{tensorpe}  For $i=1,2$, let $\B_i$ be 
  abelian \cstaralg s (perhaps not unital) and suppose $(\B_i,\D_i,\alpha_i)$ are
  essential inclusions.   Then $(\B_1\otimes\B_2,\D_1\otimes\D_2,
  \alpha_1\otimes\alpha_2)$ is an essential inclusion. 
\end{proposition}
\begin{proof}
First assume the inclusions $(\B_i,\D_i,\alpha_i)$ are unital.
For $i=1,2$, let $Y_i:=\widehat\B_i$ and $X_i:=\widehat\D_i$.    The
maps $\alpha_i$ dualize to continuous surjections $\pi_i:
Y_i\rightarrow X_i$ (thus for $d\in \D_i$,
$\pi_i(y_i)(d)=y_i(\alpha_i(d))$).  The proof in this case now follows
from Lemma~\ref{esscovpr} after noting
that the  inclusion $(\B_i, \D_i,\alpha_i)$ is essential if and only if $\pi_i$ is an essential surjection
of $Y_i$ onto $X_i$, see \cite[Lemma~4.9]{PittsIrMaIsBoReOpSeReId}.

For the general case, note that
$(\tilde\D_1\otimes\tilde\D_2,\D_1\otimes\D_2, u_{\D_1}\otimes
u_{\D_2})$ is an essential inclusion because the image of
$\D_1\otimes\D_2$ under $u_{\D_1}\otimes u_{\D_2}$ is an essential
ideal in $\tilde\D_1\otimes\tilde\D_2$.  (This can be shown directly
or one can use the facts that
$\tilde\D_1\otimes\tilde\D_2\subseteq M(\D_1)\otimes M(\D_2)\subseteq
M(\D_1\otimes\D_2)$ and any \cstaralg\ is an essential ideal in its
multiplier algebra.  A proof of the inclusion of the multiplier
algebras can be found at
\texttt{https://math.stackexchange.com/a/4458451}.)  This fact,
together with the unital case show that both of the inclusions,
\[\alpha_1(\D_1)\otimes\alpha_2(\D_2)\subseteq \tilde\alpha_1(\tilde\D_1)\otimes\tilde\alpha_2(\tilde\D_2)\subseteq
  \tilde\B_1\otimes\tilde\B_2\] have the \iip.   By
Lemma~\ref{interiip}, $\alpha_1(\D_1)\otimes\alpha_2(\D_2)\subseteq
\tilde\B_1\otimes\tilde\B_2$ has the \iip. Finally, since
$\alpha_1(\D_1)\otimes\alpha_2(\D_2)\subseteq \B_1\otimes\B_2\subseteq
\tilde\B_1\otimes\tilde\B_2$, another application of
Lemma~\ref{interiip} shows
$\alpha_1(\D_1)\otimes\alpha_2(\D_2)\subseteq \B_1\otimes\B_2$ has the
\iip, as desired.
\end{proof}

  \begin{theorem}\label{gentp}  Suppose for $i=1,2$ that $(\C_i,\D_i)$
    are \pd s (not assumed unital) and let $(\A_1,\B_i\ms \alpha_i)$ be
    Cartan envelopes for $(\C_i,\D_i)$.  Then
    $(\C_1\mintp\D_2,\D_1\otimes\D_2)$ is a \pd, and
    $(\A_1\mintp\A_2,\B_1\otimes\B_2\ms  \alpha_1\otimes\alpha_2)$ is a
    Cartan envelope for $(\C_1\mintp\D_2,\D_1\otimes\D_2)$.
  \end{theorem}
  \begin{proof}  
By    Lemma~\ref{mtpCar},  $(\A_1\mintp\A_2, \B_1\otimes\B_2)$ is a
    Cartan inclusion.  As
    $(\alpha_1\otimes\alpha_2)(\D_1\otimes\D_2)\subseteq
    \B_1\otimes\B_2$, an application of  Proposition~\ref{tensorpe}
    shows  
    $(\A_1\mintp\A_2, \B_1\otimes\B_2\ms\alpha_1\otimes\alpha_2)$ is a
    Cartan and 
    essential expansion of $(\C_1\mintp\C_2, \D_1\otimes\D_2)$.

    Since Cartan inclusions have the faithful unique
    pseudo-expectation property (see Proposition~\ref{Car+upse}), Lemma~\ref{claim2}\eqref{claim2:a}
    combined  with Lemma~\ref{tensorreg} shows that $(\C_1\mintp\C_2,
    \D_1\otimes\D_2)$ is a regular inclusion with the faithful unique
    pseudo-expectation property, that is, $(\C_1\mintp\C_2,
    \D_1\otimes\D_2)$  is a \pd.
    
      Applying Proposition~\ref{EWork} with
    $N=\N(\C_1,\D_1)\otimes\N(\C_2,\D_2)$ shows
    $\alpha_1\otimes\alpha_2$ is a regular $*$-monomorphism.
    It remains to show that $(\C_1\mintp\C_2,\D_1\otimes\D_2)$ generates
    $(\A_1\mintp\A_2, \B_1\otimes\B_2)$ in the sense of
    Definition~\ref{defCarEnv}\eqref{defCarEnv1}(ii). 

    For $i=1,2$, let $\Delta_i: \A_i\rightarrow \B_i$ be the conditional
    expectation.     Since $\B_i=\Delta_i(\alpha_i(\C_i))$, we find
    \[B_1\otimes\B_2=C^*((\Delta_1\otimes\Delta_2)\circ(\alpha_1\otimes\alpha_2)(\C_1\mintp
      \C_2)).\]  Similarly, as $\A_i=C^*(\B_i\cup \alpha_i(\C_i))$,
 the $*$-subalgebra generated by $(\B_1\otimes \B_2) \bigcup
 (\alpha_1\otimes\alpha_2)(\C_1\mintp \C_2)$ is dense in
 $\A_1\mintp \A_2$.  Therefore, $(\A_1\mintp\A_2, \B_1\otimes\B_2\ms
 \alpha_1\otimes\alpha_2)$ is a Cartan envelope for $(\C_1\mintp\C_2, \D_1\otimes\D_2)$.
\end{proof}


\mysec{Applications}\label{apps}
\numberwithin{equation}{section}

In this section we give a few applications of the results presented so
far.  We begin with an application of Theorem~\ref{!pschar} which describes
the  \cstar-envelope for an 
intermediate Banach algebras.  Several of the results in this section
apply to unital inclusions, see Definition~\ref{incdef}\eqref{innore2}.

\begin{theorem}\label{c*env}  Suppose $(\C,\D)$ is a unital \pd\ 
  and $\A\subseteq \C$ is a closed
  subalgebra satisfying
  $\D\subseteq \A\subseteq \C$  (we do not assume $\A=\A^*)$.  Let $C^*(\A)$ be the
  \cstar-subalgebra of $\C$ generated by $\A$ and let $C^*_e(\A)$ be
  the \cstar-envelope for $\A$.  Then
  \[C^*(\A)=C^*_e(\A).\]
\end{theorem}

\begin{proof}
Since $(\C,\D)$ has the faithful unique pseudo-expectation property,
the result follows from~\cite[Theorem~8.3]{PittsStReInI}.
\end{proof}

For unital
\cstaralg s $\A$ and $\B$ with $\B\subseteq \A$, when the assumption
that $\B$ is abelian and regular is dropped,
~\cite[Example~6.9]{PittsZarikianUnPsExC*In} shows that the faithful unique
pseudo-expectation property for $\B\subseteq \A$ need not
imply that $\B$ norms $\A$ in the sense
of~\cite{PopSinclairSmithNoC*Al}. Nevertheless,
\cite[Section~6.2]{PittsZarikianUnPsExC*In} argues that the faithful unique
pseudo-expectation property is ``conducive'' to norming.   Further
evidence for this statement is 
Theorem~\ref{DnC} below: for unital regular inclusions
satisfying Standing Assumption~\ref{SA1}, the faithful unique
pseudo-expectation property  implies norming.   

The proof of Theorem~\ref{DnC} requires two preparatory results which
we now present.  We shall be using the notions of topological freeness
and free points for an inclusion, see Definition~\ref{fdef}}.

\begin{lemma}\label{GNS0}  Let $(\C,\D)$ be a unital \pd. 
Suppose $\sigma\in \hat\D$ is free and let $f$ be the  unique state
extension of $\sigma$ to $\C$.
If 
$(\pi_f,\H_f,\xi_f)$ is the GNS triple  associated to $f$, then 
$\pi_f(\D)''$ is an atomic MASA in $\B(\H_f)$.

Furthermore,
if \[\O:=\{\tau\in\hat\D: \exists\, v\in \N(\C,\D) \text{ such that }
\tau=\beta_v(\sigma)\}\] is the orbit of $\sigma$ under $\N(\C,\D)$,
then 
\begin{equation}\label{GNS0k}
  \ker(\pi_f|_\D)=\bigcap_{\tau\in\O} \ker\tau.
\end{equation}
\end{lemma}
\begin{proof}
  Our first task is to show that $f$ is a compatible state on $\C$
  (\cite[Definition~4.1]{PittsStReInI}).  Let $(I(\D),\iota)$ be an
  injective envelope for $\D$ and let $E: \C\rightarrow I(\D)$ be the
  pseudo-expectation for $(\C,\D)$.  Choose
  $\sigma'\in \widehat{I(\D)}$ such that $\sigma'\circ \iota=\sigma$.
  Note that $\sigma'\circ E|_\D=\sigma'\circ\iota=\sigma$.  Since
  $\sigma$ is free, we conclude that $f=\sigma'\circ E$.
  By~\cite[Theorem~6.9]{PittsStReInII} (see Theorem~\ref{upse=>cov} below)
  $f$ is a
  compatible state.
  
  The Cauchy-Schwartz inequality implies that for any $d\in \D$ and
  $x\in\C$, $f(dx)=\sigma(d)f(x)=f(xd)$
  (\cite[Fact~1.2]{PittsNoApUnInC*Al} has the details), so $f\in\text{Mod}(\C,\D)$
  (see~\cite[Definition~2.4]{PittsStReInI}).

We  claim that if $v\in \N(\C,\D)$ satisfies 
$\sigma(v^*v)=1$ and $\sigma(v^*dv)=\sigma(d)$ for every $d\in\D$, then 
\[f(v^*xv)=f(x) \dstext{for all} x\in \C. \] Define a positive linear
functional $f_v$ on $\C$ by $f_v(x)=f(v^*xv)$; note that 
$f_v(I)=1$.   A standard fact (e.g.~\cite[Lemma~1.5.15]{BrownOzawaC*AlFiDiAp}) shows $f_v$ is a state on
$\C$. By hypothesis, $f_v$ extends $\sigma$.   Freeness of $\sigma$ gives
$f_v=f$, so the claim holds.

It follows that $f$ is a pure and  $\D$-rigid state
(see~\cite[Definition~4.11]{PittsStReInI}).   By first
applying~\cite[Proposition~4.12]{PittsStReInI} and
then~\cite[Proposition~4.8(iv)]{PittsStReInI}, we find  
$\pi_f(\D)''$ is an atomic MASA in $\B(\H_f)$.

Let $N_f=\{x\in \C: f(x^*x)=0\}$ be the left kernel of $f$.
Suppose $d\in \D$ and $\pi_f(d)\neq 0$.   Since $(\C,\D)$ is a regular
inclusion, there exists $v\in
\N(\C,\D)$ such that $\pi_f(d)(v+N_f)\neq 0$.  Then $f(v^*v)\neq 0$
because $0\neq
f(v^*d^*dv)\leq \norm{d}^2f(v^*v)$.  As
$f|_\D=\sigma$,~\cite[Proposition~4.4(v)]{PittsStReInI}, shows 
$\beta_v(\sigma)(d)\neq 0$.  As $\beta_v(\sigma)\in \O$, we find
\[\bigcap_{\tau\in \O}\ker \tau\subseteq \ker(\pi_f|_\D).\]  On the other hand,
if $d\in \ker(\pi_f|_\D)$,  then for every $v\in \N(\C,\D)$ with
$\sigma(v^*v)\neq 0$, we have $0=\pi_f(d)(v+N)$, so applying
\cite[Proposition~4.4(v)]{PittsStReInI} again,
$\beta_v(\sigma)(d)=0$.  Therefore, $d\in \bigcap_{\tau\in \O} \ker
\tau$.
Thus~\eqref{GNS0k} holds and the proof is complete.
\end{proof}

Our second preparatory result, of independent interest, shows that a
natural class of \pd s is topologically free.  
\begin{proposition}\label{cgen0}  Suppose $(\C,\D)$ is a  unital \pd\
  such that there exists
a countable subset $\fC\subseteq\N(\C,\D)$ satisfying 
 $\C=C^*(\fC\cup \D)$.  Then $(\C,\D)$ is a topologically free
 inclusion. 
\end{proposition}

\begin{proof}
Let $\fU:=\{\sigma\in \hat\D: \sigma \text{ is
free for $(\C,\D)$}\}$.  Our task is to show $\fU$ is dense in $\hat\D$.
For convenience, set \[\B:=\D^c\] and let $\pi:\hat\B\rightarrow \hat\D$ be the restriction map,
      \[\pi(\sigma_1)= \sigma_1|_\D\qquad
        (\sigma_1\in\hat\B).\] 

      Since $\fC$ is countable, so is the $*$-semigroup consisting of
      all finite products with factors taken from $\fC\cup\fC^*$.  Thus
      we  may assume without loss of generality that $\fC$ is a
      countable $*$-semigroup of normalizers.

  With the notation of~\eqref{RnS} and~\eqref{J_c}, let $\Q:=\{vh: v\in\fC,
  h\in S(v)_c\}$,  $\S:=\Q\cup\D$ and let
  \[\fM=\spn \S.\] 
  Proposition~\ref{*semigp} shows $\fM$ is
  a dense $*$-subalgebra of $\C$.

Let $(\C_1,\D_1\ms \alpha)$ be a Cartan envelope for $(\C,\D)$ and let
$\Delta: \C_1\rightarrow \D_1$ be the conditional expectation. 
Then
\[\D_1=C^*(\Delta(\alpha(\C)))=C^*(\Delta(\alpha(\fM)))\subseteq C^*(\Delta(\alpha(\fC))\cup
  \alpha(\D))\subseteq \D_1.\]  Therefore
$\D_1=C^*(\Delta(\alpha(\fC)) \cup \alpha(\D)).$

  Notice that for $v\in \fC$ and $h\in S(v)_c$,
  $\Delta(\alpha(vh))=\Delta(\alpha(v)) \alpha(h)\in
  C^*(\Delta(\alpha(\fC))\cup\alpha(\D))$,  whence
  \[\Delta(\alpha(\M))\subseteq
  C^*(\Delta(\alpha(\fC))\cup\alpha(\D)).\]

As $\Delta(\alpha(\fC))$ is countable and $(\D_1,\alpha(\D))$ has the
\iip, we may apply \cite[Theorem~2.11.16]{ExelPittsChGrC*AlNoHaEtGr}
to conclude that the set of free points for $(\D_1,\alpha(\D))$
contains a dense $G_\delta$ set.  By
Lemma~\ref{claim2}\eqref{claim2:b},
\[\alpha(\D)\subseteq \alpha(\B)\subseteq \D_1.\] As 
the set of
      free points for the inclusion $(\D_1, \alpha(\D))$ is contained the
      set of free points for $(\alpha(\B),\alpha(\D))$, we see that the set of free
      points for $(\B,\D)$  contains a
      dense $G_\delta$ set.
    
           We may therefore choose a countable family  $U_n\subseteq
           \hat\D$ of dense open sets such that 
      \[\bigcap_n U_n\] is contained in the set of free points for
      $(\B,\D)$.  As $\D\subseteq \B$ has the \iip, $\pi$ is an
      essential surjection (the paragraph
      preceding Lemma~\ref{esscovpr} has the definition of this term), so for each $n$,  $\pi^{-1}(U_n)$ is a
      dense open subset of $\hat\B$.

Let $v\in\fC$. By Lemma~\ref{relcom}\eqref{relcom1}, 
$v\in \N(\C,\B)$.  Moreover, as $(\C,\B)$ is a regular MASA
inclusion, the proof of \cite[Theorem~3.10]{PittsStReInI} shows there
is a dense open set $X_v\subseteq \hat\B$ such that each
$\sigma_1\in X_v$ is free relative to $v$.  Let
\[\Y:=\left(\bigcap_{v\in \fC} X_v\right)\cap
  \left(\bigcap_n\pi^{-1}(U_n)\right).\] As $\B$ is unital, $\hat\B$
is compact, so Baire's theorem implies
$\Y\subseteq\hat\B$ is a dense $G_\delta$ set.

 We claim that if $y\in \Y$,   then $\pi(y)$ is a free point for
 $(\C,\D)$.  To see this, suppose
$\rho_1$ and $\rho_2$ are states on $\C$ such that
$\rho_1|_\D=\rho_2|_\D=\pi(y)$. 
  Since $y\in
          \bigcap_n \pi^{-1} (U_n)$,  we find
          $\pi(y)\in\bigcap_n U_n$.  Thus
          $\pi(y)$ 
          is a free point for the inclusion $(\B,\D)$.
Next, since $y\in \bigcap_{v\in\fC} X_v$, $\rho_1(v)=\rho_2(v)$ for
 every $v\in\fC$.  Using~\cite[Fact~1.2]{PittsNoApUnInC*Al},
 for every $v\in \fC$ and $d, e\in \D$,
$\rho_1(dve)=\rho_2(dve)$.  Therefore, for every
$x\in \fM$, $\rho_1(x)=\rho_2(x)$.  As $\fM$ is
dense in $\C$, we conclude $\rho_1=\rho_2$.  Thus $\pi(y)$ is a
free point for $(\C,\D)$, so the claim holds.
          
Finally we show that $\pi(\Y)$ is dense in $\hat\D$.  Let $G\subseteq
\hat\D$ be an open set.   Then $\pi^{-1}(G)$ is open in $\hat\B$, so
there exists $y\in \Y\cap \pi^{-1}(G)$, whence $\pi(y)\in G$, as desired.
As $\pi(\Y)\subseteq \fU$, the proof is complete.
\end{proof}

We now come to our norming result.

\begin{theorem}\label{DnC}
  Suppose $(\C,\D)$ is a unital \pd.
  Then $\D$ norms $\C$.
\end{theorem}
\begin{proof}
The  argument is mostly the same as the proof of
~\cite[Theorem~8.2]{PittsStReInI}, except we use Proposition~\ref{cgen0}
instead of~\cite[Theorem~3.10]{PittsStReInI} and Lemma~\ref{GNS0} instead
of~\cite[Proposition~4.12]{PittsStReInI}.    Therefore, we shall
only outline the proof, leaving the reader to consult the proof
of~\cite[Theorem~8.2]{PittsStReInI} for additional details as desired.  (The notation here
differs somewhat from that used in~\cite[Theorem~8.2]{PittsStReInI}, but
this will present no difficulty.)

Arguing as  in the last paragraph of the proof of
~\cite[Theorem~8.2]{PittsStReInI}, it suffices to show $\D$ norms
$\C$ under 
 the additional assumption that there is a countable set
$\fC\subseteq \N(\C,\D)$ such that $\C=C^*(\fC\cup \D)$.   For the
remainder of the proof, we assume this.

Let $\F\subseteq \N(\C,\D)$ be the $*$-semigroup generated by $\fC\cup
\D$. By the additional assumption, $\spn \F$ is a dense $*$-subalgebra
of $\C$.

Let $Y\subseteq \hat\D$ be the set of free points for  $(\C,\D)$.  Proposition~\ref{cgen0}
shows $Y$ is dense in $\hat\D$.  For $y\in Y$, denote by $f_y$ the
unique state extension of $y$ to $\C$, and let $(\pi_y, \H_y, \xi_y)$
be the GNS triple for $f_y$.   Since $f_y$ is a pure state, $\pi_y$ is
an irreducible representation, and Lemma~\ref{GNS0} shows that
$\pi_y(\D)''$ is a MASA in $\B(\H_y)$. 

For $y_1, y_2\in Y$, define $y_1\sim y_2$ if and only if there exists
$v\in\F$ such that $\beta_v(y_1)=y_2$.  As in the proof
of~\cite[Theorem~8.2]{PittsStReInI},   $\sim$ is an equivalence
relation, and  if $\pi_{y_1}$ is unitarily equivalent to
$\pi_{y_2}$, then $y_1\sim y_2$.  Thus, if $y_1\not\sim y_2$,
$\pi_{y_1}$ and $\pi_{y_2}$ are disjoint representations.   

Let $\Y\subseteq Y$ be a set containing exactly one element of
each $\sim$ equivalence class.  Consider the representation  
\[\pi:=\bigoplus_{y\in\Y} \pi_y\] of $\C$ 
on $\H_\pi:=\bigoplus_{y\in\Y} \H_y$.  Since each $\pi_y(\D)''$ is an
atomic MASA and the representations $\{\pi_y:y\in \Y\}$ are pairwise
disjoint, $\pi(\D)''$ is also an atomic MASA in $\B(\H_\pi)$.

We now show $\pi$ is faithful. As
$Y$ is dense in $\hat\D$ and is the union of the $\N(\C,\D)$-orbits of the elements of $\Y$,~\eqref{GNS0k} shows that the
restriction of $\pi$ to $\D$ is faithful.  By the ideal intersection
property for $(\D^c, \D)$, $\ker\pi \cap \D^c=\{0\}$.  Then
$\ker\pi=\{0\}$ because $(\C,\D^c)$ also has the ideal intersection
property.  Therefore $\pi$ is a faithful representation of $\C$.

  A direct argument or the argument found
  on~\cite[Page~466]{CameronPittsZarikianBiCaMASAvNAlNoAlMeTh} shows
  $\pi(\D)''$ is locally cyclic
  (as defined in~\cite[Page~173]{PopSinclairSmithNoC*Al}) for
  $\B(\H_\pi)$.  Using~\cite[Lemma~2.3  and
  Theorem~2.7]{PopSinclairSmithNoC*Al}, $\pi(\D)$ norms $\B(\H_\pi)$.
  As $\pi$ is faithful, it follows that $\D$ norms $\C$.
\end{proof}

    We now extend~\cite[Theorem~8.4]{PittsStReInII} from the setting
    of virtual
    Cartan inclusions considered there to \pd s.  This is a
    significant generalization of~\cite[Theorem~8.4]{PittsStReInII}
    because it weakens the  hypothesis that $\C$ is unital and relaxes
    the condition that $\D$ be a MASA.    We remark
    that~\cite[Example~5.1]{BrownExelFullerPittsReznikoffInC*AlCaEm}
    gives an example of a Cartan inclusion $(\C,\D)$ and an
    intermediate \cstar-subalgebra $\D\subseteq \B\subseteq \C$ such
    that $(\B,\D)$ is not a regular inclusion.   Thus, in the context of
    Theorem~\ref{interAlgIso}, it is possible that the inclusions
    $(C^*(\A_i), \D_i)$ are not regular.
    \begin{theorem} \label{interAlgIso} Suppose for $i=1,2$,
      $(\C_i,\D_i)$ are 
      \pd s such that $(\tilde\C_i,\tilde\D_i)$ are \pd s, and $\A_i$ are Banach algebras satisfying
      $\D_i\subseteq \A_i\subseteq \C_i$.   Let $C^*(\A_i)$ be the
      \cstar-subalgebra of $\C_i$ generated by $\A_i$.   If
      $u:\A_1\rightarrow \A_2$ is an isometric isomorphism, then $u$
      uniquely extends to a $*$-isomorphism of $C^*(\A_1)$ onto
      $C^*(\A_2)$.
    \end{theorem}

    \begin{remark}{Remark}
We have
    included  the hypothesis
   that $(\tilde\C_i,\tilde\D_i)$ are \pd s  because we do not
    know whether regularity is preserved when units are adjoined.
When $(\tilde\C_i,\tilde\D_i)$ are  regular,
Observation~\ref{unpdiffnupd} shows $(\tilde\C_i,\tilde\D_i)$ are \pd s.
   As noted earlier, if $(\C_i,\D_i)$ have the
    AUP or if $\C_i$ are abelian, $(\tilde\C_i,\tilde\D_i)$ are \pd
    s provided that 
    $(\C_i,\D_i)$ are \pd s.   
  \end{remark}
    \begin{proof}[Proof of Theorem~\ref{interAlgIso}]

      Suppose first that $(\C_i,\D_i)$ are unital \pd s.
Theorem~\ref{c*env} shows that $C^*(\A_i)=C^*_{e}(\A_i)$.   By
Theorem~\ref{DnC}, $\D_i$ norms $\C_i$; in particular, $\D_i$ norms
$\A_i$.  Thus~\cite[Corollary~1.5]{PittsNoAlAuCoBoIsOpAl} shows $u$ uniquely extends to
a $*$-isomorphism of $C^*(\A_1)$ onto $C^*(\A_2)$.

Now suppose $(\C_i, \D_i)$ are not assumed unital and $u:
\A_1\rightarrow \A_2$ is an isometric isomorphism.   Since $(\C,\D)$
is \wnd,  
$(\tilde\C_i,\tilde\D_i)$ is a  unital \pd, so $I_{\tilde\D_i}=I_{\tilde\C_i}$.
Write $\tilde\A_i:=\A_i+\bbC I_{\tilde\C_i}$, so that 
$\tilde\D_i\subseteq \tilde\A_i\subseteq \tilde\C_i$. 
Notice that 
\[C^*(\tilde\A_i)=C^*(\A_i)+\bbC I_{\tilde\C_i}.\]
Applying~\cite[Corollary~3.3]{MeyerAdUnOpAl} with $n=1$ shows $\tilde
u:\tilde\A_1\rightarrow \tilde\A_2$  is an isometric isomorphism.
Therefore, $\tilde u$ uniquely extends to a $*$-isomorphism 
$\theta$ of $ C^*(\tilde \A_1)$ onto $C^*(\tilde\A_2)$.
Thus $\theta|_{C^*(\A_1)}$ is
a $*$-isomorphism of $C^*(\A_1)$ onto $C^*(\A_2)$ extending $u$.

If $\pi:
  C^*(\A_1)\rightarrow C^*(\A_2)$ is another $*$-isomorphism such that
  $\pi|_{\A_1}=u$, then  $\tilde\pi|_{\tilde\A_1}=\tilde u$.  As
  $\theta$ is the unique extension of $\tilde u$ to $C^*(\tilde\A_1)$, 
  $\tilde\pi=\theta$.
  Therefore, $\pi=\theta|_{C^*(\A_1)}$.  This gives uniqueness of the
  extension of $u$ to $C^*(\A_1)$
  and completes the proof.\end{proof}

If the hypothesis that $(C^*(\A_1),\D_1)$  is a regular inclusion is added to the hypotheses of Theorem~\ref{interAlgIso}, more can be said.

\begin{corollary} \label{interAlgIsoCor2} Suppose for $i=1,2$,
  $(\C_i,\D_i)$ are \pd s such that: $(\tilde\C_i,\tilde\D_i)$ are \pd
  s; $\A_i$ are Banach algebras satisfying
  $\D_i\subseteq \A_i\subseteq \C_i$; and $u:\A_1\rightarrow \A_2$ is
  an isometric isomorphism.

      If $(C^*(\A_1),\D_1)$ is regular, then $(C^*(\A_1),\D_1)$ and
      $(C^*(\A_2), u(\D_1))$ are \pd s; let $(\M_1,\N_1\ms \tau_1)$
      and $(\M_2,\N_2\ms\tau_2)$ be their Cartan envelopes.
      Furthermore, there is a unique extension of $u$ to a
      $*$-isomorphism $\mac u:\M_1\rightarrow\M_2$ such that
      the diagram~\eqref{myCD1} commutes; in particular,
      \[\mac u\circ \tau_1=\tau_2\circ u.\]
    \end{corollary}
    \begin{proof}
 Let $u':C^*(\A_1)\rightarrow C^*(\A_2)$ be the unique $*$-isomorphism extending $u$ provided by Theorem~\ref{interAlgIso}.   Since $u'|_{\D_1}=u|_{\D_1}$, $u(\D_1)$ is a \cstaralg.   Since $(C^*(\A_1),\D_1)$ is a regular inclusion, so is $(C^*(\A_2),u(\D_1))$.  
      By Corollary~\ref{heriditary}, both $(C^*(\A_1),\D_1)$ and $(C^*(\A_2), u(\D_1)$) have the faithful unique pseudo-expectation property, so they are \pd s.  Also for any $v\in \N(C^*(\A_1),\D_1)$ and $h\in \D_1$,
      \[u'(v)^*u(h)u'(v)=u'(v^*hv)\in u'(\D_1)=u(\D_1).\] Thus $u'$ is
      a regular $*$-isomorphism.   Since $(C^*(\A_2), u(\D_1)\ms u')$
      is an essential expansion of $(C^*(\A_1),\D_1)$,  an application of Theorem~\ref{cmap} completes the proof.
    \end{proof}

The following is immediate from Theorem~\ref{interAlgIso}.
\begin{corollary}\label{interAlgIsoCor}  Suppose $(\C,\D)$ is a \pd\
  such that $(\tilde\C,\tilde\D)$ is a \pd. If $\A$ is an algebra
  such that
  $\D\subseteq \A\subseteq \C$ and $C^*(\A)=\C$,  then the group of isometric
  automorphisms of $\A$ is isomorphic to the group of all
  $*$-automorphisms of $\C$ which leave $\A$ invariant.
\end{corollary}

Theorems~\ref{c*env}--\ref{interAlgIso} seem interesting even in the
commutative case, for they apply to certain function algebras.
\begin{example}\label{lh}
Let $\C$ be a unital and abelian \cstaralg, and suppose $J\idealin\C$
is an essential ideal.   Let $\A\subseteq \C$ be a Banach algebra such
that: $\A$ separates points of $\hat\C$,  $\bbC I_\C\subseteq \A$, and  $\A\cap J=\{0\}$.  Then $\A +J$ is a Banach algebra satisfying
$J\subseteq J+\A\subseteq  \C$.  By the Stone-Weierstrau\ss\ Theorem,
$C^*(\A)=\C$.     Since $(\C, J)$ is a \pd, 
$C^*_{env}(\A+J)=\C$ and any isometric automorphism of $\A+J$ uniquely
extends to a $*$-isomorphism of $\C$.  
\end{example}

\mysec[Questions]{Questions}\label{Sec:Ques}
\numberwithin{equation}{section}
In addition to Conjecture~\ref{nuclear}, we now present a few open
questions.

\begin{remark}{Question}\label{findCar}  Suppose $(\C,\D)$ is a \pd.  Must
  $\C$ contain a Cartan MASA? 

We expect the answer  is no.  Here is a possible
approach to showing this.  
Suppose $\Gamma$ is a discrete and $C^*$-simple group such
  that $C^*_r(\Gamma)$ has no Cartan MASA.
  (By~\cite[Corollary~4.11]{LiRenaultCaSuC*AlExUn}, the free group on
  $n\geq 2$ generators is such a group.)  If $\Gamma$ can be chosen so
  that $C^*_r(\Gamma)$ contains a regular MASA $\D$, then 
  $(\C^*_r(\Gamma),\D)$ is a pseudo-Cartan inclusion, giving a
  negative answer.  (To see $(\C^*_r(\Gamma),\D)$  is a
  \pd, note that it has the unique
  pseudo-expectation property by ~\cite[Theorem~3.5]{PittsStReInI}.
  But $\Gamma$ is \cstar-simple, so~\cite[Theorem~3.15]{PittsStReInI}
  shows $(C^*_r(\Gamma),\D)$ has the faithful unique
  pseudo-expectation property, whence it is a virtual Cartan
  inclusion.)  As an interesting aside, if  no MASA in $C^*_r(\Gamma)$ is
  regular, then by Fact~\ref{nopsCar0}, $C^*_r(\Gamma)$ would not
  contain a 
  pseudo-Cartan subalgebra.
\end{remark}

\begin{remark}{Question}\label{whendiag} Every \pd\ $(\C,\D)$ has a Cartan
  envelope $(\A,\B\ms\alpha)$.  What conditions on $(\C,\D)$ ensure
  that $(\A,\B)$ is a \cstar-diagonal?
\end{remark}

Recall that Barlak and Li~\cite[Corollary~1.2]{BarlakLiCaSuUCTPr}
showed that if a separable, nuclear, \cstaralg\ contains a Cartan
MASA, then it satisfies the Universal Coefficient Theorem (UCT).
Separable, nuclear \cstaralg s which satisfy the UCP are of significant
interest in the classification program.    This is in part  the motivation for
the next question.  Before stating the question, we give an example.

 Let $X$ be the
Cantor set and let $\Gamma$ be a countable discrete group having
property $(T)$.  It follows from a result of
Elek~\cite[Theorem~1]{ElekFrMiAcCoGrInPrMe},  that  there exists a
free and minimal action of $\Gamma$ on $X$ which admits an ergodic
(regular, non-atomic) invariant Borel
probability measure $\mu$.  Let $t\mapsto U_t$ be the unitary
representation of $\Gamma$ on $L^2:=L^2(X,\mu)$:
$(U_t\xi)(s)= \xi(t^{-1}s)$ and $M$ be the representation of $C(X)$ by
multiplication operators on $L^2$.

Define $\mathcal C$ to be the $C^*$-algebra generated by the images
of $U$ and $M$, and put $\mathcal D=M(C(X))$.
(It turns out $\mathcal C$ is an exotic crossed product.)
The freeness of the action implies $(\mathcal C,\mathcal D)$
is a regular MASA inclusion with the unique state extension property,
so there exists a (unique) conditional
expectation $E: \mathcal C\rightarrow \mathcal D$.

Let $J:=\{x\in \mathcal C: E(x^*x)=0\}$ be the left kernel of $E$.
By~\cite[Theorem~3.15]{PittsStReInI}, $J$ is a ideal of $\mathcal C$
having trivial intersection with $\mathcal D$ and
\cite[Theorem~3.7]{ExelPittsZarikianExIdFrTrGrC*Al} shows
$J\neq\{0\}$.  By~\cite[Lemma~3.1]{ArchboldBunceGregsonExStC*AlII},
$(\mathcal C/J, \mathcal D)$ is a regular inclusion having the unique
state extension property.  It follows $(\mathcal C/J, \mathcal D)$ is
a regular inclusion with the extension property which has a faithful
conditional expectation; thus $(\mathcal C/J, \mathcal D)$ is a
$C^*$-diagonal.  We do not know whether $\C$ contains a subalgebra
$\D'$ such that $(\C,\D')$ is a pseudo-Cartan inclusion (or Cartan
inclusion).

\begin{remark}{Question}\label{quotCar}   Suppose $\C$ is a
\cstaralg, $J\idealin\C$ and $\C/J$ admits a subalgebra $\E$ such
that  $(\C/J, \E)$ is a \pd\ (or
  a Cartan inclusion).  When can one conclude that  exists $\D\subseteq \C$ such that
  $(\C,\D)$ is a \pd\ (or Cartan inclusion)?
\end{remark}

Next, the uniqueness statement of Theorem~\ref{!pschar} suggests the
possibility of finding a canonical method of replacing certain 
non-Hausdorff twists with  Hausdorff ones.

\begin{remark}{Question}\label{wCar+pd}   Suppose $(\Sigma, G)$ is a twist
  over the second countable, \' etale and topologically free groupoid
  $G$ (see~\cite[Definition~3.4.6]{ExelPittsChGrC*AlNoHaEtGr} for the
  definition of a topologically free groupoid).  We assume $\unit{G}$
  is Hausdorff, but we do not assume $G$ is Hausdorff. Then
  $(C^*_{ess}(\Sigma,G), C_0(\unit{G}))$ is a weak-Cartan inclusion
  (see~\cite[Corollary~3.9.5]{ExelPittsChGrC*AlNoHaEtGr}) 
  and hence a \pd.  Let $(\A,\B\ms \alpha)$ be its Cartan envelope,
  and let $(\Sigma_{\A} , G_{\A})$ be the twist associated to
  $(\A,\B)$.  Then $G_{\A}$ is Hausdorff.  What is the relationship
  between $(\Sigma_{\A} , G_{\A})$ and $(\Sigma,G)$?
\end{remark}

We close with a technical question.  Several of the results in
Section~\ref{apps} assume that the unitization of a \pd\ is again a
\pd\ and it would be interesting to know if that hypothesis can be
removed.   The issue is whether $(\tilde\C,\tilde\D)$ is regular.
\begin{remark}{Question}\label{unpdq}   If $(\C,\D)$ is a \pd, must
  $(\tilde\C,\tilde\D)$ be a \pd?
\end{remark}  

\appendix

\mysec{}

\label{AppendixA} 

\begin{remark*}{Changes in Notation and Terminology}
  \it Throughout Appendix~\ref{AppendixA}, we depart from the notation
  and terminology used in Sections~\ref{1intro}--\ref{apps} above;
  instead we use the notation and terminology found
  in~\cite{PittsStReInII}.  In particular, all inclusions below are
  unital, and for a unital \cstaralg\ $B$, an essential extension for
  $\B$ will mean a pair $(\A,\alpha)$ consisting of a unital
  \cstaralg\ $\A$ and a unital $*$-monomorphism
  $\alpha: \B\rightarrow \A$ such that $\alpha(\B)\subseteq \A$ has
  the \iip.
\end{remark*}

Lemma~2.3 of~\cite{PittsStReInII} is
presented without proof, and there is an error in its statement.
The error 
is the assertion that the map $\ropen(\hat\A)\ni G\mapsto r^{-1}(G)$
is a Boolean algebra isomorphism of $\ropen(\hat\A)$ onto
$\ropen(\hat\B)$: it should have said
$\ropen(\hat\A)\ni G\mapsto \left(\overline{r^{-1}(G)}\right)^\circ\in
\ropen (\hat\B)$ is a Boolean algebra isomorphism.  Here is the full
and 
corrected statement.
\begin{lemma}[{Corrected \cite[Lemma~2.3]{PittsStReInII}}]\label{isolattice}
  Suppose $\A$ and $\B$ are abelian, unital
  \cstaralg s, $(\B,\alpha)$ is an essential
  extension of $\A$, and  $r:\hat\B\rightarrow \hat\A$ is the
  continuous surjection,
  $\rho\in\hat\B\mapsto \rho\circ \alpha$.  Then the maps
  \begin{align}\rideal(\B)\ni J\mapsto \alpha^{-1}(J) &\dstext{and}
    \ropen(\hat\A) \ni G\mapsto
    \left(\overline{r^{-1}(G)}\right)^\circ\label{A1.1}\\
\intertext{are Boolean algebra
  isomorphisms of $\rideal(\B)$ onto $\rideal(\A)$ and
  $\ropen(\hat\A)$ onto $\ropen(\hat\B)$ respectively.  The inverses
  of these maps are}
  \rideal(\A)\ni K\mapsto \alpha(K)^\dperp&\dstext{and}
    \ropen(\hat\B)\ni H\mapsto (r(\overline{H}))^\circ\label{A1.2}\\
\intertext{respectively.  
Furthermore, for $J\in\rideal(\B)$ and $G\in \ropen(\hat\A)$,} 
  \supp(\alpha^{-1}(J))=
    (r(\overline{\supp(J)}))^\circ
    &\dstext{and}\ideal(\left(\overline{r^{-1}(G)}\right)^\circ)=\alpha(\ideal(G))^\dperp.\label{A1.3}
  \end{align}
\end{lemma}

A complete proof of Lemma~\ref{isolattice} may be found in
\cite{PittsIrMaIsBoReOpSeReId}.   The proofs
of~\eqref{A1.1} and \eqref{A1.2} are found in
\cite[Propositions~3.2 and~4.17]{PittsIrMaIsBoReOpSeReId} and
\eqref{A1.3} is \cite[Lemma~4.12]{PittsIrMaIsBoReOpSeReId} combined
with~\cite[Lemma~4.7]{PittsIrMaIsBoReOpSeReId}.    
To make the notational
transition from \cite{PittsIrMaIsBoReOpSeReId} to the notation of
Lemma~\ref{isolattice}, take $Y:=\hat\B$, $X=\hat{A}$ and $\pi:=r$.

Unfortunately, the misstatement in~\cite[Lemma~2.3]{PittsStReInII}
leads to a gap in the proof of \cite[Theorem~6.9]{PittsStReInII}.
The statement
of~\cite[Theorem~6.9]{PittsStReInII} is correct, but its proof is
insufficient, as we now describe.

Let $(\C,\D)$ be a regular inclusion with the unique pseudo-expectation
property, let $(I(\D),\iota)$ be an injective envelope for $\D$ and let
$E:\C\rightarrow I(\D)$ be the pseudo-expectation. In the
proof of~\cite[Theorem~6.9]{PittsStReInII} we  claimed that if
$v\in \N(\C,\D)$ and $\rho\in\widehat{I(\D)}$ satisfies
$\rho(E(v))\neq 0$, then $\rho\circ\iota\in (\fix\beta_v)^\circ$.
That claim is then used show that $\rho\circ E$ is a compatible
state.  Our proof of this claim is insufficient because it uses the
incorrectly stated portion of \cite[Lemma~2.3]{PittsStReInII}.

We now give a proof of \cite[Theorem~6.9]{PittsStReInII}, whose
statement we have reproduced in Theorem~\ref{upse=>cov}.  The proof
uses the correction to~\cite[Lemma~2.3]{PittsStReInII} given as
Lemma~\ref{isolattice} above  to show that
$\rho\circ E$ is a compatible state.

\begin{theorem}[{\cite[Theorem~6.9]{PittsStReInII}}]\label{upse=>cov} Suppose $(\C,\D)$ is a regular
  inclusion with the unique pseudo-expectation property.  Then
  $(\C,\D)$ is a covering inclusion and $\fS_s(\C,\D)$ is a compatible
  cover for $\hat\D$.  Furthermore, if $F$ is any closed subset of
  $\Mod(\C,\D)$ which covers $\hat\D$, then $\fS_s(\C,\D)\subseteq F$.
  \end{theorem}
  \begin{proof}

 Denote by $r$ the ``restriction'' map, $\widehat{I(\D)} \ni
\rho\mapsto \rho\circ\iota\in \hat\D$.  Since $(I(\D),\iota)$ is an
injective envelope for $\D$, it is in particular an essential
extension of $\D$.  

 For $v\in \N(\C,\D)$,  let $$U_v:=\{\rho\in
\widehat{I(\D)}: \rho(E(v))\neq 0\}.$$   If $\rho\in U_v$, the
Cauchy-Schwartz inequality gives,
\[|\rho(E(v))|^2\leq \rho(E(v^*v))=(\rho\circ\iota)(v^*v),\] so
$r(\rho)\in \dom\beta_v$.   By~\cite[Lemma~2.5]{PittsStReInI},
$\rho\circ E|_\D=r(\rho)\in \fix\beta_v$.  Thus,
\begin{equation}\label{fix}r(U_v)\subseteq
    \fix(\beta_v).
  \end{equation}
  
Let us  show that every $\tau\in\fS_s(\C,\D)$ is a compatible
state. 
 Fix $\tau\in \fS_s(\C,\D)$ and suppose $\tau(v)\neq 0$ for some $v\in
 \N(\C,\D)$.  Our task is to show that $|\tau(v)|^2=\tau(v^*v)$.   Write $\tau=\rho\circ E$ for some $\rho\in
 \widehat{I(\D)}$.  
Let
 \[X:=\{\rho'\in \widehat{I(\D)}: \rho'(E(v))> 
   \rho(E(v))/2\}.\]  By construction, $X$ is an open subset of
 $\widehat{I(\D)}$ and $\rho\in X$.  Since 
 $\overline X\subseteq U_v$, ~\eqref{fix} gives,
 \[r(\overline X)^\circ\subseteq (\fix\beta_v)^\circ.\]
 Put $G:=r(\overline X)^\circ$.  Since $\widehat{I(\D)}$ is a Stonean
   space, $\overline X$ and $\overline{r^{-1}(G)}$ are clopen sets and
   hence are regular open subsets of $\widehat{I(\D)}$.  
 Then
 \[\overline{r^{-1}(G)}=(\overline{r^{-1}(G)})^\circ=\left(\overline{r^{-1}(
     r(\overline
X)^\circ)}\right)^\circ =\overline X,\]
with the last equality following from the portions of
Lemma~\ref{isolattice} concerning regular open sets.   Therefore, we may
 find a net $(\rho_\lambda)$ in $r^{-1}(G)$ such that
 $\rho_\lambda\rightarrow \rho$.   Since
 $r(\rho_\lambda)=\rho_\lambda\circ\iota\in(\fix\beta_v)^\circ$, we
 may find $d_\lambda\in \ideal((\fix\beta_v)^\circ)$ with
 $\rho_\lambda(\iota(d_\lambda))=1$.   Since $v d_\lambda\in \D^c$
 (see~\cite[Lemmas~2.14 and~2.15]{PittsStReInII}),
 and 
 $E|_{\D^c}$ is a homomorphism (by~\cite[Lemma~6.8]{PittsStReInII})
 \[|\rho_\lambda(E(v))|^2=|\rho_\lambda(E(v)\iota(d_\lambda))|^2=|\rho_\lambda(E(vd_\lambda))|^2
   =|\rho_\lambda(E(d_\lambda^*v^*vd_\lambda))|
   =\rho_\lambda(E(v^*v)).\]  As $\tau=\rho\circ E$,
 \[|\tau(v)|^2=\lim_\lambda
   |\rho_\lambda(E(v))|^2=\lim_\lambda\rho_\lambda(E(v^*v)) =\tau(v^*v).\]
It follows that $\tau$ is a compatible state.

The remainder of the proof now follows exactly as in the proof
of~\cite[Theorem~6.9]{PittsStReInII};  we include the details for
convenience.

Let us show the invariance of $\fS_s(\C,\D)$.  Choose $\tau\in
\fS_s(\C,\D)$ and write $\tau=\rho\circ
E$ for some $\rho\in \widehat{I(\D)}$.  Suppose $v\in \N(\C,\D)$ is
such that $\rho(\iota(v^*v))\neq 0$.  Let $P$ and $Q$ be the support
projections in $I(\D)$ for the ideals $\overline{vv^*\D}$ and $\overline{v^*v\D}$
respectively.    Then $\tilde\theta_v$ is a partial automorphism with
domain $PI(\D)$ and range $QI(\D)$.   Define $\tau'\in
\widehat{I(\D)}$ by \[\tau'(h)=\rho(\tilde\theta_v(Ph)), \qquad h\in
  I(\D).\]
  For $x\in \C$,
\cite[Proposition~6.2]{PittsStReInII} gives,
\begin{align*}
  \rho(E(v^*xv))&= \rho(\tilde\theta_v(E(vv^*x)))=
                  \rho(\tilde\theta_v(\iota(vv^*) PE(x)))\\
  &=\rho(\iota(v^*v))\rho(\tilde\theta_v(PE(x)))=
    \rho(\iota(v^*v))(\tau'\circ E)(x).
\end{align*}
Thus $\fS_s(\C,\D)$ is invariant.

If $\sigma\in\hat\D$, choose any $\rho\in \widehat{I(\D)}$ such that
$\rho\circ\iota=\sigma$.  Then $\sigma=(\rho\circ E)|_\D$, so
$\fS_s(\C,\D)$ covers $\hat\D$.  Thus, $\fS_s(\C,\D)$ is a compatible
cover for $\hat\D$ and $(\C,\D)$ is a covering inclusion.

Finally, if $F\subseteq \Mod(\C,\D)$ is closed and covers $\hat\D$, then $\fS_s(\C,\D)\subseteq F $ by \cite[Theorem~6.1(b)]{PittsStReInII}.

\end{proof}


\def\cprime{$'$}

\end{document}